\documentclass[12pt]{amsart}

\newtheorem{theorem}{Theorem}[section] 
\newtheorem{corollary}[theorem]{Corollary} 
\newtheorem{conjecture}[theorem]{Conjecture} 
\newtheorem{lemma}[theorem]{Lemma} 
\newtheorem{proposition}[theorem]{Proposition}

\theoremstyle{definition}
\newtheorem{definition}{Definition}

\theoremstyle{remark}

  {\end{list}}

\usepackage{amssymb} 
\usepackage{verbatim}

\def\deg{{\mathop{\rm deg}}}

\def\CPA{{\mathop{\rm CPA}}}
\def\CR{{\mathop{\rm CR}}}
\def\Zh{{\mathop{\rm Zh}}}
\def\Meas{{\mathop{\rm Meas}}}
\def\Re{{\mathop{\rm Re}}}
\def\Im{{\mathop{\rm Im}}}
\def\VEC{{\mathop{\rm Vec}}}
\def\IN{{\mathop{\rm In}}}
\def\OUT{{\mathop{\rm Out}}}

\def\BDV{{\mathop{\rm BDV}}}
\def\Ker{{\mathop{\rm Ker}}}
\def\Dir{{\mathop{\rm Dir}}}
\def\can{{\mathop{\rm can}}}
\def\Tr{{\mathop{\rm Tr}}}

\def\<{\langle }
\def\>{\rangle }

\def\({(\!(}
\def\){)\!)}

\def\BF1{{\mathbf 1}}

\def\ff{{\mathbf f}}

\def\CC{{\mathbb C}}

\def\RR{{\mathbb R}}
\def\ZZ{{\mathbb Z}}

\def\cA{{\mathcal A}}

\def\cC{{\mathcal C}}

\def\cH{{\mathcal H}}

\def\cQ{{\mathcal Q}}
\def\cD{{\mathcal D}}

\def\cF{{\mathcal F}}

\def\vv{\vec{v}}

\def\span{{\mathop{\rm span}}}

\newcommand{\la}{\langle}               
\newcommand{\ra}{\rangle}               

\DeclareMathSymbol{\varnothing} {\mathord}{AMSb}{"3F} 
 
\theoremstyle{definition} 
\theoremstyle{remark} 
\numberwithin{equation}{section}

\begin{document}
\title[Harmonic Analysis]
{Harmonic Analysis on Metrized Graphs}


\author{Matt Baker}
\address{Matt Baker\\
School of Mathematics\\
Georgia Institute of Technology\\
Atlanta, Georgia 30332-0160\\
USA}
\email{mbaker@math.gatech.edu}

\author{Robert Rumely}
\address{Robert Rumely\\ 
Department of Mathematics\\
University of Georgia\\
Athens, Georgia 30602\\
USA}
\email{rr@math.uga.edu}

\thanks{The authors' research was supported in part by NSF grant DMS
  0070736.  The authors would like to thank the participants of the Summer 2003 REU program ``Analysis on
  Metrized Graphs'' for their enthusiasm, and for doing a number of useful computations.
  They would also like to thank Russell Lyons, Xander Faber, Zubeyir Cinkir, and the anonymous referee for helpful feedback 
  on earlier drafts of this paper.}

\begin{abstract} 
This paper studies the Laplacian operator on a metrized graph, and its
spectral theory.
\end{abstract} 

\maketitle

\section{Introduction}

\vskip .1 in

A metrized graph $\Gamma$ is a finite connected graph 
equipped with a distinguished parametrization of each of its edges.  
(A more rigorous definition can be found in \S\ref{MetrizedGraphSection}).  
In particular, $\Gamma$ is a one-dimensional manifold except at
finitely many ``branch points'', where it looks locally like an
$n$-pointed star.  A metrized graph should be thought 
of as an analytic object,  not just a combinatorial one. 
In particular there are a Laplacian operator
$\Delta$ on $\Gamma$ and corresponding ``Green's functions'' used to
invert it.  In this paper, we will study the spectral theory of the
Laplacian on a metrized graph.  

\vskip .1 in
While metrized graphs occur in fields as diverse as chemistry,
physics, and mathematical biology, 
our motivation comes from arithmetic geometry.  
The eigenfunction decomposition of the Arakelov-Green's function
on a Riemann surface 
\begin{equation} \label{SpecDec}
g_{\mu}(x,y) \ = \ \sum_{n=1}^{\infty} 
         \frac{1}{\lambda_n} f_n(x) \overline{f_n(y)} 
\end{equation}
plays an important role in Arakelov theory.  For example, it is the key 
ingredient used by Faltings and Elkies to give lower bounds for 
$g_{\mu}(x,y)$-discriminant sums  (see \cite[\S VI, Theorem 5.1]{Lang}).  
Estimates for these sums can then be used to prove the
effectivity of certain Arakelov divisors on an arithmetic surface
(see \cite[\S V, Theorem 4.1]{Lang}).  

Metrized graphs were introduced
as a nonarchimedean analogue of a Riemann surface
in (\cite{RumelyBook}, \cite{CR}, and \cite{Zh1}).  
The Laplacian on a metrized graph 
and the kernel of integration that inverts it,
the Arakelov-Green's function $g_\mu(x,y)$ associated to a measure 
$\mu$ on $\Gamma$, were studied in \cite{CR} and \cite{Zh1}.   
Those papers made it natural to ask if there is an 
eigenfunction decomposition of $g_{\mu}(x,y)$ 
analogous to (\ref{SpecDec}) on a metrized graph, 
keeping in mind that for applications to nonarchimedean Arakelov theory,
it is necessary to allow $\mu$ to be a singular measure.  
This paper therefore addresses the question 
``What are the eigenfunctions of the Laplacian 
on a metrized graph $\Gamma$ relative to an arbitrary bounded,
signed measure $\mu$?''  

\vskip .1 in
Initially, we discovered decompositions of $g_{\mu}(x,y)$ 
for certain special examples by ad hoc means (see \S\ref{BriefExamplesSection} below).  
This paper is the outcome of our attempt to understand those examples.   
For an arbitrary measure $\mu$ of total mass one on $\Gamma$, we obtain 
an eigenfunction expansion converging uniformly to $g_{\mu}(x,y)$.  
We give three characterizations of eigenfunctions and study their properties.  
We give an algorithm which computes the eigenvalues 
in certain cases, and report some examples for interesting measures and graphs 
found by students in an REU held at the University of Georgia in
summer 2003.  As applications of the general theory, 
we prove the positivity of an ``energy pairing'' on 
the space of measures of total mass $0$,
we obtain Elkies-type bounds for $g_{\mu}(x,y)$-discriminant sums generalizing
those found for the circle graph in \cite{HS}, and we
show that $\mu = dx$ is the unique measure of total mass one for which the integral operator
associated to $g_{\mu}(x,y)$ has minimal trace.  We also raise several questions
for future investigation.

\vskip .1 in
In the course of writing this paper, 
and from the referee's report, we learned of an extensive 
mathematical literature concerning spectral analysis and 
differential equations on metrized graphs, 
which are also known in the literature as ``networks'', ``metric graphs'', 
and ``quantum graphs''.  
In particular we note the papers of J.~P.~Roth \cite{Roth}, 
C.~Cattaneo \cite{Cattaneo}, S.~Nicaise \cite{Nicaise1}, \cite{Nicaise2}, 
J.~von Below \cite{Be1}, \cite{Be2}, F.~A.~Mehmeti \cite{AM1}--\cite{AM3}, 
B.~Gaveau, M.~Okada, and T.~Okada \cite{GOO}, T.~Okada \cite{Okada},
M.~Solomyak \cite{Solomyak}, and L.~Friedlander~\cite{Friedlander},
\cite{Friedlander2}.
These works all define the Laplacian somewhat differently than we have
done in the present paper.  
Typically, they deal with Laplacians of functions $f$ satisfying 
certain natural boundary conditions at the nodes (vertices) 
which make the ``discrete part'' of the Laplacian (as defined in the present
paper) equal to zero.  They then study (in our terminology) 
the spectral theory of the graph relative to the measure $\mu = dx$.
The results in the papers cited above, while perhaps not entirely 
superseding ours in the case $\mu = dx$, go much further in a
number of directions.  In particular, most of these papers treat 
infinite networks, whereas we have restricted our attention to finite ones.
So to a large extent, it is the emphasis on eigenfunctions with respect 
to general measures $\mu$ (of total mass one) which is new in this paper.
We have included Section~\ref{CanonicalTauSection} below partly to 
illustrate why one might be interested in the theory for
measures other than $dx$.  The spectral theory of more general measures 
is at least partly treated in \cite{AM2}, but we are unable to draw 
a precise comparison because the terminology and definitions are so different.

For more discussion of the relationship between these works and the present
paper, see \S\ref{RelatedWorkSection} below.

%

\vskip .1 in

We now give an overview of the contents of this paper.
First, we will recall the notions of Laplacian operators on $\Gamma$ given in
\cite{CR} and \cite{Zh1}.

\vskip .1 in

\subsection{The Laplacian on the space $\CPA(\Gamma)$}

Let $\CPA(\Gamma)$ be the space of continuous, 
piecewise affine complex-valued
functions on $\Gamma$.  

For each point $p \in \Gamma$, consider the set
$\VEC(p)$ of ``formal unit vectors $\vv$ emanating from $p$''. 
(A rigorous definition is given in \S\ref{MetrizedGraphSection}).
The cardinality of $\VEC(p)$ is just the ``valence'' of $p$ (which is
equal to 2 for all but finitely many points $p \in \Gamma$).
If $\vv \in \VEC(p)$, we define the directional derivative $d_{\vv} f(p)$ of
$f$ at $p$ (in the direction $\vv$) to be
$$
d_{\vv} f(p) = \lim_{t \to 0^+} \frac{f(p + t\vv) - f(p)}{t}.
$$

For $f \in \CPA(\Gamma)$, we define the Laplacian $\Delta f$ of $f$ to
be the discrete measure
$\Delta f = -\sum_{p \in \Gamma} \sigma_p(f) \delta_p$,
where $\sigma_p(f)$ is the ``sum of the slopes of $f$ in all
directions pointing away from $p$''; more precisely, 
$\sigma_p(f) = \sum_{\vv \in \VEC(p)} d_{\vv} f(p)$.
(Here $\delta_p$ is the discrete probability measure on 
$\Gamma$ with the property that $\int_\Gamma f(x) \delta_p(x) = f(p)$
for all $f \in \cC(\Gamma)$.)

\vskip .1 in

This is essentially the definition of the Laplacian $\Delta_{\CR}$ given 
by Chinburg and Rumely in \cite{CR}, 
except that our Laplacian is the negative of theirs.  
We also remark that the definition of $\Delta$ on $\CPA(\Gamma)$ 
is closely related to the
classical definition of the Laplacian matrix for a finite
weighted graph (see \cite{BF} for details). 


\vskip .1 in

\subsection{The Laplacian on the Zhang space}

Following Zhang \cite{Zh1}, one can define the Laplacian on a more
general class of functions than $\CPA(\Gamma)$.  The resulting
operator combines aspects of the ``discrete'' Laplacian from the previous
section and the ``continuous'' Laplacian $-f''(x)dx$ on $\RR$.  

We define the `Zhang space' $\Zh(\Gamma)$ 
to be the set of all continuous functions
$f : \Gamma \rightarrow \CC$ such that $f$ is piecewise $\cC^2$ and
$f^{\prime \prime}(x) \in L^1(\Gamma)$.
(In general, when we say that $f$ is piecewise $\cC^k$, we mean
that there is a finite set of points $X_f \subset \Gamma$
such that $\Gamma \backslash X_f$ is a finite union of open intervals,
and the restriction of $f$ to each of those intervals is $\cC^k$.)

For $f \in \Zh(\Gamma)$, we define $\Delta f$ to be the complex
Borel measure on $\Gamma$ given by
$$
\Delta(f) = -f''(x)dx - \sum_{p \in X_f} \left( \sum_{\vv \in
    \VEC(p)} d_{\vv} f(p) \right) \delta_p,
$$
where if $x = p + t \vv \in \Gamma \backslash X_f$, we set
$f''(x) = \frac{d^2}{dt^2} f(p+t \vv)$.

If $f \in \CPA(\Gamma) \subset \Zh(\Gamma)$, then clearly the two
definitions of the Laplacian which we have given agree with each other.

\medskip

Let $\la \; , \; \ra$ denote the $L^2$ inner product on $\Zh(\Gamma)$,
i.e., $\la f,g \ra = \int_\Gamma f(x) \overline{g(x)} dx$.
The following result shows that the Laplacian on $\Zh(\Gamma)$
is ``self-adjoint'', and justifies the choice of sign in the definition
of $\Delta$.  

\begin{proposition}
\label{SelfAdjointProp}
For every $f,g \in \Zh(\Gamma)$, $\int_\Gamma \overline{g} \, d\Delta f =
\int_\Gamma f \, d\overline{\Delta g} = \la f', g' \ra$.
%
\end{proposition}

\begin{proof}  This formula is given in Zhang (\cite{Zh1}, proof of
Lemma~a.4), and is a simple consequence of integration by parts.
We repeat the proof for the convenience of the reader.

Since $f, g \in \Zh(\Gamma)$, there is a vertex set $X$ for $\Gamma$ such that
$\Gamma \backslash X$ consists of a finite union of open intervals,
$f^{\prime \prime}(x)$ and $g^{\prime \prime}(x)$ are continuous on 
$\Gamma \backslash X$, and
$f^{\prime \prime}(x), g^{\prime \prime}(x) \in L^1(\Gamma)$.  
In particular $f^{\prime}$ and $g^{\prime}$
are continuous on $\Gamma \backslash X$.  Using that $f''(x)$ and $g''(x)$
belong to $L^1(\Gamma)$ one sees that  
$f^{\prime}$ and $g^{\prime}$ have boundary limits
on each of the open segments making up $\Gamma \backslash X$, and these limits
coincide with the directional derivatives at those points.  

Identifying $\Gamma \backslash X$ with a union of open intervals $(a_i,b_i)$,
$i = 1, \ldots, N$, and identifying the directional derivatives at the points
in $X$ with one-sided derivatives at the endpoints, we find by integration by
parts that 
\begin{equation*}
\begin{aligned}
\<f,g\>_{\Dir} 
   & \ = \ \sum_{i=1}^N 
     \big(-\int_{a_i}^{b_i} f(x) \overline{g^{\prime \prime}(x)} \, dx 
          + f(b_i) \overline{g^{\prime}(b_i)}
          - f(a_i) \overline{g^{\prime}(a_i)} \big) \\
   & \ = \ \sum_{i=1}^N \int_{a_i}^{b_i} 
               f^{\prime}(x) \overline{g^{\prime}(x)} \, dx 
   \ = \ \int_{\Gamma} f^{\prime}(x) \overline{g^{\prime}(x)} \, dx.
\end{aligned}
\end{equation*}   
\end{proof}


It follows that the Laplacian has total mass 0:

\begin{corollary}
For any $f \in \Zh(\Gamma)$, we have $\int_\Gamma \Delta f = 0$.
\end{corollary}

\begin{proof}
$\int_\Gamma \Delta f = 
\int_\Gamma 1 \, d\Delta f = 
\int_\Gamma f \, d\overline{\Delta 1} = 0$.
\end{proof}

This generalizes the fact (see \cite[Corollary~2.9]{CR}) that if $f \in \CPA(\Gamma)$
and $\Delta f = \sum_i c_i \delta_{q_i}$, then $\sum c_i = 0$.

\vskip .1 in

\subsection{The Dirichlet inner product and eigenfunctions of the Laplacian}
\label{PreliminaryDirichletSection}

Motivated by Proposition~\ref{SelfAdjointProp},
we define the {\it Dirichlet inner product} on
$\Zh(\Gamma)$ by
$$
\la f, g \ra_{\Dir} = \int_\Gamma f'(x) \overline{g'(x)} dx.
$$

We define a corresponding {\it Dirichlet seminorm} by
$\| f \|_{\Dir}^2 = \la f,f \ra_{\Dir}$.

Note that
$\| f \|_{\Dir} \geq 0$ for all $f \in \Zh(\Gamma)$, with equality if and only if $f$
is constant.  Thus $\| \; \|_{\Dir}$ is indeed a seminorm on $\Zh(\Gamma)$.

To obtain an honest {\it norm}, 
let $\mu$ be a real-valued signed Borel measure on $\Gamma$ with
$\mu(\Gamma) = 1$ and $|\mu|(\Gamma) < \infty$, and define
$\Zh_\mu(\Gamma) = \{ f \in \Zh(\Gamma) \; : \; \int_\Gamma f d\mu = 0 \}$.
Then for $f \in \Zh_\mu(\Gamma)$,
$\|f\|_{\Dir} = 0$ if and only if $f = 0$ (i.e., $\| \, \|_{\Dir}$ defines a norm on $\Zh_\mu(\Gamma)$).

\medskip

\medskip

%
%


Let $\Dir_\mu(\Gamma)$ denote the Hilbert space completion of $\Zh_\mu(\Gamma)$ with
respect to the Dirichlet norm.  
By Corollary~\ref{Cor26} below, 
$\Dir_\mu(\Gamma)$ can be naturally identified with a subspace of
$\cC(\Gamma)$, the space of continuous complex-valued functions on
$\Gamma$.

We define a nonzero function $f \in \Zh_\mu(\Gamma)$ to be an 
{\it eigenfunction of the Laplacian $($with respect to $\mu$$)$} if
there exists $\lambda \in \CC$ such that
$$
\int_\Gamma \overline{g(x)} \; d(\Delta f) = \lambda \int
\overline{g(x)} f(x) \, dx,
$$
or equivalently $\la f, g \ra_{\Dir} = \lambda \la f, g \ra$,
for all $g \in \Dir_\mu(\Gamma)$.

It is important to note that we only test this identity for $g$
in the space $\Dir_\mu(\Gamma)$, so in particular the eigenvalues and
eigenfunctions depend on the choice of $\mu$.  

The condition $\la f, g \ra_{\Dir} = \lambda \la f, g \ra$ makes sense
for $f,g \in \Dir_\mu(\Gamma)$ as well, so there is a natural
extension to $\Dir_\mu(\Gamma)$ of the notion of an eigenfunction of
the Laplacian.  (It is irrelevant whether we require $g$ to be in
$\Zh_\mu(\Gamma)$ or $\Dir_\mu(\Gamma)$.)  Later we will see that the space 
$\Zh_\mu(\Gamma)$ is not always large enough to contain a complete 
orthonormal set of eigenfunctions, and we will describe a space which does.

Setting $f=g$ in the definition of an eigenfunction, it is easy to see
that the eigenvalues of the Laplacian on $\Zh_\mu(\Gamma)$ (or $\Dir_\mu(\Gamma)$)
are real and positive.

\medskip

{\bf Example 1:} Suppose $\mu = \frac{1}{2} \delta_0 + \frac{1}{2}
\delta_1$ on $\Gamma = [0,1]$.  We will see later (Proposition~\ref{Prop33}) 
that every eigenfunction of
the Laplacian on $\Dir_\mu(\Gamma)$ is piecewise $\cC^\infty$.
It follows that $f$ is an eigenfunction
of $\Delta$ with respect to $\mu$ iff there exists $\lambda \in \RR$ such that
$$
\la f, g \ra_{\Dir} = f'(1)\overline{g(1)} - f'(0)\overline{g(0)} - 
\int_0^1 f''(x)\overline{g(x)} \, dx
= \lambda \int_0^1 f(x) \overline{g(x)} \, dx
$$
for all $g \in \Zh_\mu(\Gamma)$.

As $g(0) = -g(1)$ for all $g \in \Zh_\mu(\Gamma)$, this simplifies to
$$
\overline{g(1)}[f'(1) + f'(0)] 
- \int_0^1 f''(x)\overline{g(x)}
\, dx = \lambda \int_0^1 f(x) \overline{g(x)} \, dx.
$$
It is not hard to show that this identity holds for all $g \in
\Zh_\mu(\Gamma)$ if and only if $f'(0) = -f'(1)$ and 
$f''(x) = -\lambda f(x)$.  As $f \in \Zh_\mu(\Gamma)$, we must also have
$f(0) = -f(1)$.  These conditions yield the eigenfunctions
$a_n \cos(n \pi x), b_n \sin(n \pi x)$ for all {\it odd} integers
$n\geq 1$, and corresponding eigenvalues $n^2\pi^2$ (each occuring
with multiplicity 2), again for $n\geq 1$ odd.  


\vskip .1 in

{\bf Example 2:} Suppose $\mu = dx$ is Lebesgue measure on $\Gamma = [0,1]$. 
It is shown in (\S\ref{ExampleSection}, Example 1) that eigenfunctions of the
Laplacian in $\Dir_{\mu}(\Gamma)$ satisfy $f''(x) = -\lambda f(x)$ with the
boundary conditions $f'(0) = f'(1) = 0$.  These conditions yield the 
eigenfunctions $a_n \cos(n \pi x)$ for $n = 1, 2, 3, ...$,  with eigenvalues
$\lambda_n = n^2 \pi^2$. 

\vskip .1 in
Comparing Examples (1) and (2) 
illustrates the fact that the eigenfunctions and eigenvalues 
depend not only on $\Gamma$, but on the choice of $\mu$.
\vskip .1 in

In \S\ref{EigenfunctionViaPairingsSection}, 
we will discuss how the classical Rayleigh-Ritz method can be used 
to obtain the following result:

\begin{theorem}
\label{EigenvalueTheorem1}
Each eigenvalue $\lambda$ of $\Delta$ on $\Dir_\mu(\Gamma)$ occurs with finite
multiplicity.  If we write the eigenvalues as $0 < \lambda_1 \leq
\lambda_2 \leq \cdots$ $($with corresponding eigenfunctions
$f_1,f_2,\ldots$, normalized so that $\| f_n \|_2 = 1$$)$, then
$\lim_{n\to\infty} \lambda_n = \infty$ and
$\{ f_n \}$ forms a basis for the Hilbert space $\Dir_\mu(\Gamma)$.
\end{theorem}

It follows that each $f \in \Dir_\mu(\Gamma)$ has a generalized Fourier expansion
\begin{equation}
\label{GeneralizedFourierExpansion}
f = \sum_{n=1}^\infty c_n f_n.
\end{equation}
We will see later that
the series in (\ref{GeneralizedFourierExpansion}) converges
uniformly to $f$.
However, the techniques of \S\ref{EigenfunctionViaPairingsSection}
are not strong enough to prove even pointwise convergence.

To establish stronger convergence results, we study 
eigenfunctions of the Laplacian via an integral
operator $\varphi_\mu$ which (in a suitable sense) is inverse to the Laplacian
on $\Dir_\mu$.  Working with integral operators, rather than differential
operators, is a well-known method for improving convergence.

\subsection{The general measure-valued Laplacian}

An important observation is that the Laplacian exists as a measure-valued operator on a larger
space than $\Zh(\Gamma)$.  In \S\ref{MeasureSection}, we will define a
space $\BDV(\Gamma) \subset \cC(\Gamma)$ which is in a precise sense the largest space of
continuous functions $f$ for which $\Delta f$ exists as a complex
Borel measure.  (The abbreviation BDV stands for ``bounded
differential variation''.  We refer to \S\ref{MeasureSection} for
the definition of $\BDV(\Gamma)$, and of $\Delta f$ for $f \in
\BDV(\Gamma)$.)  We also define 
$\BDV_\mu(\Gamma) = \{ f \in \BDV(\Gamma) \; : \; \int_\Gamma f \,
d\mu = 0 \}$.

Consideration of the measure-valued Laplacian on $\BDV(\Gamma)$ 
is important for our integral operator approach to the
spectral theory of $\Delta$, and the space $\BDV(\Gamma)$ plays a role
in our theory similar to that of a Sobolev space in classical analysis.

\subsection{The $j$-function and the Arakelov-Green's function $g_\mu$}
\label{jfunctionSection}

In \cite{CR} (see also \cite{BF}), a kernel $j_{\zeta}(x,y)$ giving a fundamental solution
of the Laplacian is defined and studied as a function of $x, y, \zeta
\in \Gamma$.
For fixed $\zeta$ and $y$ it has the following physical interpretation:
when $\Gamma$ is viewed as a resistive electric
circuit with terminals at $\zeta$ and $y$, with the resistance in each edge
given by its length, then $j_{\zeta}(x,y)$ is the voltage difference
between $x$ and $\zeta$, when unit current enters at $y$ and exits at
$\zeta$ (with reference voltage 0 at $\zeta$).
For fixed $y$ and $\zeta$, the function $j_\zeta(x,y)$ is in
$\CPA(\Gamma)$ as a function of $x$ and 
satisfies the differential equation
\begin{equation}
\Delta_x(j_{\zeta}(x,y)) \ = \ \delta_y(x) - \delta_{\zeta}(x) 
\end{equation}
(see \cite{CR}, formula (2), p.13;  recall that our $\Delta = - \Delta_{CR}$).


As before, 
let $\mu$ be a real-valued signed Borel measure of total mass $1$ on $\Gamma$.
We define the Arakelov-Green's function $g_\mu(x,y)$ associated to $\mu$ to be 
$$
g_\mu(x,y) = \int_\Gamma j_\zeta(x,y) \, d\mu(\zeta) - C_{\mu} 
$$
where $C_{\mu} = \int_\Gamma j_\zeta(x,y) \, d\mu(\zeta) d\mu(x) d\mu(y)$.  

One can characterize $g_\mu(x,y)$ as the unique function on $\Gamma \times \Gamma$
such that

\begin{itemize}
\item[$(1)$] 
$g_\mu(x,y)$ is an element of $\BDV_\mu(\Gamma)$ as a
  function of both $x$ and $y$.
\item[$(2)$] 
For fixed $y$, $g_\mu$ satisfies the identity
$$
\Delta_x g_\mu(x,y) = \delta_y(x) - \mu(x).
$$
\item[$(3)$] 
$\iint_{\Gamma \times \Gamma} g_\mu(x,y) d\mu(x) d\mu(y) = 0$.
\end{itemize}

\subsection{The integral transform $\varphi_\mu(x,y)$}

For $f \in L^2(\Gamma)$, define
$$
\varphi_\mu(f) = \int_\Gamma g_\mu(x,y) f(y) dy.
$$

We will see that $\varphi_\mu : L^2(\Gamma) \to \BDV_\mu(\Gamma)$ and
that 
$$
\Delta(\varphi_\mu(f)) = f(x)dx - (\int_\Gamma f(x) dx) \mu.
$$

In particular, if $\int_\Gamma g \, d\mu = 0$ then
$$
\int_\Gamma \overline{g(x)} \Delta(\varphi_\mu(f))  = \int_\Gamma
\overline{g(x)} f(x) dx.
$$

We will see that the operator $\varphi_\mu : L^2(\Gamma) \to
L^2(\Gamma)$ is compact and self-adjoint, and
then an application of the spectral theorem will yield:

\begin{theorem} 
\label{EigenvalueTheorem2}
The map $\varphi_{\mu}$ has countably many eigenvalues $\alpha_i$,
each of which is real and occurs with finite multiplicity.  
The nonzero eigenvalues can be ordered so that 
$|\alpha_1| \ge |\alpha_2| \ge \ldots$ and
$\lim_{i \rightarrow \infty} |\alpha_i| = 0$.
Each eigenfunction corresponding to a nonzero eigenvalue belongs
to $\BDV_{\mu}(\Gamma)$.   $L^2(\Gamma)$ has an orthonormal
basis consisting of eigenfunctions of $\varphi_{\mu}$, and each eigenfunction
corresponding to a nonzero eigenvalue belongs to $\BDV_{\mu}(\Gamma)$.  
\end{theorem}

We will also see that a nonzero function $f \in \Dir_\mu(\Gamma)$ is an eigenfunction for
$\Delta$ on $\Dir_\mu(\Gamma)$ with eigenvalue $\lambda$ if and only if
it is an eigenfunction of $\varphi_\mu$ with eigenvalue $\alpha = 1/\lambda$.

\vskip .1 in

This gives us an alternative way to
understand the spectral theory of the Laplacian.  In particular, it
will allow us to prove that the generalized Fourier series expansion
in (\ref{GeneralizedFourierExpansion}) converges uniformly on $\Gamma$ 
for every $f \in \Dir_\mu(\Gamma)$.

\subsection{Eigenfunction expansion of $g_\mu(x,y)$}

The integral operator approach will also enable us to 
prove the following result concerning the eigenfunction
expansion for $g_\mu(x,y)$:

\begin{proposition} 
\label{Prop18a} 
Let the $f_n$ be as in the statement of Theorem~\ref{EigenvalueTheorem1}.
Then the series $\sum_{n=1}^{\infty} \frac{f_n(x) \overline{f_n(y)}}{\lambda_n}$ 
converges uniformly to $g_{\mu}(x,y)$ for all $x, y \in \Gamma$.  
\end{proposition} 

Proposition~\ref{Prop18a} can be used to prove the following results, 
which are closely related to results in \cite{BRequidistr}
and served as part of the motivation for this paper.  For the
statement of Theorem~\ref{Energy1}, let $\Meas(\Gamma)$ denote the
space of all bounded signed Borel measures on $\Gamma$.

\begin{theorem} 
\label{Energy1}
Among all  $\nu \in \Meas(\Gamma)$ with $\nu(\Gamma) = 1$,  
$\mu$ is the unique measure minimizing the energy integral
\begin{equation}
I_{\mu}(\nu) \ = \  \iint_{\Gamma \times \Gamma} g_{\mu}(x,y) \, 
                   d\nu(x) \overline{d\nu(y)} \ .
\end{equation}
\end{theorem}

\begin{proposition}
\label{Energy2}
There exists a constant $C > 0$ such that if $N\geq 2$ and
$x_1,\ldots,x_N \in \Gamma$, then
\begin{equation*}
\frac{1}{N(N-1)} \sum_{i \neq j} g_\mu(x_i,x_j) \geq -\frac{C}{N}.
\end{equation*}
\end{proposition}

We remark that Theorem~\ref{Energy1} is motivated by classical results from potential
theory, and Proposition~\ref{Energy2} specializes to
\cite[Lemma~2.1]{HS} when $\Gamma$ is a circle.

\vskip .1 in

In the latter part of the paper, we will study the regularity and boundedness
of the eigenfunctions, and give an algorithm to compute them for 
a useful class of measures $\mu$.  

\vskip .1 in

\subsection{Examples of eigenfunction expansions}
\label{BriefExamplesSection}

\vskip .1 in

Here are some examples of eigenfunction expansions for $g_\mu(x,y)$.

\vskip .1 in

\begin{itemize}
\item[1.] $\Gamma = [0,1]$, $\mu = dx$.  
$$
\begin{array}{lll}
g_\mu(x,y) & = & \left\{ 
\begin{array}{ll}
\frac{1}{2} x^2 + \frac{1}{2}(1-y)^2 - \frac{1}{6} & {\rm if \;} x<y
\\
\frac{1}{2} (1-x)^2 + \frac{1}{2}y^2 - \frac{1}{6} & {\rm if \;} x\geq
y \\
\end{array}
\right.
\\
           & = & 2 \sum_{n\geq 1} \frac{\cos(n \pi x)\cos(n \pi
           y)}{\pi^2 n^2}. \\
\end{array}
$$
\item[2.] $\Gamma = [0,1]$, $\mu = \frac{1}{2} \delta_0
  + \frac{1}{2}\delta_1$.  
$$
\begin{array}{lll}
g_\mu(x,y) & = & \frac{1}{4} - \frac{1}{2} |x-y| \\
           & = & \sum_{n\geq 1 \; {\rm odd}} \left( \frac{\cos(n \pi
           x)\cos(n \pi y)}{\pi^2 n^2} + \frac{\sin(n \pi x)\sin(n \pi
           y)}{\pi^2 n^2} \right). \\
\end{array}
$$
\item[3.] $\Gamma = [0,1]$, $\mu = \delta_0$.  
$$
\begin{array}{lll}
g_\mu(x,y) & = & \min \{ x,y \} \\
           & = & 8 \sum_{n\geq 1 \; {\rm odd}} \frac{\sin(\frac{n \pi
           x}{2})\sin(\frac{n \pi y}{2})}{\pi^2 n^2}. \\
\end{array}
$$
\item[4.] $\Gamma = \RR / \ZZ$, $\mu = dx$.  
$$
\begin{array}{lll}
g_\mu(x,y) & = & \frac{1}{2} |x-y|^2 - \frac{1}{2}|x-y| + \frac{1}{12}
       \\  & = & 2\sum_{n\geq 1} 
\left( \frac{\cos(2\pi 
           nx)\cos(2\pi ny)}{4\pi^2 n^2} + 
\frac{\sin(2\pi 
           nx)\sin(2\pi ny)}{4\pi^2 n^2}
\right). \\

\end{array}
$$
\item[5.] $\Gamma = \RR / \ZZ$, $\mu = \delta_0$. 
$$
\begin{array}{lll}
g_\mu(x,y) & = & \left\{ 
\begin{array}{ll}
x(1-y) & {\rm if \;} x<y \\
y(1-x) & {\rm if \;} x\geq y \\
\end{array}
\right.
\\
           & = & 2 \sum_{n\geq 1} \frac{\sin(n \pi x)\sin(n \pi
           y)}{\pi^2 n^2}. \\
\end{array}
$$
\end{itemize}

\subsection{Related Work}
\label{RelatedWorkSection}

The term ``metrized graph'' was introduced by the second author in \cite{RumelyBook}, 
and was further developed by Chinburg and Rumely in \cite{CR} and by Zhang in \cite{Zh1}.  
In those papers, metrized graphs are seen to arise naturally from considerations
in arithmetic geometry, and the function $g_\mu(x,y)$ is related to
certain pairings which occur in Arakelov theory and arithmetic capacity theory.  
The present paper was originally motivated by the study of reduction graphs and
Berkovich spaces associated to algebraic curves, with an eye toward
proving general adelic equidistribution theorems for points of small canonical height
(see \cite{BRequidistr}).  It has been further developed in \cite{RumelyBaker} to yield a 
detailed study of the Laplacian on the Berkovich projective line
over a complete and algebraically closed nonarchimedean field.

\vskip .1 in

It seems clear, however, that the theory developed in the present work should have applications
beyond its origins in arithmetic geometry.  There is already a vast  
literature concerning objects which are more or less the same as metrized graphs.  
Other names one encounters for these objects include networks, metric graphs, and quantum graphs.
We are not expert enough to survey the entire literature and comment on its relation to the
present paper, so out of necessity we will restrict our comments to a few specific related works.  
Metrized graphs have applications to the physical sciences; for example, the papers \cite{Nicaise1} and \cite{Nicaise2}
are written in the context of mathematical biology, and deal with
the spread of potential along the dendrites of a neuron.
The survey paper \cite{Kuchment} of Kuchment, on the other hand, is motivated by considerations in quantum physics.
And the results from \cite{AM1}-\cite{AM3} and \cite{Be1}-\cite{Be2} have applications to various aspects of 
physics and electrical and mechanical engineering.

\vskip .1 in
There is a well-developed mathematical theory dealing with 
differential equations and spectral problems on metrized graphs.   
We note especially the following works.  
In \cite{Roth}, J.~P.~Roth solves the heat equation on a metrized graph, 
and uses this to show the existence of nonisomorphic isospectral 
metrized graphs (``Can you hear the sound of a graph?'').  
In \cite{Be1},  J.~ von Below studies the characteristic equation 
associated to the eigenvalue problem for heat transmission on a metrized
graph, and in the case corresponding to our $\mu = dx$, gives a precise 
description of the eigenvalues and their multiplicities (\cite{Be1}, p.320).  
In \cite{Be2} he studies very general problems of  
Sturm-Liouville type on graphs.  
He shows they have discrete eigenspectra and that the
eigenfunctions yield uniformly convergent series expansions of 
well-behaved functions.  In \cite{AM2}, F.~A.~Mehmeti studies nonlinear
wave equations on CW-complexes of arbitrary dimension,
together with boundary conditions coming from 
higher-dimensional analogues of Kirckhoff's laws.   
By casting the theory in abstract functional-analytic form,  
he obtains very general results about asymptotics of eigenvalues.      
In \cite{Friedlander}, L.~Friedlander shows that for a ``generic'' metrized graph, 
all eigenvalues of the Laplacian (relative to the measure $\mu = dx$) have
multiplicity one.  And in \cite{Friedlander2}, Friedlander gives a sharp lower bound
for the $dx$-eigenvalues of the Laplacian on a metrized graph.

\vskip .1 in

We should also mention here the recent book \cite{FJ} by Favre and Jonsson.  They develop a detailed theory of Laplacians
on (not necessarily finite) metrized trees which is closely related to our approach in the special case where the graph is both finite
and a tree.

\vskip .1 in

Finally, we would like to point out that the present paper is more or less self-contained, and uses less
sophisticated tools from analysis than most of the papers cited above.  In particular, we make no use of the theory of 
Sobolev spaces (which are omnipresent in the literature on networks and quantum graphs), 
although their introduction would probably shorten some of the proofs here.
For the reader seeking an even less technical introduction to the Laplacian on a
metrized graph, a leisurely survey can be found in \cite{BF}.

\vskip .1 in

\section{Metrized Graphs}
\label{MetrizedGraphSection}

\vskip .1 in

\subsection{Definition of a metrized graph}

We begin by giving a rigorous definition of a metrized graph,
following \cite{Zh1}.  
We will then explain how this definition 
relates to the intuitive definition given in the introduction.

\vskip .1 in 

\begin{definition} \label{Def1}
A metrized graph $\Gamma$ is a compact, connected
metric space such that for each $p \in \Gamma$, there exist
a radius $r_p > 0$ and an integer $n_p \ge 1$ such that $p$
has a neighborhood $V_p(r_p)$ isometric to the star-shaped set 
\begin{equation*}
S(n_p,r_p) \ = \ \{ z \in \CC : z = t e^{k \cdot 2 \pi i/n_p} \
     \text{for some $0 \le t < r_p$ and some $k \in \ZZ$} \}, \ 
\end{equation*}
equipped with the path metric.
\end{definition}

There is a close connection between metrized graphs and finite
weighted graphs.  A {\it finite weighted graph} $G$ is a connected,
weighted combinatorial graph equipped with a collection $V = \{
v_1,\ldots,v_n \}$ of vertices, $E = \{ e_1,\ldots,e_m\}$ of edges, 
and $W = \{ w_{ij} \}$ of nonnegative weights ($1\leq i,j \leq n$) satisfying:
\begin{itemize}
\item If $w_{ij}=0$ then there is no edge connecting $v_i$ and $v_j$.
\item If $w_{ij}>0$ then there is a unique edge $e_k$ connecting $v_i$
  and $v_j$.
\item $w_{ij} = w_{ji}$ for all $i,j$.
\item $w_{ii} = 0$ for all $i$.
\end{itemize}

We let $w(e_k)$ denote the weight of edge $e_k$, and define the {\it
  length} of $e_k$ to be $L(e_k) = 1/w(e_k)$.

\vskip .1 in

There is a natural partial ordering on the collection of finite weighted graphs,
where $G' \leq G$ if we can refine $G$ to $G'$ by a sequence of
length-preserving (as opposed to weight-preserving) subdivisions.

We write $G \sim G'$ if there exists a finite weighted graph $G''$ with $G''
\leq G$ and $G'' \leq G'$, i.e., if $G$ and $G'$ have a common
subdivision.  It is easy to see that this defines an equivalence
relation on the collection of finite weighted graphs.



\begin{lemma}
\label{lemma:equivalence}
There is a bijective correspondence between metrized graphs and
equivalence classes of finite weighted graphs.
\end{lemma}

\begin{proof}
Suppose first that $\Gamma$ is a metric space satisfying Definition
\ref{Def1}.  By compactness, it can be covered by a finite number of neighborhoods $V_p(r_p)$.  Hence
there are only finitely many points $p \in \Gamma$ for which $n_p \ne 2$.
We will call this collection of points, or any finite set of points in $\Gamma$
containing it, a vertex set for $\Gamma$.  Let $V$ be a nonempty vertex
set for $\Gamma$.  Then $\Gamma \backslash V$ has finitely many connected
components.  Each of these can be covered by a finite number of open intervals
(if $p$ is a vertex, then $V_p(r_p) \backslash \{ p \}$ is isometric
to a set composed of $n_p$ disjoint open intervals), and so is
isometric to an open interval $(0,L_i)$.  The closure of such a component in
$\Gamma$ is an edge;  it will be isometric to either a segment or a loop.
After adjoining finitely many points to $V$, we can assume that each edge
has $2$ distinct boundary points, and hence is isometric to a closed interval.
The resulting set of vertices and edges defines a finite weighted graph $G$
with weights given by $w(e_k) = 1/L_k$.  It is easy to see that
the equivalence class of $G$ is independent of all choices made.

\vskip .1 in

Conversely, given an equivalence class of 
finite weighted graphs, let $G$ be a representative of this class.
We define a corresponding metrized graph $\Gamma$ as follows.
The points of $\Gamma$ are just the points of
$G$, where each edge of $G$ is considered as a line segment of length
$L_i$.  It is clear from the definition of equivalence that
this set is independent of the choice of $G$.

Choose a parametrization $\psi_i : [a_i,b_i] \rightarrow e_i$ of each
edge $e_i$ of $G$ by a
line segment of length $L_i$, and let the `chart' 
$\varphi_i : e_i \rightarrow [a_i,b_i]$ be the inverse to the
parametrization $\psi_i$.  We will call a function
$\gamma : [a,b] \rightarrow \Gamma$ a {\it piecewise isometric path} if there is a
finite partition $a = c_0 < c_1 < \cdots < c_m = b$ such that for each
$j = 1, \ldots, m$ there is an edge $e_{i(j)}$ such that
$\gamma([c_{j-1},c_j]) \subset e_{i(j)}$, and $\varphi_{i(j)} \circ \gamma(t)$
has constant derivative $\pm 1$ on $[c_{j-1},c_j]$.  For any two points
$p, q \in \Gamma$, put
\[
d(p,q) \ = \ \inf_{\substack{\gamma : [a,b] \rightarrow \Gamma \\ 
                             \gamma(a) = p, \ \gamma(b) = q}} |b-a|,
\]
where the infimum is taken over all piecewise isometric paths from $p$ to $q$.
It is easy to see that $d(p,q)$ is independent of the choice of
parametrizations, that $d(p,q)$ defines a metric on $\Gamma$, and that
for each $p$ there is an $r_p > 0$ such that for any $q \in \Gamma$ with
$d(p,q) < r_p$ there is a unique shortest piecewise isometric path from
$p$ to $q$.  If $p$ belongs to the vertex set of $G$,
let $n_p$ be the valence of $G$ at $p$ (the number of edges of $G$
incident to $p$).  Otherwise, let $n_p = 2$.  Then 
the neighborhood $V_p(r_p) = \{q \in \Gamma : d(p,q) < r_p \}$ is
isometric to the star $S(r_p,n_p) \subset \CC$.
\end{proof}

\vskip .1 in

Henceforth, given a metrized graph $\Gamma$, we will 
frequently choose without comment a corresponding finite weighted
graph $G$ as in Lemma~\ref{lemma:equivalence}, 
together with distinguished parametrizations of the
edges of $G$.  The reader will easily verify that all important definitions below (including the definitions of 
directional derivatives and the Laplacian) are independent of these implicit choices.

By an `isometric path' in $\Gamma$, we will mean a $1-1$ piecewise isometric path.
Since $\Gamma$ is connected, it is easy to see that for all $x,y \in
\Gamma$, there exists an isometric path from $x$ to $y$.
We will say that an isometric path $\gamma : [0,L] \rightarrow \Gamma$
emanates from $p$, and ends at $q$, if $\gamma(0) = p$ and $\gamma(L) = q$.
If $p \in \Gamma$, there are $n_p$ distinguished isometric paths $\gamma_i$
emanating from $p$ such that any other isometric path emanating from $p$
has an initial segment which factors through an initial segment of one 
of the distinguished paths.  We will formally introduce $n_p$ ``unit vectors'' 
$\vv_i$ to describe these paths, and
write $p + t \vv_i $ for $\gamma_i(t)$.  By abuse of language,
we will refer to the distinguished paths $\gamma_i(t) = p + t \vv_i$
as `edges emanating from $p$'.  If $f : \Gamma \rightarrow \CC$ is a function,
and $\vv = \vv_i$ a unit vector at $p$, then we will define
the `(one-sided) derivative of $f$ in the direction $\vv$' to be
\begin{equation*}
  d_{\vv}f(p) \ = \ \lim_{t \rightarrow 0^+}
                 \frac{ f(p + t \vv) - f(p) }{t}
\end{equation*}
provided the limit exists as a finite number.  Write $\VEC(p)$
for the collection of directions emanating from $p$.

\vskip .1 in

\vskip .1 in

\section{Eigenfunctions via comparison of Dirichlet and $L^2$ norms}
\label{EigenfunctionViaPairingsSection}

\vskip .1 in

In this section, we discuss a ``Rayleigh-Ritz''-type method for constructing
eigenfunctions of the Laplacian.  The method is motivated by
a classical procedure for solving the Laplace-Dirichlet eigenvalue
problem (see e.g. \cite[Chapter 4]{PDE}).  

\vskip .1 in

Recall that the $L^2$ inner product for $f, g \in L^2(\Gamma)$ is
$\<f,g\> = \<f,g\>_{L^2} = \int_{\Gamma} f(x) \overline{g(x)} \, dx$, 
with associated norm $\| \cdot \|_2 = \<f,f\>^{1/2}$.  

\medskip

Let $\mu$ be a real-valued signed Borel measure on $\Gamma$ with
$\mu(\Gamma) = 1$ and $|\mu|(\Gamma) < \infty$.  
Recall from \S\ref{PreliminaryDirichletSection} that 
$\Zh(\Gamma)$ denotes the space of all continuous
functions $f : \Gamma \to \CC$ such that 
$f$ is piecewise $\cC^2$ and $f'' \in L^1(\Gamma)$, and
$\Zh_\mu(\Gamma) = \{ f \in \Zh(\Gamma) \; : \; \int_\Gamma
f \, d\mu = 0 \}.$  
Also, $\Dir_\mu(\Gamma)$ denotes the completion of $\Zh_\mu(\Gamma)$ with
respect to the Dirichlet norm
$\| \cdot \|_{\Dir}$, and $L^2_\mu(\Gamma)$ denotes the completion of $\Zh_\mu(\Gamma)$ with
respect to the $L^2$-norm.  

\medskip

\begin{lemma}
\label{lemma:PI}
For all $f \in \Zh_\mu(\Gamma)$, we have the estimates
$$
\| f \|_{\infty} \leq \ell(\Gamma)^{\frac{1}{2}} \cdot |\mu|(\Gamma)
\cdot \| f \|_{\Dir},
$$
and the \emph{Poincar{\'e} Inequality}
$$
\| f \|_2 \leq \ell(\Gamma) \cdot |\mu|(\Gamma) \cdot \| f \|_{\Dir},
$$
where $\ell(\Gamma) = \int_\Gamma dx$ is the total length of $\Gamma$.

In particular, there exist constants $C_\infty,C_2 > 0$ $($depending
only on $\mu$ and $\Gamma$$)$ such that
$\| f \|_\infty \leq C_\infty \| f \|_{\Dir}$ and 
$\| f \|_2 \leq C_2 \| f \|_{\Dir}$ for all $f \in \Zh_\mu(\Gamma)$.
\end{lemma}

\begin{proof}
Fix a point $x_0 \in \Gamma$, let $x \in \Gamma$ be another point,
and let $\gamma$ be an isometric path from $x_0$ to $x$. 

Since $f$ is continuous and piecewise $\cC^1$, it follows from the fundamental theorem of calculus
that $f(x) - f(x_0) = \int_\gamma f'(t) \, dt$,
so that 
\[
|f(x) - f(x_0)| \leq \int_\gamma |f'(t)| \, dt \leq \int_\Gamma
|f'(t)| \, dt \leq \ell(\Gamma)^{1/2} \| f \|_{\Dir}
\]
by the Cauchy-Schwarz inequality.
Since $f \in \Zh_{\mu}(\Gamma)$ we have $\int_{\Gamma} f(x) \, d\mu(x) = 0$,
so $\int_\Gamma \left( f(x) - f(x_0) \right) \, d\mu(x) = -f(x_0)$, and
\begin{equation*}
\begin{aligned}
|f(x_0)| = |\int_\Gamma \left( f(x) - f(x_0) \right) \, d\mu(x)|
& \leq \int_\Gamma |f(x) - f(x_0)| \, d|\mu|(x) \\ 
& \leq \ell(\Gamma)^{1/2} |\mu|(\Gamma) \| f \|_{\Dir}.
\end{aligned}
\end{equation*}
This holds for each $x_0 \in \Gamma$, so
$\| f \|_{\infty} \leq \ell(\Gamma)^{\frac{1}{2}} |\mu|(\Gamma)  \| f
\|_{\Dir}$ and
\[
\| f \|_2^2 =  \int_\Gamma |f(x)|^2 dx 
\leq  \ell(\Gamma)^2  |\mu|(\Gamma)^2 \| f \|_{\Dir}^2 
\]
as desired.
\end{proof}

\begin{corollary}
A sequence in $\Zh_\mu(\Gamma)$ which is Cauchy with respect to the Dirichlet
norm is also Cauchy with respect to the $L^2$-norm.
In particular, there is a natural inclusion map $\Dir_\mu(\Gamma) \subseteq L^2_\mu(\Gamma)$.
\end{corollary}

\vskip .1 in

\begin{corollary} \label{Cor26}
We can uniquely identify $\Dir_\mu(\Gamma)$ with a linear subspace
of $\cC(\Gamma)$ containing $\Zh_{\mu}(\Gamma)$
in such a way that for each $f \in \Dir_\mu(\Gamma)$,
$\| f \|_{\infty} \ \le \ C_\infty \| f \|_{\Dir}$.
Under this identification, each $f \in \Dir_\mu(\Gamma)$ satisfies
$\int_{\Gamma} f(x) \, d\mu(x) = 0$.
\end{corollary}             

\begin{proof}
Take $f \in \Dir_\mu(\Gamma)$.
Let $f_1, f_2, \ldots \in \Zh_{\mu}(\Gamma)$
be a sequence of functions with $\|f-f_n\|_{\Dir} \rightarrow 0$.
Then for each $\varepsilon > 0$, there is an $N$ 
such that for all $n, m \ge N$,  $\|f_n-f_m\|_{\Dir} < \varepsilon$.

By Lemma~\ref{lemma:PI},
$\|f_n-f_m\|_{\infty} \ \le \ C_\infty \|f_n-f_m\|_{\Dir}$,
so the functions $f_n$ converge uniformly to a function $F(x)$.
Since each $f_n(x)$ is continuous,  $F(x)$ is continuous.
Furthermore $\int_{\Gamma} F(x) \, d\mu(x) = 0$
since $\int_{\Gamma} f_n(x) \, d\mu(x) = 0$ for each $n$.  

It is easy to see that $F(x)$ is independent of the sequence $\{f_n\}$
in $\Zh_{\mu}(\Gamma)$  converging to $f$.
(If $\{h_n\}$ is another such sequence, apply the argument above to
the sequence $f_1, h_1, f_2, h_2, \ \ldots$.)  In particular,
$F(x) = 0$ if $f = 0$ in $\Dir_\mu(\Gamma)$.
Thus, there is a natural injection
$\iota_{\mu} : \Dir_{\mu}(\Gamma) \rightarrow \cC(\Gamma)$.
We use this to identify $f$ with $F(x) = \iota_{\mu}(f)$.

We claim that $\| F \|_{\infty} \le C_\infty \| f \|_{\Dir}$.
Indeed, if $\{f_n\}$ is a sequence of functions in $\Zh_{\mu}(\Gamma)$
converging to $f$, then
\[
\|F\|_{\infty} =  \lim_{n \rightarrow \infty} \|f_n\|_{\infty} 
       \le  \lim_{n \rightarrow \infty} C_\infty \|f_n \|_{\Dir} 
      \ = \ C_\infty \| f \|_{\Dir}.
\]

Under the identification of $f$ with $\iota_{\mu}(f)$,
functions in $\Zh_{\mu}(\Gamma)$
are taken to themselves, since $f \in \Zh_{\mu}(\Gamma)$
is the limit of the constant sequence $\{f,f,f,\ldots\}$.

We leave the uniqueness assertion to the reader.
\end{proof}

\vskip .1 in

Henceforth we will identify elements of $\Dir_{\mu}(\Gamma)$ with continuous
functions as in Corollary \ref{Cor26}.
Since $\Zh_\mu(\Gamma)$ is dense in $\Dir_\mu(\Gamma)$,
we easily deduce the following stronger version of
Lemma~\ref{lemma:PI}.

\begin{corollary}
\label{cor:strongerPI}
If $f \in \Dir_\mu(\Gamma)$, then 
$\| f \|_\infty \leq C_\infty \| f \|_{\Dir}$ 
and 
$\| f \|_2 \leq C_2 \| f \|_{\Dir}$.
\end{corollary}

The next lemma will be used in the proof of Lemma~\ref{lemma:weaktostrong}.
For the statement, recall that a sequence $\{ f_n \}$ in a Hilbert space $H$ 
{\it converges weakly} to $f \in H$ if 
$$
\< f_n,g \> \to \< f,g \>
$$
for all $g \in H$.

\begin{lemma}
\label{lemma:weak-j}
Suppose $f_n \in \Zh_\mu(\Gamma)$ and $f_n \to 0$ weakly in $\Dir_\mu(\Gamma)$.  Then for every subset
$S$ of $\Gamma$ isometric to a closed interval, we have
$$
\int_S f'_n(x) dx \to 0.
$$
\end{lemma}

\begin{proof}
Assume that $S$ is parametrized by the segment $[a,b]$.  
Then by the fundamental theorem of calculus and the defining relation
$\Delta_x j_\zeta(x,y) = \delta_y(x) - \delta_\zeta(x)$
for $j_\zeta(x,y)$, we have
$$
\begin{array}{lll}
\int_S f_n'(x) dx & = & f_n(b) - f_n(a) \\
& = & \int_\Gamma f_n(x) \, d\left( \delta_b(x) - \delta_a(x) \right) \\
& = & \int_\Gamma f_n(x) \, d\left( \Delta_x j_a(x,b) \right) \\
& = & \int_\Gamma f_n'(x) \left( \frac{d}{dx} j_a(x,b) \right) dx, \\
\end{array}
$$
where the last equality follows from Proposition~\ref{SelfAdjointProp}.

Since $j_a(x,b) + c_\mu \in \Zh_\mu(\Gamma)$ for some constant $c_\mu
\in \RR$, the hypothesis of weak convergence now gives us what we want.
\end{proof}

\vskip .1 in

{\bf Remark:} 
When $\Gamma$ contains a cycle, there need not be a function $g_S$ in 
$\Zh_\mu(\Gamma)$ whose derivative is the characteristic function of
$S$.  When such a function $g_S$ exists, however, one can argue more
simply that
$$
\int_S f'_n(x) \, dx = \int_\Gamma f'_n(x) g_S'(x) \, dx \to \int_\Gamma 0 \, dx = 0.
$$

\vskip .1 in

The following lemma is analogous to Rellich's Theorem in classical
analysis:

\begin{lemma}
\label{lemma:weaktostrong}
If $f_n \to f$ weakly in $\Dir_\mu(\Gamma)$ and the sequence $\| f_n
\|_{\Dir}$ is bounded, then 
$f_n \to f$ strongly in $L^2_\mu(\Gamma)$.
In particular, $\| f_n \|_2 \to \| f \|_2$.
\end{lemma}

\begin{proof}

Choose a sequence $g_n$ of functions in $\Zh_\mu(\Gamma)$ such that 
$\| g_n - f_n \|_{\Dir} \to 0$.  If we can prove the lemma for the
sequence $g_n$, then the corresponding result for $f_n$ follows easily
using Corollary~\ref{cor:strongerPI}.  We may therefore assume, without loss of
generality, that each $f_n \in \Zh_\mu(\Gamma)$.

We want to show that 
$\int_\Gamma |f_n(x) - f(x)|^2 dx \to 0$.
By Lebesgue's dominated convergence theorem, 
it is enough to show that $f_n \to f$ pointwise, 
and that there exists $M>0$ such that $| f_n(x) | \leq M$ for all $x \in \Gamma, n \geq 1$.
The latter assertion is clear from Lemma~\ref{lemma:PI}, since the
Dirichlet norms of the $f_n$'s are assumed to be bounded.  So it
remains to prove that $f_n$ converges pointwise to $f$.  

Let $h_n = f_n - f$.  
As in the proof of Lemma~\ref{lemma:PI}, fix a point $x_0 \in
\Gamma$, let $x\in \Gamma$, and choose an isometric path $\gamma$ from
$x_0$ to $x$.
Then by the fundamental theorem of calculus, we have 
$h_n(x) = h_n(x_0) + \int_\gamma h'_n(t) \, dt$
for all $n$.  

As $h_n \to 0$ weakly in $\Dir_\mu(\Gamma)$, we see by Lemma~\ref{lemma:weak-j} 
that
$$
|h_n(x) - h_n(x_0)| = |\int_\gamma h'_n(t) \, dt | \to 0,
$$
so that $h_n(x) - h_n(x_0) \to 0$ as $n\to\infty$.

On the other hand, 
since the functions $h_n(x)$ are uniformly bounded on $\Gamma$, 
it follows by Lebesgue's Dominated Convergence Theorem that
$$
\int_\Gamma \left( h_n(x) - h_n(x_0) \right) 
d\mu(x) \to 0.
$$

As $\int_\Gamma h_n(x) \, d\mu(x) = 0$, this shows that $h_n(x_0) \to
0$ as $n\to \infty$ as desired.
\end{proof}

\medskip

%

\medskip

We now define what it means to be an eigenfunction of the 
Laplacian on $\Dir_\mu(\Gamma)$.  

\begin{definition} \label{DefDirDF} 
A nonzero function $f$ is an {\it eigenfunction} for $\Delta$ acting on
$\Dir_\mu(\Gamma)$ if there exists a constant $\lambda \in \CC$ (the {\it
  eigenvalue}) such that for all $g \in \Dir_\mu(\Gamma)$ 
\begin{equation}
\label{eigenequation}
\< f,g \>_{\Dir} = \lambda \< f,g \> \ .
\end{equation}
\end{definition}

It is clear (setting $g = f$) that every eigenvalue of $\Delta$ is
real and positive, and we will write $\lambda = \gamma^2 > 0$.


\begin{theorem}
\label{MainPairingTheorem} { \ }

$A)$ Each eigenvalue $\lambda$ of $\Delta$ on $\Dir_\mu(\Gamma)$ 
is nonzero and occurs with finite
multiplicity.  If we write the eigenvalues as $0 < \lambda_1 \leq
\lambda_2 \leq \cdots$ (with corresponding eigenfunctions
$f_1,f_2,\ldots$), then $\lim_{n\to\infty} \lambda_n = \infty$.

$B)$ Suppose the eigenfunctions $f_n$ are normalized so that $\|f_n\|_2 =
1$.  Then $\{ f_n \}$ forms a basis for $\Dir_\mu(\Gamma)$, and an orthonormal basis for $L^2_\mu(\Gamma)$.

\end{theorem}

We will not give a proof of Theorem~\ref{MainPairingTheorem}, 
since the ideas involved are well-known (see e.g. \cite[Chapter 4]{PDE}, \cite{Kuchment}).
However, it seems worthwhile to at least point out that the ``Rayleigh-Ritz''
approach to proving Theorem~\ref{MainPairingTheorem} is based on
the following useful characterization of the leading eigenvalue $\lambda_1$.

\begin{theorem}
\label{RayleighRitzTheorem}
Let $S_1 = \{ f \in \Dir_\mu(\Gamma) \; : \; \| f \|_2 = 1 \}$
be the $L^2$ unit sphere in $\Dir_\mu(\Gamma)$, let
$\gamma_1 = \inf_{f \in S_1} \| f \|_{\Dir}$,
and let $\lambda_1 = \gamma_1^2$. 
Then $\lambda_1$ is the smallest eigenvalue of $\Delta$ on $\Dir_\mu(\Gamma)$.
\end{theorem}

The key ingredient needed to prove Theorem~\ref{RayleighRitzTheorem} along
the lines of \cite[Chapter 4]{PDE} is the estimate found in Lemma~\ref{lemma:weaktostrong}.

\vskip .1 in
Before moving on, it may be useful to note how eigenvalues 
and eigenfunctions behave under scaling in the underlying graph.  
If $\Gamma$ is a metrized graph, and $0 < \beta \in \RR$, we can define
a new metrized graph $\Gamma(\beta)$ whose underlying point set is the 
same as $\Gamma$, but in which all distances are multiplied by $\beta$.  
Write $d(x,y)$ for the metric on $\Gamma$, and $d_{\beta}(x,y)$ for the metric
on $\Gamma(\beta)$.  There is a natural isomorphism 
$\sigma_{\beta} : \Gamma(\beta) \rightarrow \Gamma$ such that for all 
$x, y \in \Gamma(\beta)$, 
\begin{equation*}
d_{\beta}(x,y) \ = \ \beta \cdot d(\sigma_{\beta}(x),\sigma_{\beta}(y)).
\end{equation*}
If $\mu$ is a measure on $\Gamma$, 
then $\mu_{\beta} = \mu \circ \sigma_{\beta}$ is 
a measure on $\Gamma(\beta)$.  

\begin{proposition} \label{ScalingProp}
Let $\Gamma$ be a metrized graph, 
and let $\mu$ be a real, signed Borel measure on $\Gamma$ of total mass $1$.
Take $\beta > 0$.  Then there is a one-one
correspondence between the eigenvalues $\lambda_n$ of $\Delta$ 
on $\Dir_{\mu}(\Gamma)$ and the eigenvalues $\lambda_n(\beta)$ of $\Delta$
on $\Dir_{\mu_{\beta}}(\Gamma(\beta))$, given by 
\begin{equation*}
\lambda_n(\beta) \ = \ \frac{1}{\beta^2}  \lambda_n.
\end{equation*}
Furthermore, $f \in \Dir_{\mu}(\Gamma)$ is an eigenfunction for $\Delta$ with
eigenvalue $\lambda$ if and only if 
$f \circ \sigma_{\beta} \in \Dir_{\mu_{\beta}}(\Gamma(\beta))$ 
is an eigenfunction for $\Delta$ with eigenvalue $\frac{1}{\beta^2}\lambda$.
\end{proposition}

\begin{proof}
This follows immediately from the definitions, if one notes that the Dirichlet
and $L^2$ pairings on $\Gamma$ and $\Gamma(\beta)$ are related by 
\begin{eqnarray*}
\<f \circ \sigma_{\beta},g \circ \sigma_{\beta}\>_{\Dir(\Gamma(\beta))}  
        & = & \frac{1}{\beta} \cdot \<f,g\>_{\Dir(\Gamma)} \ ,\\
\<f \circ \sigma_{\beta},g \circ \sigma_{\beta}\>_{L^2(\Gamma(\beta))}  
        & = & \beta \cdot \<f,g\>_{L^2(\Gamma)} \ ,
\end{eqnarray*}
as follows from the chain rule and the change of variables formula.  
\end{proof} 


\vskip .1 in

\section{The Classes $\cD(\Gamma)$ and $\BDV(\Gamma)$}
\label{MeasureSection}

\vskip .1 in

Let $\cD(\Gamma)$ be
the class of all functions on $\Gamma$ whose one-sided derivatives
exist everywhere, i.e.,
\begin{equation*}
\cD(\Gamma) = \{ f : \Gamma \rightarrow \CC : \
          \text{ $d_{\vv}f(p)$ exists for each $p \in \Gamma$
                                    and $\vv \in \VEC(p)$} \} \ .
\end{equation*}
It is easy to see that
each $f \in \cD(\Gamma)$ is continuous.

Let $\cA$ be the Boolean algebra of subsets of $\Gamma$ generated by the
{\emph{connected}} open sets.  
Each $S \in \cA$ is a finite disjoint union of sets
isometric to open, half-open, or closed intervals (we consider isolated points
to be degenerate closed intervals);  conversely, all such sets belong to $\cA$.
Define the set $b(S)$ of boundary points of $S$ to be the collection of  
points belonging to the closures of both $S$ and $\Gamma \backslash S$.
It is easy to see that each $S \in \cA$ has only finitely many boundary points.
Note that under this definition, if $\Gamma = [0,1]$ and 
$S = [0,\frac{1}{2}]$, for example,
then the left endpoint $0$ is not a boundary point of $S$.

For each $p \in b(S)$, let $\IN(p,S)$ be the set of
`inward-directed unit vectors at $p$': the set of all 
$\vv \in \VEC(p)$ for which
$p+t \vec v$ belongs to $S$ for all sufficiently small $t > 0$.  
Let $\OUT(p,S) = \VEC(p) \backslash \IN(p,S)$ be the collection of 
`outward-directed unit vectors at $p$'.  If $p$ is an isolated point of $S$,
then $\IN(p,S) = \emptyset$ and $\OUT(p,S) = \VEC(p)$.

If $f \in \cD(\Gamma)$, then we can define a finitely additive set
function $m_f$ on $\cA$ by requiring that for each $S \in \cA$
\begin{equation*}
m_f(S) \ = \
   \sum_{ \substack{ b \in b(S) \\ b \notin S}}
                 \sum_{\vv \in \IN(b,S)} d_{\vv}f(p)
   - \sum_{\substack{b \in b(S) \\ b \in S }}
                 \sum_{\vv \in \OUT(b,S)} d_{\vv}f(p) \ .
\end{equation*}
Thus, for an open set $S \in \cA$,
\begin{equation*}
m_f(S) \ = \
      \sum_{b \in b(S)} \sum_{\vv \in \IN(b,S)} d_{\vv}f(p) \ ,
\end{equation*}
for a closed set $S \in \cA$,
\begin{equation*}
m_f(S) \ = \
   - \sum_{b \in b(S)} \sum_{\vv \in \OUT(b,S)} d_{\vv}f(p) \ ,
\end{equation*}
and for a set $S = \{p\}$ consisting of a single point,
\begin{equation*}
m_f(\{p\}) \ = - \sum_{\vv \in \VEC(p)} d_{\vv}f(p) \ .
\end{equation*}
The finite additivity is clear if one notes that each $S \in \cA$ can be
written as a finite disjoint union of open intervals and points,
and that for each boundary point $p$ of $S$, the set $\IN(p,S)$
coincides with $\OUT(p,\Gamma \backslash S)$.

Note that
\begin{equation*}
m_f(\emptyset) \ = \ m_f(\Gamma) \ = \ 0 \ ,
\end{equation*}
since by our definition, 
both the empty set and the entire graph $\Gamma$ have no boundary points.  
It follows for any $S \in \cA$
\begin{equation*}
m_f(\Gamma \backslash S) \ = \ -m_f(S) \ .
\end{equation*}
If $f_1, f_2 \in \cD$ and $c_1, c_2 \in \CC$, then it is easy to see that
the set function corresponding to $c_1 f_1 + c_2 f_2$ is
\begin{equation*}
m_{c_1 f_1 + c_2 f_2} = c_1 m_{f_1} + c_2 m_{f_2} \ .
\end{equation*}

\vskip .1 in
We will say that a function $f \in \cD(\Gamma)$ is of 
``bounded differential variation'', and write $f \in \BDV(\Gamma)$, 
if there is a constant $B > 0$ such that for any countable collection
$\cF$ of pairwise disjoint sets in $\cA$,
\begin{equation*}
\sum_{S_i \in \cF} |m_f(S_i)| \ \le \ B \ .
\end{equation*}
It is easy to see that $\BDV(\Gamma)$ is a $\CC$-linear subspace
of $\cD(\Gamma)$.

\begin{proposition} \label{Prop1}  Let $\Gamma$ be a metrized graph.

A)  If $f \in \BDV(\Gamma)$, then there are only countably many points
     $p_i \in \Gamma$ for which $m_f({p_i}) \ne 0$, and $\sum_i |m_f(p_i)|$ 
     converges.

B)  In the definition of $\BDV(\Gamma)$ one can restrict to families
    of pairwise disjoint connected open sets, or connected closed sets.
    More precisely, a function $f \in \cD(\Gamma)$ belongs to $\BDV(\Gamma)$
    if and only if

\quad $(1)$ there is a constant $B_1$ such that for any countable family $\cF$
     of pairwise disjoint connected open sets $V_i \in \cA$,
\begin{equation*}
     \sum_{\substack{ V_i \in \cF \\ \text{$V_i$ connected, open} }}
          |m_f(V_i)| \ \le \ B_1 \ ; \quad \text{or}
\end{equation*}

\quad $(2)$ there is a constant $B_2$ such that for any countable family $\cF$
     of pairwise disjoint connected closed sets $E_i \in \cA$,
\begin{equation*}
     \sum_{\substack{ E_i \in \cF \\ \text{$E_i$ connected, closed} }}
          |m_f(E_i)| \ \le \ B_2 \ .
\end{equation*}

\end{proposition}

\begin{proof} For A), note that for each $0 < n \in \ZZ$ there can be
at most $n B = B/(1/n)$ points $p \in \Gamma$ with $|m_f(\{p\})| \ge 1/n$,
so there are at most countably many points with $m_v(\{p\}) \ne 0$.
Moreover,
\begin{equation*}
\sum_{\substack{ p \in \Gamma \\ m_f(\{p\}) \ne 0 }}
               |m_f(\{p\})| \ < \ B \ .
\end{equation*}

For B), note that if $f \in \BDV(\Gamma)$ then $B(1)$
and $B(2)$ hold trivially.  Conversely, first
suppose $B(1)$ holds.  We begin by showing that
there are only countably many points $p \in \Gamma$ for which
$m_f(p) \ne 0$, and that
\begin{equation} \label{FGH1}
\sum_{\substack{ m_f(\{p\}) \ne 0 }} |m_f(\{p\})| \ < \ 2 B_1 \ .
\end{equation}
To see this, suppose that for some $0 < n \in \ZZ$,
there were more than $2 n B_1$ points $p$ with $|m_v(\{p\})| \ge 1/n$.
Then either there would be more than $n B_1$ points with $m_v(\{p\}) \ge 1/n$,
or more than $n B_1$ points with $m_v(\{p\}) \le -1/n$.  Suppose for example
that $p_1, \ldots, p_r$ satisfy $m_v(\{p_i\}) \le -1/n$, where $r > nB_1$.
Put $K = \{p_1, \ldots, p_r\}$, and let $V_1, \ldots, V_s$ be the connected
components of $\Gamma \backslash K$.  Then
\begin{eqnarray*}
B_1 & < & -(m_f(\{p_1\}) + \cdots + m_f(\{p_r\})) \ = \ -m_f(K) \\
    & = & m_f(V_1 \cup \cdots \cup V_s) \ = \ m_v(V_1) + \cdots + m_v(V_s) \\
    & \le & |m_f(V_1)| + \cdots + |m_f(V_r)| \ \le \ B_1 \ .
\end{eqnarray*}
Hence there can be only countably many points $p_i$ with
$m_f(\{p_i\}) \ne 0$, and (\ref{FGH1}) holds.

Now let $\cF$ be any countable collection of pairwise disjoint sets
$S_i \in \cA$.  Each $S_i$ can be decomposed as a finite disjoint union
of connected open sets and sets consisting of isolated points, and if
\begin{equation*}
S_i \ = \  V_{i1} \cup \cdots \cup V_{i,r_i}
              \cup \{p_{i1}\} \cup \cdots \cup \{p_{i,s_i}\}
\end{equation*}
is such a decomposition, then
\begin{equation*}
\sum_{S_i \in \cF} |m_f(S_i)| \ \le \
     \sum_{i,j} |m_f(V_{ij})| + \sum_{i,k} |m_f(\{p_{ik}\})| \ \le \ 3 B_1 \ .
\end{equation*}
It follows that $f \in \BDV(\Gamma)$.

Next suppose that $B(2)$ holds;  we will show that
$B(1)$ holds as well.  Consider a countable
collection $\cF$ of pairwise disjoint connected open sets $V_i$.  Decompose
$\cF = \cF_+ \cup \cF_-$, where $V_i \in \cF_+$ iff $m_f(V_i) \ge 0$,
and $V_i \in \cF_-$ iff $m_f(V_i) < 0$.  Relabel the sets
so that $\cF_+ = \{V_1, V_3, V_5, \ldots \}$.  For each $n$,
put $K_n = \Gamma \backslash (V_1 \cup V_3 \cup \cdots \cup V_{2n+1})$,
and decompose $K_n$ as a finite disjoint union of connected closed sets
$E_1, \cdots, E_r$.  Then
\begin{eqnarray*}
\sum_{i=0}^n m_f(V_{2i+1}) & = & m_f(\Gamma \backslash K_n) \ = \ |m_f(K_n)| \\
   & \le & \sum_{j=1}^r |m_f(E_j)| \ \le \ B_2 \ .
\end{eqnarray*}
Letting $n \rightarrow \infty$ we see that
$\sum_{V_i \in \cF_+} m_f(V_i) \le B_2$.
Similarly $\sum_{V_i \in \cF_-} |m_f(V_i)| \le B_2$,
so $\sum_{V_i \in \cF} |m_f(V_i)| \le 2 B_2$, and we are done.
\end{proof}

\vskip .1 in
We now come to the main result of this section.

\begin{theorem} \label{Thm2}
If $f \in \BDV(\Gamma)$, then the finitely additive set function $m_f$
extends to a bounded complex measure $m_f^*$, with total mass $0$,
on the $\sigma$-algebra of Borel sets of $\Gamma$.
\end{theorem}

\begin{proof}
We begin with a reduction.  It suffices to show that the restriction
of $m_f$ to each edge $e_i$ extends to a Baire measure on $e_i$.
Identifying $e_i$ with its parametrizing interval,
we can assume without loss of generality that $\Gamma = [a,b]$ is a closed interval
and that $f: [a,b] \rightarrow \CC$ is in $\BDV([a,b])$.

We next decompose $m_f$ into an atomic and an atomless part.

Let $\{p_1, p_2, \ldots \}$ be the points in $[a,b]$
for which $m_f(\{p_i\}) \ne 0$.  For brevity, write $c_i = m_f(\{p_i\})$;
by hypothesis, $\sum_i |c_i|$ converges.  Define $g : \RR \rightarrow \CC$ by
\begin{equation*}
g(t) \ = \ -\frac{1}{2} \sum_i  c_i \, |t-p_i| \ ,
\end{equation*}
and let $m_g^* = \sum_i c_i \delta_{p_i}(x)$ be the Borel measure
giving the point mass $c_i$ to each $p_i$.

By a direct computation, one checks that for each $x \in \RR$,
both one-sided derivatives $g_{\pm}^{\prime}(x)$ exist, with
\begin{equation*}
g_{-}^{\prime}(x) = \frac{1}{2}(\sum_{x \le p_i} c_i - \sum_{x > p_i} c_i),
\qquad
g_{+}^{\prime}(x) = \frac{1}{2}(\sum_{x < p_i} c_i - \sum_{x \ge p_i} c_i) \ .
\end{equation*}
For any closed subinterval $[c,d] \subset [a,b]$,
\begin{eqnarray*}
m_g([c,d]) & = & g_{-}^{\prime}(c) - g_{+}^{\prime}(d) \\
& = & \sum_{c \le p_i \le d} c_i \ = \ m_g^*([c,d]).
\end{eqnarray*}
In particular the set function $m_g$ has precisely the same point masses
as $m_f$.  Similar computations apply to open and half-open intervals.
Thus the measure $m_g^*$ extends $m_g$.

Replacing $f(t)$ by $h(t) = f(t)-g(t)$, we are reduced to the situation
where $h \in \BDV([a,b])$ has no point masses.  This means that for each
$p \in (a,b)$,
\begin{equation*}
0 \ = \ m_h(\{p\}) \ = \ h_{-}^{\prime}(p) - h_{+}^{\prime}(p),
\end{equation*}
or in other words that $h^{\prime}(p)$ exists.  By hypothesis,
both $h_{+}^{\prime}(a)$ and $h_{-}^{\prime}(b)$ exist, so $h(t)$ is
differentiable on $[a,b]$.  The fact that $h \in \BDV([a,b])$ means
that $h^{\prime}(t)$ is of bounded total variation.

We claim that $h^{\prime}(t)$ is continuous on $[a,b]$.
Suppose it were discontinous at some point $p$.  By Rudin (\cite{Rud}, p.109) 
the existence of $h^{\prime}(t)$ for all $t$ means that $h^{\prime}(t)$ 
cannot have a jump discontinuity, so the discontinuity must be due to
oscillation.  (Rudin states the result for real-valued functions;
apply it to the real and imaginary parts of $h(t)$ separately.) 
Hence there would be an $\varepsilon > 0$, 
and a sequence of points $t_1, t_2, \ldots$, either monotonically 
increasing or monotonically decreasing to $p$, such that
$|h^{\prime}(t_i)-h^{\prime}(t_{i+1})| > \varepsilon$  for each $i$.
Assume for convenience that the $t_i$ are monotonically increasing.
Considering the intervals $[t_1,t_2], [t_3,t_4], \ldots$ we see that
\begin{equation*}
\sum_{i=1}^{\infty} |m_h([t_{2i-1},t_{2i}])| \ = \ \infty \ ,
\end{equation*}
contradicting the fact that $h \in \BDV([a,b])$.

For $a \le t \le b$, let $T(t)$ be the cumulative total variation function 
of $h^{\prime}(t)$:  letting $\cQ$ vary over all finite partitions 
$a = q_0 < q_1 < \cdots < q_m = t$ of $[a,t]$, 
\begin{equation*}
T(t) \ = \ \sup_{\cQ} \sum_{j=1}^m |h^{\prime}(q_j)-h^{\prime}(q_{j-1})| \ .
\end{equation*}
Then $T : [a,b] \rightarrow \RR$ is monotone increasing,
and is continuous since $h^{\prime}(t)$ is continuous.  By Royden
(\cite{Roy}, Proposition 12, p.301), there is a unique bounded Borel measure
$\nu$ on $[a,b]$ such that $\nu((c,d]) = T(d)-T(c)$ for each half-open
interval $(c,d] \subset [a,b]$.  Since $T(t)$ is continuous, $\nu$ has
no point masses.  

Next put  
\begin{equation*}
T_1(t) = T(t) - \Re(h^{\prime}(t)), 
\qquad T_2(t) = T(t) - \Im(h^{\prime}(t)) \ .
\end{equation*}
Then $T_1(t)$ and $T_2(t)$ are also monotone increasing and continuous, 
so by the same Proposition there are bounded Borel measures $\nu_1$, $\nu_2$ 
on $[a,b]$ such that for each half-open interval $(c,d] \subset [a,b]$,
\begin{equation*}
\nu_1((c,d]) = T_1(d)-T_1(c), \qquad \nu_2((c,d]) = T_2(d)-T_2(c) \ .
\end{equation*}
As before, $\nu_1$ and $\nu_2$ have no point masses.

Now define the complex Borel measure $m_h^* = (\nu_1-\nu) + (\nu_2-\nu) i$.  
By construction, $m_h^*$ has finite total mass.  
For any half-open interval $(c,d] \subset [a,b]$,
\begin{equation*}
m_h^*((c,d]) \ = \ h^{\prime}(c) - h^{\prime}(d) \ = \ m_h((c,d]) \ . 
\end{equation*}
Since the measures $\nu$, $\nu_1$, and $\nu_2$ have no point masses,
\begin{equation*}
m_h^*([c,d]) \ = \ m_h^*((c,d)) \ = \ m_h^*([c,d)) \ = \ m_h^*((c,d]) \ .  
\end{equation*}
Thus $m_h^*$ extends the finitely additive set function $m_h$.   
  
Finally, adding the measures $m_g^*$ extending $m_g$ 
and $m_h^*$ extending $m_h$,
we obtain the desired measure $m_f^* = m_g^* + m_h^*$ extending $m_f$.
Since $m_f(\Gamma) = 0$, also $m_f^*(\Gamma) = 0$.
\end{proof}

\vskip .1 in

\section{Laplacians}

\vskip .1 in

In this paper, for a function $f \in \BDV(\Gamma)$
we define the Laplacian $\Delta(f)$ to be the measure
given by Theorem \ref{Thm2}
\begin{equation*}
\Delta(f) \ = \ m_f^* \ .
\end{equation*}
Sometimes we will need to apply the Laplacian 
to a function of several variables, fixing all but one of them.  
In this case we write $\Delta_x(f(x,y))$ 
for the Laplacian of the function $F_y(x) = f(x,y)$.

\vskip .1 in
We will now show that the Laplacian on $\BDV(\Gamma)$ agrees (up to sign) 
with the Laplacians defined by Chinburg and Rumely (\cite{CR}) 
and Zhang (\cite{Zh1}) on their more restricted classes of functions.  
If $\Delta_{\CR}$ denotes the Chinburg-Rumely Laplacian on the space 
$\CPA(\Gamma)$ of continuous,
piecewise affine functions on $\Gamma$, then it is easy to verify that
if $f \in \CPA(\Gamma)$, then $f \in \BDV(\Gamma)$ 
and $\Delta_{\CR}(f) = - \Delta(f)$.

\vskip .1 in

Recall that the Zhang space $\Zh(\Gamma)$ is the set of all continuous functions
$f : \Gamma \rightarrow \CC$ such that

\quad a) there is a finite set of points $X_f \subset \Gamma$
      such that $\Gamma \backslash X_f$ is a finite union of open intervals,
      the restriction of $f$ to each of those intervals is $\cC^2$, and

\quad b) $f^{\prime \prime}(x) \in L^1(\Gamma)$.

\noindent{Here} we write $f^{\prime \prime}(x)$ for 
$\frac{d^2}{dt^2} f(p+t\vv)$, if $x = p+t\vv \in \Gamma \backslash X_f$.  

\begin{lemma}
\label{ZhangBDVLemma}
The space $\Zh(\Gamma)$ is a subset of $\BDV(\Gamma)$.  
\end{lemma}

\begin{proof}
We first prove that the directional derivatives $d_{\vv}f(p)$ exist for each
$p \in X_f$ and each $\vv \in \VEC(p)$.  Fix such a $p$ and $\vv$, and let
$t_0$ be small enough that $p + t \vv \in \Gamma \backslash X_f$ for
all $t \in (0,t_0)$.  By abuse of notation, we will write $f(t)$ for
$f(p+t \vv)$.  By hypothesis, $f \in \cC^2((0,\delta))$.
By the Mean Value Theorem, $d_{\vv}f(p)$ exists if and only if 
$\lim_{t \rightarrow 0^+} f^{\prime}(t)$ exists (in which
case the two are equal).  
Given $\varepsilon > 0$, let $0 < \delta < t_0$ be small enough that
$\int_{(0,\delta)} |f^{\prime \prime}(t)| dt < \varepsilon$;
this is possible since $f^{\prime \prime} \in L^1(\Gamma)$.
Then for all $t_1, t_2 \in (0,\delta)$,
\begin{equation}
\label{ZhangBDVequation}
|f^{\prime}(t_2) - f^{\prime}(t_1)|
\ \le \ \int_{t_1}^{t_2} |f^{\prime \prime}(t)| dt \ < \ \varepsilon,
\end{equation}
which proves that $\lim_{t \rightarrow 0^+} f^{\prime}(t) =
d_{\vv}f(p)$ exists.


Equation (\ref{ZhangBDVequation}) also implies (in the notation of
\S\ref{MeasureSection})
that for every countable family $\cF$ of pairwise disjoint connected closed sets $E_i
\in \cA$, we have
\begin{equation*}
     \sum_{E_i \in \cF} |m_f(E_i)| \ \le \ \sum_{p \in X_f} |m_f(\{ p
     \})| + \int_\Gamma |f^{\prime \prime}(t)| dt < \infty,
\end{equation*}
so that $f \in \BDV(\Gamma)$ as desired.
\end{proof}


\vskip .1 in

For $f \in \Zh(\Gamma)$, Zhang (\cite[Appendix]{Zh1}) defined the Laplacian
$\Delta_{\Zh}$ to be the measure
\begin{equation*}
\Delta_{\Zh}(f) \ = \ - f^{\prime \prime}(x) dx -
    \sum_{p \in X_f} (\sum_{\vv \in \VEC(p)} d_{\vv}f(p)) \, \delta_p(x) \ .
\end{equation*}
We will now see that for $f \in \Zh(\Gamma)$, $\Delta(f) = \Delta_{\Zh}(f)$.  

\begin{proposition} \label{Prop3} { \ }

$A)$  If $f \in \Zh(\Gamma)$, then $\Delta_{\Zh}(f) = \Delta(f)$.

$B)$  If $f \in \BDV(\Gamma)$, and $\Delta(f)$ has the form
$\Delta(f) = g(x) dx + \sum_{p \in X} c_p \delta_{p_i}(x)$ for a
piecewise continuous function $g \in L^1(\Gamma)$ and a finite set $X$,
then $f \in \Zh(\Gamma)$.  
Furthermore, let $X_g$ be a finite set of points containing $X$,
a vertex set for $\Gamma$, and the finitely many points where $g(x)$ 
is not continuous;
put $c_p = 0$ for each $p \in X_g \backslash \{p_1, \ldots, p_m\}$.
Then $f^{\prime \prime}(x) = -g(x)$ for each $x \in \Gamma \backslash X_g$,
$f^{\prime}(x)$ is continuous on the closure of each segment of 
$\Gamma \backslash X_g$ $($interpreting $f^{\prime}(x)$ as a one-sided 
derivative at each endpoint$)$, and $\Delta_p(f) = -c_p$ for each $p \in X_g$.
\end{proposition}

\begin{proof}

A)  Suppose $f \in \Zh(\Gamma)$.  For each $p \in X_f$,
we have $\Delta(f)(\{p\}) = \Delta_{\Zh}(f)(\{p\})$.
To show that $\Delta(f) = \Delta_{\Zh}(f)$  it suffices to show that
$\Delta(f)((c,d)) = \Delta_{\Zh}(f)((c,d))$ for each open interval $(c,d)$
contained in an edge of $\Gamma \backslash X_f$.  But
\begin{eqnarray*}
\Delta(f)((c,d))
    & = & m_f((c,d)) \ = \ f^{\prime}(c) - f^{\prime}(d) \\
    & = & - \int_c^d f^{\prime \prime}(x) \, dx \ = \ \Delta_{\Zh}((c,d)) \ .
\end{eqnarray*}

B)  Now suppose $f \in \BDV(\Gamma)$, and $\Delta(f)$ has the given form.
Consider an edge in $\Gamma \backslash X_g$, and identify it with an
interval $(a,b)$ by means of the distinguished parametrization.
For each $x \in (a,b)$, we have $\Delta(f)({x}) = 0$, so $f^{\prime}(x)$ exists.
If $h > 0$ is sufficiently small, then
\begin{equation*}
f^{\prime}(x+h)-f^{\prime}(x) \ = \ -\Delta(f)([x,x+h])
                              \ = \ -\int_x^{x+h} g(t) \ dt \ ,
\end{equation*}
while if $h < 0$ then
\begin{equation*}
f^{\prime}(x+h)-f^{\prime}(x) 
         \ = \ \Delta(f)([x+h,x]) \ = \ \int_{x+h}^x g(t) dt 
         \ = \ - \int_x^{x+h} g(t) \ dt \ .
\end{equation*}
Hence
\begin{equation*}
f^{\prime \prime}(x)
    =  \lim_{h \rightarrow 0} \frac{f^{\prime}(x+h)-f^{\prime}(x)}{h} 
    =  \lim_{h \rightarrow 0} \left( -\frac{1}{h} \int_x^{x+h} g(t) \
    dt \right)
   \ = \ - g(x) \ .
\end{equation*}
The assertion that $\Delta_p(f) = -c_p$ for each $p \in X_g$ is clear.
\end{proof}

\vskip .1 in
\begin{corollary}  \label{Cor4} { \ }

$A)$ If $f \in \BDV(\Gamma)$ and $\Delta(f) = \sum_{i=1}^k c_i
\delta_{p_i}$ is a discrete measure, then $f \in \CPA(\Gamma)$.

$B)$ If $f \in \BDV(\Gamma)$, and $\Delta(f) = 0$, then $f(x) = C$ is constant.  
\end{corollary}

\begin{proof}

A) Since $\Delta(f)$ is discrete, it follows
by Proposition \ref{Prop3} that $f(x) \in \Zh(\Gamma)$
and $\Delta(f) = \Delta_{\Zh}(f)$.  Fixing a vertex set $X$ for $\Gamma$, 
we see that $f^{\prime \prime}(x) = 0$ on $\Gamma \backslash X$, so 
$f(x)$ is affine on each segment of $\Gamma \backslash X$, which means 
that $f \in \CPA(\Gamma)$.  

B)  If $\Delta(f) = 0$, then by part A),  $f \in \CPA(\Gamma)$.  Suppose $f(x)$
were not constant, and put $M = \max_{x \in \Gamma} f(x)$.  
Let $\Gamma_0 = \{ x \in \Gamma : f(x) = M \}$ and let $x_0$ be a boundary 
point of $\Gamma_0$.  Since $f \in \CPA(\Gamma)$, one sees easily that
$\Delta(f)(\{x_0\}) > 0$.  This contradicts $\Delta(f) = 0$.  
\end{proof}

\vskip .1 in

\section{The kernels $j_{\zeta}(x,y)$ and $g_{\mu}(x,y)$}
\label{KernelSection}

\vskip .1 in

In (\cite{CR}), there is a detailed study of a ``harmonic kernel'' function 
$j_{\zeta}(x,y)$ having the following properties
(see Lemma 2.10 of \cite{CR}, the Appendix to \cite{Zh1}, 
or the expository article \cite{BF}):

\vskip .05 in
A)  It is jointly continuous as a function
of three variables.  

B)  It is non-negative, with $j_{\zeta}(\zeta,y) = j_{\zeta}(x,\zeta) = 0$  
for all $x, y, \zeta \in \Gamma$.

C)  It is symmetric in $x$ and $y$:  for fixed $\zeta$, 
\begin{equation*}
j_{\zeta}(x,y) = j_{\zeta}(y,x) \ .
\end{equation*}

D)  For fixed $\zeta$ and $y$, the function $j_{\zeta,y}(x) = j_{\zeta}(x,y)$
is in $\CPA(\Gamma)$ and satisfies the Laplacian equation
\begin{equation*}
\Delta_x(j_{\zeta}(x,y)) \ = \ \delta_y(x) - \delta_{\zeta}(x) \ .
\end{equation*}


\vskip .1 in
For any real-valued, signed Borel measure $\mu$ on $\Gamma$ with $\mu(\Gamma)=1$ and $|\mu|(\Gamma) < \infty$,
put
\begin{equation*}
j_{\mu}(x,y) \ = \ \int_{\Gamma} j_{\zeta}(x,y) \, d\mu({\zeta}) \ .
\end{equation*}
Clearly $j_{\mu}(x,y)$ is symmetric, and is jointly continuous
in $x$ and $y$.  It has the following properties
which are less obvious :

\begin{proposition} \label{Prop5} { \ }

$A)$  There is a constant $c_{\mu}(\Gamma)$ such that for each $y \in \Gamma$
\begin{equation*}
\int_{\Gamma} j_{\mu}(x,y) \, d\mu(x) \ = \ c_{\mu}(\Gamma) \ .
\end{equation*}

$B)$  For each $y \in \Gamma$, the function $F_y(x) = j_{\mu}(x,y)$
belongs to the space $\BDV(\Gamma)$, and satisfies
\begin{equation*}
\Delta_x(F_y) \ = \ \delta_y(x) - \mu \ .
\end{equation*}
\end{proposition}

\begin{proof} For A) see \cite{CR}, Lemma 2.16, p.21.

For B), fix $y$ and note that the existence of the directional derivatives
$d_{\vv} F_y(x)$  follows from (\cite{CR}, Lemma 2.13, p.19).
Hence $F_y \in \cD(\Gamma)$.  Also, noting that $j_y(y,x) = 0$ for all $x$,
it follows that for the measure $\omega = \mu - \delta_y(x)$ of
total mass $0$, we have $j_{\omega}(x,y) = j_{\mu}(x,y)$ for all $x$.
Then, (\cite{CR}, Lemma 2.14, p.20) tells us that for any subset
$D \subset \Gamma$ which is a finite union of closed intervals,
the finitely additive set function $m_{F_y}$ satisfies
\begin{equation*}
m_{F_y}(D) \ = \ - \omega(D) \ = \ (\delta_y(x) - \mu)(D)
\end{equation*}
(again recall that  $\Delta = -\Delta_{\CR}$).  For any countable
collection $\cF$ of pairwise disjoint sets $D_i$ of the type above,
it follows that
\begin{equation*}
\sum_{i=1}^{\infty} |m_{F_y}(D_i)| \ \le \ \sum_{i=1}^{\infty} |\omega|(D_i)
\ \le \ |\omega|(\Gamma) \ ,
\end{equation*}
so by Proposition \ref{Prop1}, $F_y \in \BDV(\Gamma)$.  The measure
$\Delta(F_y) = m_{F_y}^*$ attached to $F_y$ is determined by its
value on closed intervals,
so it coincides with $\delta_y(x) - \mu$.
\end{proof}

\vskip .1 in
We now define a kernel $g_{\mu}(x,y)$ which will play a key role in
the rest of the paper:
\begin{equation*}
g_{\mu}(x,y) \ = \ j_{\mu}(x,y) - c_{\mu}(\Gamma) \ .
\end{equation*}
It follows from Proposition~\ref{Prop5} that 
$g_{\mu}(x,y)$ continuous, symmetric, and for each $y$
\begin{equation} \label{xf1}
\int_{\Gamma} g_{\mu}(x,y) \, d\mu(x) \ = \ 0 \ .
\end{equation}

\section{The operator $\varphi_{\mu}$}

\vskip .1 in

{\bf Convention.}
For the rest of this paper, $\mu$ will denote a real, signed Borel
measure with $\mu(\Gamma) = 1$ and $|\mu|(\Gamma) < \infty$.

\vskip .1 in

For any complex Borel measure $\nu$ on $\Gamma$ with $|\nu|(\Gamma) < \infty$,
define the integral transform
\begin{equation*}
\varphi_{\mu}(\nu)(x) \ = \ \int_{\Gamma} g_{\mu}(x,y) \, d\nu(y) \ .
\end{equation*}
Write 
\begin{equation*}
\BDV_{\mu}(\Gamma) 
    = \{ f \in \BDV(\Gamma) : \int_{\Gamma} f(x) \, d\mu(x) = 0 \} \ .  
\end{equation*}    
    
\begin{proposition} \label{Prop6}
The function $F(x) = \varphi_{\mu}(\nu)(x)$ belongs to $\BDV_{\mu}(\Gamma)$,
and satisfies
\begin{equation*}
\Delta(F) \ = \ \nu - \nu(\Gamma) \, \mu \ .
\end{equation*}
\end{proposition}

\begin{proof}
First note the following identity:
for each $p \in \Gamma$,
\begin{equation*}
j_{\zeta}(x,y) \ = \ j_{\zeta}(x,p) - j_y(x,p) + j_y(\zeta,p) \ .
\end{equation*}
To see this, fix $\zeta$, $y$ and $p$; put $H(x) = j_{\zeta}(x,y)$,
$h(x) = j_{\zeta}(x,p) - j_y(x,p)$.  Then
\begin{equation*}
\Delta_x(H) \ = \ \Delta_x(h) \ = \ \delta_y(x) - \delta_{\zeta}(x) \ ,
\end{equation*}
so by (\cite{CR}, Lemma 2.6, p.11) $H(x)-h(x)$ is a constant $C$.
To evaluate $C$, take $x = \zeta$;
using $j_{\zeta}(\zeta,y) = j_{\zeta}(\zeta,p) = 0$
we obtain $C = j_y(\zeta,p)$.

Inserting this into the definition of $\varphi_{\mu}(\nu)$, we find
\begin{eqnarray*}
F(x) & = & \int_{\Gamma} g_{\mu}(x,y) \, d\nu(y) \\
     & = & \int_{\Gamma} \left( \int_{\Gamma}
               \left\{ j_{\zeta}(x,p) - j_y(x,p) + j_y(\zeta,p)
     \right\} \, d\mu(\zeta)
                       - c_{\mu}(\Gamma) \right) \, d\nu(y) \\
     & = & \nu(\Gamma) j_{\mu}(x,p) - j_{\nu}(x,p) + C_p,
\end{eqnarray*}
where $C_p$ is a constant.  Hence Proposition~\ref{Prop5}(B) shows
that $F \in \BDV(\Gamma)$ and $\Delta(F) = \nu - \nu(\Gamma) \mu$.

Finally, the fact that $F \in \BDV_\mu(\Gamma)$ follows from 
Fubini's theorem:
\begin{equation*}
\int_{\Gamma} F(x) \, d\mu(x)
     \ = \ \int_{\Gamma} 
             \Big( \int_{\Gamma} g_{\mu}(x,y) \, d\mu(x) \Big) \, d\nu(y)
     \ = \ 0 \ .
\end{equation*}
\end{proof}

\section{Eigenfunctions of $\varphi_{\mu}$ in $L^2(\Gamma)$}
\label{IntegralOperatorEigenfunctions}

\vskip .1 in

Given a Borel measurable function $f: \Gamma \rightarrow \CC$,
write $\|f\|_1$, $\|f\|_2$, and $\|f\|_{\infty}$ for
its $L^1$, $L^2$ and $\sup$ norms.
Write $L^1(\Gamma)$, $L^2(\Gamma)$
for the spaces of $L^1$, $L^2$ functions on $\Gamma$ relative to the
measure $dx$.  
It follows from the Cauchy-Schwarz inequality that $L^2(\Gamma) \subset L^1(\Gamma)$.
Write $\<f,g\> = \int_{\Gamma} f(x) \overline{g(x)} \, dx$
for the inner product on $L^2(\Gamma)$.


\vskip .1 in

Given $f \in L^2(\Gamma)$, let $\nu$ be the bounded Borel measure $f(x) dx$, 
and define
\begin{equation*}
\varphi_{\mu}(f) \ = \ \varphi_{\mu}(\nu)
                 \ = \ \int_{\Gamma} g_{\mu}(x,y) f(y) \, dy \ .
\end{equation*}
By Proposition \ref{Prop6}, $\varphi_{\mu}(f) \in \BDV_{\mu}(\Gamma)$. 

Equip $\BDV_{\mu}(\Gamma)$ with the $L^2$ norm, and view it as a subspace
of $L^2(\Gamma)$.  Write $L^2_{\mu}(\Gamma)$ for the $L^2$ completion
of $\BDV_{\mu}(\Gamma)$.   

\vskip .1 in

For the next proposition, recall that a sequence $\nu_n$ of measures on a
compact space $X$ {\em converges weakly} to a measure $\nu$ if $\int f
d\nu_n \to \int f d\nu$ for all continuous function $f : X \to \RR$.

\begin{proposition} \label{Prop7}
Suppose $\nu_1, \nu_2, \nu_3, \ldots$ is a bounded sequence of 
complex Borel measures on $\Gamma$ which converge weakly to a 
Borel measure $\nu$.  Then the functions $\varphi_{\mu}(\nu_i)$ converge
uniformly to $\varphi_{\mu}(\nu)$ on $\Gamma$.
\end{proposition}

\begin{proof}
First, note that $\nu$, being a weak limit of bounded measures, 
must be bounded.
Write $F_i(x) = \varphi_{\mu}(\nu_i)$ and $F(x) = \varphi_{\mu}(\nu)$.
For each $x$, $g_{\mu}(x,y)$ is a continuous function
of $y$.  Hence, by the definition of weak convergence,
the functions $F_i(x)$ converge pointwise to $F(x)$. 

We claim that the convergence is uniform.
Fix $\varepsilon > 0$, and let $d(x,y)$ be the metric on the graph $\Gamma$.
Since $g_{\mu}(x,y)$ is continuous and $\Gamma$ is compact,
$g_{\mu}(x,y)$ is uniformly continuous.
Hence there is a $\delta > 0$ such that for any $x_1, x_2 \in \Gamma$
satisfying $d(x_1,x_2) < \delta$,
\begin{equation*}
|g_{\mu}(x_1,y) - g_{\mu}(x_2,y)| \ < \ \varepsilon
\end{equation*}
for all $y \in \Gamma$.
Write $B_{x}(\delta) = \{ z \in \Gamma : d(x,z) < \delta \}$.
Since $\Gamma$ is compact, it can be covered by finitely many balls
$B(x_j,\delta)$, $j = 1, \ldots, m$.  Let $N$ be large enough so that
$|F_i(x_j)-F(x_j)| < \varepsilon$ for all $j = 1, \ldots, m$
and all $i \ge N$.  Then for each $x \in \Gamma$, and each $i \ge N$,
if $x_j$ is such that $d(x,x_j) < \delta$ then
\begin{equation*}
|F_i(x)-F(x)| \ \le \ |F_i(x)-F_i(x_j)| + |F_i(x_j)-F(x_j)|+|F(x_j)-F(x)| \ .
\end{equation*}
But $|F_i(x)-F_i(x_j)|
   \le \int_{\Gamma} |g_{\mu}(x,y)-g_{\mu}(x_j,y)| \, d|\nu_i|(y)
   \le |\nu_i|(\Gamma) \varepsilon$,
and similarly 
$|F(x)-F(x_j)| \le |\nu|(\Gamma) \varepsilon$.
Hence $|F_i(x)-F(x)| < (\sup_i |\nu_i|(\Gamma) + |\nu|(\Gamma)+1) \varepsilon$.
\end{proof}

\begin{proposition} \label{Prop8}  The operator  
$\varphi_{\mu} : L^2(\Gamma) \rightarrow L^2(\Gamma)$
is compact and self-adjoint, and maps $L^2(\Gamma)$ into
$\BDV_{\mu}(\Gamma)$.

\end{proposition}

\begin{proof}
Since $\mu$ is real, $g_{\mu}(x,y)$ is real-valued and symmetric, so 
$\varphi_{\mu}$ is self-adjoint.

It remains to see that $\varphi_{\mu}$ is compact.
Let $f_1, f_2, \ldots \in L^2(\Gamma)$ be a sequence of functions
with bounded $L^2$ norm.  By the Cauchy-Schwarz inequality, the $L^1$ norms of the $f_i$
are also bounded, so the sequence of measures $f_i(x) \, dx$ is
bounded and has a subsequence
converging weakly to a bounded measure $\nu$ on $\Gamma$.  After passing
to this subsequence we can assume it is the entire sequence $\{ f_i \}$.
By Proposition \ref{Prop7}, the functions $\varphi_{\mu}(f_i)$ converge
uniformly to $\varphi_{\mu}(\nu)$, which belongs to $\BDV_{\mu}(\Gamma)$.
Clearly they converge to $\varphi_{\mu}(\nu)$ under the $L^2$ norm as well.
\end{proof}

\vskip .1 in 

\begin{definition} \label{DefVarphiEV} 
A nonzero function $f \in L^2(\Gamma)$ is an eigenfunction for $\varphi_{\mu}$, 
with eigenvalue $\alpha$, if $\varphi_{\mu}(f) = \alpha \cdot f$.
\end{definition}

Applying well-known results from spectral theory, we have

\begin{theorem} \label{Thm9}
The map $\varphi_{\mu}$ has countably many eigenvalues $\alpha_i$,
each of which is real and occurs with finite multiplicity.  
The nonzero eigenvalues can be ordered so that 
$|\alpha_1| \ge |\alpha_2| \ge \ldots$ and
$\lim_{i \rightarrow \infty} |\alpha_i| = 0$.
Each eigenfunction corresponding to a nonzero eigenvalue belongs
to $\BDV_{\mu}(\Gamma)$.   Finally, the space $L^2(\Gamma)$ has an orthonormal
basis consisting of eigenfunctions of $\varphi_{\mu}$.
\end{theorem}

\begin{proof}  See Lang (\cite{Lang}):  all the assertions are contained in the
Spectral Theorem for compact hermitian operators (\cite{Lang}, p.165), 
together with (\cite{Lang}, Theorem 10, p.208), 
and (\cite{Lang}, Theorem 11, p.210).  
\end{proof}

\vskip .1 in
We now investigate the eigenfunctions of $\varphi_{\mu}$ with eigenvalue $0$, 
that is, the kernel of $\varphi_{\mu}$ in $L^2(\Gamma)$.  We begin
by reformulating Proposition \ref{Prop6} in the present context:

\begin{lemma} \label{Lem10}  For each $f \in L^2(\Gamma)$, 
\begin{equation*}
\Delta(\varphi_{\mu}(f)) \ = \ f(x) \, dx - (\int_{\Gamma} f(x) dx) \, \mu \ .
\end{equation*}
\end{lemma}

\begin{proposition} \label{Prop11}
If $\mu = g(x) dx$ for a real-valued function $g(x) \in L^2(\Gamma)$
with $\int_{\Gamma} g(x) dx = 1$, then $\Ker(\varphi_{\mu}) = \CC \cdot g(x)$ 
and the $L^2$-closure of $\varphi_{\mu}(\BDV_{\mu}(\Gamma))$ 
$($and of $\varphi_{\mu}(L^2(\Gamma))$$)$ is $\CC \cdot g(x)^{\perp}$.  
Otherwise, $\Ker(\varphi_{\mu}) = \{0\}$ and the $L^2$-closure of 
the image is all of $L^2(\Gamma)$.  

In either case, $\Ker(\varphi_{\mu}) \cap \BDV_{\mu}(\Gamma) = \{0\}$.  
\end{proposition}

\begin{proof}
Suppose $0 \ne f \in L^2(\Gamma)$ belongs to $\Ker(\varphi_{\mu})$.  
By Lemma \ref{Lem10}, 
\begin{equation*}
0 \ = \ \Delta_x(\varphi_{\mu}(f)) 
  \ = \ f(x) dx - (\int_{\Gamma} f(x) dx) \mu \ . 
\end{equation*}
Write $C_f = \int_{\Gamma} f(x) dx$.  Since $f \ne 0$, the equation above
shows that $C_f \ne 0$.  Hence $\mu = (1/C_f) f(x) dx = g(x) dx$, 
with $g \in L^2(\Gamma)$ and $\int_{\Gamma} g(x) dx = 1$, 
and $f \in \CC \cdot g(x)$.  Conversely if
$\mu$ has this form, then  
\begin{eqnarray*}
\varphi_{\mu}(g)(x) \ = \ \int_{\Gamma} g_{\mu}(x,y) g(y) dy  
\ = \ \int_{\Gamma} g_{\mu}(x,y) \, d\mu(y)  \ = \ 0  
\end{eqnarray*}
by formula (\ref{xf1}) and the symmetry of $g_{\mu}(x,y)$. 

When $\mu = g(x) dx$, if $g(x) \notin \BDV(\Gamma)$ then clearly
$\Ker(\varphi_{\mu}) \cap \BDV_{\mu}(\Gamma) = \{0\}$.  If 
$g(x) \in \BDV(\Gamma)$, then since $g(x)$ is real-valued, 
$\int_{\Gamma} g(x) d\mu(x) = \int_{\Gamma} g(x)^2 dx > 0$, so 
$g(x) \notin \BDV_{\mu}(\Gamma)$. 

The assertions about the $L^2$-closure of the image of $\varphi_{\mu}$ follow
from our description of the kernel and from Theorem \ref{Thm9}.
\end{proof}



\vskip .1 in

\section{Eigenfunctions from the point of view of Laplacians}
\label{EFLaplacianSection} 

In this section we will give an integro-differential
characterization of eigenfunctions of $\varphi_{\mu}$.

\begin{proposition} \label{Prop32}
A function $0 \ne f \in \BDV(\Gamma)$ is an eigenfunction of $\varphi_{\mu}$
with nonzero eigenvalue if and only if

$A)$ $\int_{\Gamma} \, f(x) d\mu(x) = 0$;  and

$B)$ for some constants $\lambda, C \in \CC$, 
\quad $\Delta(f) = \lambda \cdot (f(x) dx - C \mu)$.

\noindent{If $A)$ and $B)$ hold,} 
then necessarily $\lambda \in \RR$ and $C = \int_{\Gamma} f(x) \, dx$, 
with $\varphi_{\mu}(f) = \frac{1}{\lambda} \cdot f$.   
\end{proposition}

\noindent{\bf Remark:}  In Theorem \ref{Thm16} below, we will see that in 
fact $\lambda > 0$.  
\vskip .1 in

\begin{proof}
If $0 \ne f \in \BDV(\Gamma)$ is an eigenfunction of $\varphi_{\mu}$ with
nonzero eigenvalue $\alpha$, then $\varphi_{\mu}(f) = \alpha  f$, so
$0 = \int \varphi_{\mu}(f)(x) \, d\mu(x) 
= \alpha \cdot \int_{\Gamma} f(x) \, d\mu(x)$ so $f \in \BDV_{\mu}(\Gamma)$.
By Lemma \ref{Lem10}
\begin{equation*}
\Delta(f) \ = \ \frac{1}{\alpha} \cdot \Delta(\varphi_{\mu}(f)) \ = \ 
\frac{1}{\alpha} \cdot (f(x) dx - C_f \mu)
\end{equation*}
where $C_f = \int_{\Gamma} f(x) dx$.

Conversely, if A) holds, then $f \in \BDV_{\mu}(\Gamma)$;  and if B)
also holds then by Lemma \ref{Lem10} 
\begin{equation} \label{FNV1}
\Delta(\lambda \cdot \varphi_{\mu}(f) - f) 
\ = \ \lambda \cdot(\int_{\Gamma} f(x) dx - C) \cdot \mu \ . 
\end{equation}
By Theorem \ref{Thm2}, 
\begin{eqnarray}
0  & = & \int_{\Gamma} d\Delta(\lambda \varphi_{\mu}(f) - f) \label{FNV2}\\
   & = & \lambda \cdot (\int_{\Gamma} f(x) dx - C) \cdot \int_{\Gamma} d\mu  
   \ = \ \lambda \cdot (\int_{\Gamma} f(x) dx - C) \ . \notag
\end{eqnarray} 
By (\ref{FNV1}) we have $\Delta( \lambda \varphi_{\mu}(f) - f)) = 0$, 
so Corollary \ref{Cor4} shows $\lambda \varphi_{\mu}(f) - f$ 
is a constant function. 
Since $f$ and $\varphi_{\mu}(f)$ both belong to $\BDV_{\mu}(\Gamma)$, 
we conclude $\lambda \varphi_{\mu}(f) = f$.  Here $f \ne 0$ by hypothesis,
so $\lambda \ne 0$.  It follows that $f$ is an eigenvalue of $\varphi_{\mu}$
with nonzero eigenvalue $\alpha = 1/\lambda$.  By Theorem \ref{Thm9},
$\alpha$, and hence $\lambda$, is real.  Since $\lambda \ne 0$, equation
(\ref{FNV2}) gives $C = \int_{\Gamma} f(x) dx$.
\end{proof}

\begin{definition} \label{EFDef}
A nonzero function $f \in \BDV(\Gamma)$ will be called an eigenfunction
of $\Delta$ in $\BDV_{\mu}(\Gamma)$ if it satisfies conditions (A) and (B)
in Proposition \ref{Prop32} 
\end{definition}  

\vskip .1 in

\section{Positivity and the Dirichlet Semi-norm}

\vskip .1 in

In this section we will show that the nonzero eigenvalues of $\varphi_{\mu}$ 
are positive.  This is intimately connected with the
positivity of the Dirichlet semi-norm on $\BDV_{\mu}(\Gamma)$, 
which we also investigate.


\vskip .1 in
 On the space $\Zh(\Gamma)$ the Dirichlet inner product $\<f,g\>_{\Dir}$
is given by the equivalent formulas 
\begin{equation} \label{FBN1}
\<f,g\>_{\Dir} \ = \ \int_{\Gamma} f'(x) \overline{g'(x)} \, dx  
              \ = \ \int_{\Gamma} f(x) \, d\overline{\Delta(g)}
               \ = \ \int_{\Gamma} \overline{g(x)} \, d\Delta(f) \ .
\end{equation}
We have seen that $\Zh(\Gamma) \subset \BDV(\Gamma)$.  We claim that for 
all $f, g \in \BDV(\Gamma)$  
\begin{equation} \label{FNA1}
\int f(x) \, d\overline{\Delta(g)}
               \ = \ \int \overline{g(x)} \, d\Delta(f).
\end{equation}    
            
To see this, take $f,g \in \BDV(\Gamma)$ and write 
$\nu = \Delta(f)$ and $\omega = \Delta(g)$.  
Then $\nu(\Gamma) = \omega(\Gamma) = 0$. 
For any real Borel measure $\mu$ on $\Gamma$, the function  
$f_{\mu}(x) := \varphi_{\mu}(\nu) = \ \int_{\Gamma} g_{\mu}(x,y) \, d\nu(y)$
belongs to $\BDV_{\mu}(\Gamma)$.  
By Proposition \ref{Prop6}, $\Delta(f_{\mu}) = \nu$.    
Hence $\Delta(f-f_{\mu}) = 0$, which means that $f = f_{\mu}+C$
for some  constant $C$.  Thus for any $g \in \BDV(\Gamma)$,
\[
\int_{\Gamma} f(y) \, \overline{d\Delta(g)(y)} 
     =  \int_{\Gamma} f_{\mu}(y) \, \overline{d\omega(y)} 
     =  \iint_{\Gamma \times \Gamma} g_{\mu}(x,y) 
                      \, d\nu(x) \, \overline{d\omega(y)}.
\]

In light of this we define the Dirichlet pairing 
on $\BDV(\Gamma) \times \BDV(\Gamma)$ by 
\begin{eqnarray}
\<f,g\>_{\Dir} & = & \int_{\Gamma} f(y) \, \overline{d\Delta(g)(y)} 
        \ = \ \int_{\Gamma} \overline{g(x)} \, d\Delta(f)(x) \label{fgz2} \\
        & = & \iint_{\Gamma \times \Gamma} g_{\mu}(x,y) \, d\Delta(f)(x) 
                           \, \overline{d\Delta(g)(y)} \ . \label{fgz3}
\end{eqnarray}   

Note that 
if $f(x)$ and $g(x)$ are replaced by $f(x)+C$ and $g(x)+D$ for constants
$C$ and $D$, the value of $\<f,g\>_{\Dir}$ remains unchanged, since 
$\Delta(C)= \Delta(D) = 0$. 
 Finally, observe that the choice of $\mu$ in (\ref{fgz3}) 
is arbitrary, since it does not occur in (\ref{fgz2}).    
In particular (taking $\mu = \delta_{\zeta}(x)$), 
the kernel $g_{\mu}(x,y)$ in (\ref{fgz3})
can be replaced by $j_{\zeta}(x,y)$ for any point $\zeta$.  

\vskip .1 in
We will now show that $\<f,g\>_{\Dir}$ defines a semi-norm on $\BDV(\Gamma)$
whose kernel consists of precisely the constant functions.  For this, we will
need some preliminary lemmas.  The first is well-known:  

\begin{lemma} \label{LemProductConvergence}
Let $X,Y$ be compact Hausdorff topological spaces, and let $\mu_i$
(resp. $\nu_i$) be a bounded sequence of signed Borel measures on $X$
(resp. $Y$).  If $\mu_i$ converges weakly to $\mu$ on $X$ and $\nu_i$
converges weakly to $\nu$ on $Y$, then $\mu_i \times \nu_i$ converges
weakly to $\mu \times \nu$ on $X \times Y$.
\end{lemma}

\begin{proof}
By hypothesis, there exist $M,N>0$ such that $|\mu_i|(\Gamma) \leq M$
and $|\nu_i(\Gamma)| \leq N$ for all $i$.
According to the Stone-Weierstrass theorem, the set 
\begin{equation*}
S = \{ \sum c_i f_i(x) g_i(y) \; : \; c_i \in \RR, f_i \in \cC(X), g_i
\in \cC(Y) \}
\end{equation*}
is dense in $\cC(X \times Y)$, and it follows easily from the
hypotheses of weak convergence that 
\begin{equation*}
\int_{X \times Y} F \, d(\mu_i \times \nu_i) \to \int_{X \times Y}
F \, d(\mu \times \nu)
\end{equation*}
for all $F \in S$.  
Now suppose $F \in \cC(X \times Y)$ is arbitrary, and let $\varepsilon
> 0$.  Since $S$ is dense, we may choose $H(x,y) \in S$ so that
the difference $G(x,y) = F(x,y) - H(x,y)$ satisfies $\| G \|_\infty <
\varepsilon$.
Since $|\mu \times \nu| \leq |\mu| \times |\nu|$ (which follows from the
Hahn decomposition theorem), we obtain
\begin{equation*}
|\int G \, d(\mu_i \times \nu_i) - \int G \, d(\mu \times \nu)| \leq
2MN \| G \|_\infty < 2MN\cdot\varepsilon .
\end{equation*}
This gives what we want.
\end{proof}

\begin{lemma} \label{Lem14}
Suppose $f, g \in \BDV(\Gamma)$ and that $f_1, f_2, \ldots$ and
$g_1, g_2, \ldots$ are sequences of functions in $\BDV(\Gamma)$
such that

$A)$  $|\Delta(f)|(\Gamma)$, $|\Delta(g)|(\Gamma)$,
and the sequences $|\Delta(f_i)|(\Gamma)$ and $|\Delta(g_j)|(\Gamma)$
are bounded;

$B)$  $\Delta(f_1), \Delta(f_2), \ldots$ converges weakly to $\Delta(f)$;

$C)$  $\Delta(g_1), \Delta(g_2), \ldots$ converges weakly to $\Delta(g)$.

\noindent{Then} 
$\lim_{i \rightarrow \infty} \<f_i,g_i\>_{\Dir} = \<f,g\>_{\Dir}$.
\end{lemma}

\begin{proof}  
Write $\nu = \Delta(f)$, $\omega = \Delta(g)$,
$\nu_i = \Delta(f_i)$, and $\omega_i = \Delta(g_i)$.  
By Lemma~\ref{LemProductConvergence},
the sequence of measures $\nu_i \times \omega_i$ converges 
weakly to $\nu \times \omega$ on $\Gamma \times \Gamma$.

Hence
\begin{eqnarray*}
\lim_{i \rightarrow \infty} \<f_i,g_i\>_{\Dir}
   & = & \lim_{i \rightarrow \infty}
               \iint_{\Gamma \times \Gamma} g_{\mu}(x,y)
                          \, d\nu_i(x) \overline{d\omega_i(y)} \\
   & = &   \iint_{\Gamma \times \Gamma} g_{\mu}(x,y)
                          \, d\nu(x) \overline{d\omega(y)}
   \ = \ \<f,g\>_{\Dir}.
\end{eqnarray*}
\end{proof}

We now come to the main result of this section.

\begin{theorem} \label{Thm16} { \ }

$A)$  For each $f \in \BDV(\Gamma)$ the Dirichlet pairing satisfies
         $\<f,f\>_{\Dir} \ge 0$, with $\<f,f\>_{\Dir} = 0$ if and only if
         $f(x) = C$ is constant.  Thus $\<f,g\>_{\Dir}$ is a semi-norm
         on $\BDV(\Gamma)$ with kernel $\CC \cdot 1$, and
         its restriction to $\BDV_{\mu}(\Gamma)$ is a norm.

$B)$  The nonzero eigenvalues $\alpha$ of
         $\varphi_{\mu}$ are positive, and
         for each $f \in L^2(\Gamma)$, 
\begin{equation*}
         \<\varphi_{\mu}(f),f\> \ \ge \ 0 \ .
\end{equation*}
\end{theorem}

\begin{proof}
Let $f \in \BDV(\Gamma)$ be arbitrary, and put $\nu = \Delta(f)$.
It is easy to construct a sequence of discrete measures $\nu_1, \nu_2, \ldots$
which converge weakly to $\nu$,
with $\nu_i(\Gamma) = 0$ and $|\nu_i|(\Gamma) \le |\nu|(\Gamma)$ for each $i$.
(Choose a sequence of numbers $\eta_i \rightarrow 0$.  For each $i$
subdivide $\Gamma$ into finitely many segments $e_{ij}$ with length at most
$\eta_i$.  For each $i,j$ choose a point $p_{ij} \in e_{ij}$, and put
$\nu_i = \sum_j \nu(e_{ij}) \delta_{p_{ij}}(x)$.)

Put $f_i = \varphi_{\mu}(\nu_i)$ for each $i$.
Then $\Delta(f_i) = \nu_i$ and $f_i \in \CPA(\Gamma) \subset \Zh(\Gamma)$.
It follows that  $\<f_i,f_i\>_{\Dir} = \int |f_i'(x)|^2 \, dx \ge 0$ 
for each $i$. By Lemma \ref{Lem14},
\begin{equation*}
\<f,f\>_{\Dir}
\ = \ \lim_{i \rightarrow \infty} \<f_i,f_i\>_{\Dir} \ \ge 0 \ .
\end{equation*}

Thus, at least $\<f,f\>_{\Dir}$ is non-negative.
We will use this to show that the eigenvalues of
$\Delta$ in $\BDV_{\mu}(\Gamma)$ are positive.
Suppose $0 \ne f \in \BDV_{\mu}(\Gamma)$ is an eigenfunction of $\varphi_\mu$
with eigenvalue $\alpha$.  By Proposition \ref{Prop11}, $\alpha \ne 0$.  

Since $\varphi_{\mu}(f) = \alpha \cdot f$, 
Lemma \ref{Lem10} shows that
\begin{equation*}
\alpha \cdot \Delta(f) \ = \ \Delta(\varphi_{\mu}(f)) 
          \ = \  f(x) dx - C_f \cdot \mu 
\end{equation*}
where $C_f = \int f(x) \, dx$.  Recalling that $\alpha$ is real, 
and using that $\int f(x) \, d\mu(x) = 0$, we have 
\begin{eqnarray}
0 \ \le \ \<f,f\>_{\Dir}
      & = & \int_{\Gamma} f(x) \, \overline{d\Delta(f)(x)} \\
      & = & \int_{\Gamma} f(x) \, \frac{1}{\alpha} \overline{f(x)} dx
               \ = \ \frac{1}{\alpha} \<f,f\>_{L^2} \ .
\end{eqnarray}
Since $\<f,f\>_{L^2} > 0$, we must have $\alpha > 0$.  

It follows from this that $\<\varphi_{\mu}(f),f\>_{L^2} \ge 0$ 
for all $f \in L^2(\Gamma)$.  

We can now show the positivity of $\<f,f\>_{\Dir}$.

Since the Dirichlet pairing is positive semi-definite on $\BDV(\Gamma)$, 
it follows from the Cauchy-Schwarz inequality that if $\<f,f\>_{\Dir}=0$ 
then $\< f,g \>_{\Dir} = 0$ for all $g \in \BDV(\Gamma)$, i.e.,
\[
\int_\Gamma \overline{g(x)} \, d\Delta(f)(x) = 0
\]
for all $g \in \BDV(\Gamma)$.  
But then $\Delta(f)=0$, so that $f$ is constant by Corollary~\ref{Cor4}.
In particular, $\<f,g\>_{\Dir}$ is positive definite on $\BDV_{\mu}(\Gamma)$. 
\end{proof}

\vskip .1 in
\noindent{\bf Application:  Positivity of the energy pairing.}
\vskip .1 in

Let $\Meas(\Gamma)$ be the space of  bounded Borel measures on $\Gamma$,
and let $\Meas_0(\Gamma)$ be the subspace consisting of measures 
with $\nu(\Gamma) = 0$.
For $\nu, \omega \in \Meas(\Gamma)$, define the `$\mu$-energy pairing'
\begin{equation*}
\<\nu,\omega\>_{\mu} 
   \ = \ \iint_{\Gamma \times \Gamma} g_{\mu}(x,y) \, 
                   d\nu(x) \overline{d\omega(y)} \ . 
\end{equation*}
Note that if $\nu,\omega \in \Meas_0(\Gamma)$, then the energy pairing  
\[
\< \nu,\omega \> \ := \ \< \nu,\omega \>_{\mu}
\]
is independent of $\mu$.
Indeed, let $f, g \in \BDV(\Gamma)$ be functions
with $\nu = \Delta(f)$, $\omega = \Delta(g)$;    
such functions exist by Proposition \ref{Prop6}.  
By (\ref{fgz3}), $\< \nu,\omega \>_{\mu} = \< f,g \>_{\Dir}$,
which is independent of $\mu$.   

\vskip .1 in

\begin{theorem} \label{Thm17} { \ }

$A)$  The energy pairing $\<\nu,\omega\>$ 
is positive definite on $\Meas_0(\Gamma)$.

$B)$  For each $\mu$, the $\mu$-energy pairing $\<\nu,\omega\>_{\mu}$ 
is positive semi-definite on $\Meas(\Gamma)$, with kernel $\CC \mu$.

$C)$ Among all  $\nu \in \Meas(\Gamma)$ with $\nu(\Gamma) = 1$,  
$\mu$ is the unique measure which minimizes the energy integral
\begin{equation*}
I_{\mu}(\nu) \ = \  \iint_{\Gamma \times \Gamma} g_{\mu}(x,y) \, 
                   d\nu(x) \overline{d\nu(y)} \ .
\end{equation*}
\end{theorem}

\begin{proof}  Fix $\mu$.  
Given $\nu, \omega \in \Meas_0(\Gamma)$, 
put $f = \varphi_{\mu}(\nu)$ and $g = \varphi_{\mu}(\omega)$.  As noted above,  
$\<\nu,\omega\> = \<\nu,\omega\>_{\mu} = \<f,g\>_{\Dir}$.  
By Theorem \ref{Thm16}, $\<\nu,\omega\>$ is positive definite.  

Since $\int_{\Gamma} g_{\mu}(x,y) \, d\mu(y) = 0$ for each $x$, clearly
$I_{\mu}(\mu) = 0$ and $\<\omega,\mu\>_{\mu} = \<\mu,\omega\> = 0$ 
for each $\omega \in \Meas(\Gamma)$.  
Since $\Meas(\Gamma) = \Meas_0(\Gamma) + \CC \mu$, it follows that 
$\<\nu,\omega\>_{\mu}$ is positive semidefinite on $\Meas(\Gamma)$ 
with kernel $\CC \mu$.   

Now let $\nu$ be any measure with $\nu(\Gamma) = 1$.
Then $\nu-\mu \in \Meas_0(\Gamma)$, so
\begin{eqnarray*}   
I_{\mu}(\nu) & = & \<\nu,\nu\>_{\mu}  
             \ = \ \<\mu + (\nu-\mu),\mu + (\nu-\mu)\>_{\mu} \\ 
             & = & I_{\mu}(\mu) + \<\nu-\mu,\nu-\mu\>_{\mu} 
             \ \ge \ I_{\mu}(\mu) \ = \ 0 \ ,
\end{eqnarray*}
with equality if and only if $\nu=\mu$.
\end{proof} 

\vskip .1 in

\section{The space $\Dir_{\mu}(\Gamma)$ revisited}

\vskip .1 in
In this section we will use the fact that $\|f\|_{\Dir}$ is a norm on
$\BDV_{\mu}(\Gamma)$ to show that $\BDV_{\mu}(\Gamma)$ is a subspace of
$\Dir_{\mu}(\Gamma)$.  This is a key step towards relating the notions of 
eigenfunctions in Theorem \ref{MainPairingTheorem} and Theorem \ref{Thm9}.

\vskip .1 in
Recall that $\Zh(\Gamma)$ is the space of continuous functions 
$f : \Gamma \rightarrow \CC$ which are piecewise $\cC^2$, 
with $f''(x) \in L^1(\Gamma)$.
Let $V(\Gamma)$ be the space of continuous functions 
$f : \Gamma \rightarrow \CC$ 
such that $f$ is piecewise $\cC^1$ and $f^{\prime} \in L^2(\Gamma)$.
It easy check that $\Zh(\Gamma) \subset V(\Gamma)$.    
For $f, g \in V(\Gamma)$, we define 
\begin{equation} \label{FMX1}
\<f,g\>_{\Dir} = \int_{\Gamma} f^{\prime}(x) \overline{g^{\prime}(x)} \, dx.
\end{equation} 
It is easy to see that for $f \in V(\Gamma)$, 
$\<f,f\>_{\Dir} = 0$ if and only if $f(x)$ is a constant function.  

Given a real Borel measure $\mu$ of total mass $1$, 
put $V_{\mu}(\Gamma) = \{f \in V(\Gamma) : \int f(x) \, d\mu(x) = 0 \}$.
Then the restriction of $\<f,g\>_{\Dir}$ to $V_{\mu}(\Gamma)$ 
is an inner product.   

The space $\Zh(\Gamma)$ is contained in $\BDV(\Gamma)$ and $V(\Gamma)$, 
and both of these spaces are contained in $\cC(\Gamma)$, 
though neither is contained in the other.  Putting this another way, 
there is a natural extension of $\<f,g\>_{\Dir}$ on $\Zh(\Gamma)$ 
to each of $\BDV(\Gamma)$ and $V(\Gamma)$, and their restrictions to 
$\BDV_{\mu}(\Gamma)$ and $V_{\mu}(\Gamma)$ induce norms which coincide
on $\Zh_{\mu}(\Gamma)$.  Note that $\CPA_{\mu}(\Gamma)$, 
the space of continuous, piecewise affine functions 
with $\int_{\Gamma} f(x) \, d\mu(x) = 0$, is a subset of $\Zh_{\mu}(\Gamma)$.  
We will now see that both $\BDV_{\mu}(\Gamma)$ and $V_{\mu}(\Gamma)$ 
are contained in $\Dir_{\mu}(\Gamma)$, the completion of $\Zh_{\mu}(\Gamma)$
under the Dirichlet norm.  

\begin{proposition} \label{Prop23}
Under the Dirichlet norm $\|f\|_{\Dir} = \<f,f\>_{\Dir}^{1/2}$,

$A)$ $\CPA_{\mu}(\Gamma)$ is dense in $\BDV_{\mu}(\Gamma)$\ , 

$B)$ $\CPA_{\mu}(\Gamma)$ is dense in $V_{\mu}(\Gamma)$\ .  
\end{proposition}

\vskip .1 in
\noindent{Before giving the proof, we need the following lemma.}

\begin{lemma} \label{Lem24}
Fix $A, B \in \CC$, and let $[a,b]$ be a closed interval.  Let $\cH$
be the set of continuous functions $f:[a,b] \rightarrow \CC$ such
that $f(a) = A$, $f(b) = B$, $f$ is differentiable on $(a,b)$,
and $f^{\prime}(x) \in L^2([a,b])$.  Let $h(x) \in \cH$ be the unique affine
function with $h(a) = A$, $h(b) = B$.  Then for each $f \in \cH$,
\begin{equation*}
\int_a^b |f^{\prime}(x)|^2 \, dx \ge \ 
\int_a^b |f^{\prime}(x)-h^{\prime}(x)|^2 \, dx \ .
\end{equation*}
\end{lemma}

\begin{proof}  Put $C = (B-A)/(b-a)$.  Then $h^{\prime}(x) = C$
for all $x$.  For each $f \in \cH$,
\begin{eqnarray}
\int_a^b f^{\prime}(x) \overline{h^{\prime}(x)} \, dx
& = & \int_a^b f^{\prime}(x) \overline{C} dx
\ = \ (f(b)-f(a)) \cdot \overline{C} \notag  \\
& = & (B-A) \overline{C} \ = \ (b-a) |C|^2 \ \in \RR \ . \label{hhx1}
\end{eqnarray}
In particular
\begin{equation}
\int_a^b h^{\prime}(x) \overline{h^{\prime}(x)} \, dx
\ = \ (b-a) |C|^2 \ . \label{hhx2}
\end{equation}
Fix $f \in \cH$.  Using (\ref{hhx1}) and (\ref{hhx2}) we see that
\begin{eqnarray*}
\int_a^b |f^{\prime}(x) -h^{\prime}(x)|^2
   & = & \int_a^b (f^{\prime}(x)-h^{\prime}(x))
                (\overline{f^{\prime}(x)-h^{\prime}(x)}) \, dx \\
   & = & \int_a^b |f^{\prime}(x)|^2 \, dx - (b-a)|C|^2
\end{eqnarray*}
Since $(b-a)|C|^2 > 0$, the lemma follows.
\end{proof}

\vskip .1 in
\begin{proof} {of Proposition \ref{Prop23}.}

A) Put $\nu = \Delta(f)$, so $f = \varphi_{\mu}(\nu)$ and
$\<f,f\>_{\Dir} = \iint g_{\mu}(x,y) \, d\nu(x) \overline{d\nu(y)}$.
Let $\nu_1, \nu_2, \ldots$ be a bounded sequence of discrete measures
of total mass zero 
which converge weakly to $\nu$.
Put $f_n = \varphi_{\mu}(\nu_n)$.  Then $f_n \in \CPA_{\mu}(\Gamma)$, 
and $\<f-f_n,f-f_n\>_{\Dir}$ converges to $0$ by 
Lemma~\ref{LemProductConvergence}, since the measures 
$\nu_n \times \nu_n$, $\nu \times \nu_n$, and $\nu_n \times \nu$ 
converge weakly to $\nu \times \nu$ on $\Gamma \times \Gamma$.


\vskip .1 in

B)  Take $f \in V_{\mu}(\Gamma)$, and fix $\varepsilon > 0$.   
Let $X \subset \Gamma$
be the finite set of points where $f^{\prime}(x)$ is not defined.
After enlarging $X$, we can assume that it contains a vertex set for $\Gamma$.
For each $p \in X$ and $\delta > 0$, 
let $N(p,\delta) = \{ z \in \Gamma : d(p,z) < \delta \}$ be the
$\delta$-neighborhood of $p$ under the canonical metric on $\Gamma$.

Since $f^{\prime} \in L^2(\Gamma)$, there is
a $\delta_1 > 0$ such that
\begin{equation*}
\int_{\bigcup_{x \in X} N(x,\delta_1)}
         f^{\prime}(x) \overline{f^{\prime}(x)} \, dx  \ < \ \varepsilon \ .
\end{equation*}
Put $U  = \bigcup_{x \in X} N(x,\delta_1)$, and put $K = \Gamma \backslash U$.
Then $K$ is compact.  Since $f^{\prime}(x)$ is continuous on $K$, it is
uniformly continuous.  Let $0 < \delta_2 < \delta_1$ be small enough that
$|f^{\prime}(x)-f^{\prime}(y)| < \sqrt{\varepsilon}$ for all $x,y \in K$
with $d(x,y) < \delta_2$.  Let $Y = \{y_1, \ldots, y_M\} \subset K$
be a finite set such that for each $x \in K$, there is some $y_i \in Y$
with $d(x,y_i) < \delta_2$.  After enlarging $Y$ we can assume it contains
all the boundary points of $K$.

The set $X \cup Y$ partitions $\Gamma$ into a finite union of segments.
Let $F : \Gamma \rightarrow \CC$ be the unique continuous function 
such that $F(z) = f(z)$ for each $z \in X \cup Y$, and which is affine
on $\Gamma \backslash (X \cup Y)$.  

Let $(a,b)$ be one of the segments in $\Gamma \backslash (X \cup Y)$.
If $(a,b) \subset K$, then necessarily $d(a,b) < \delta_2$.
For each $x \in (a,b)$ we have $F^{\prime}(x) = (f(b)-f(a))/(b-a)$.
On the other hand, by the Mean Value Theorem there is a point
$x_0 \in (a,b)$ with $f^{\prime}(x_0) = (f(b)-f(a))/(b-a)$.
Since $x, x_0 \in [a,b]$ it follows that $d(x,x_0) < \delta_2$.
Hence
\begin{equation*}
|f^{\prime}(x) - F^{\prime}(x)| \ = \ |f^{\prime}(x)-f^{\prime}(x_0)|
                                \ < \ \sqrt{\varepsilon} \ .
\end{equation*}
It follows that 
\begin{equation*}
\int_a^b |f^{\prime}(x)-F^{\prime}(x)|^2 \, dx \ \le \ (b-a) \varepsilon \ .
\end{equation*}
On the other hand, if $(a,b) \subset U$, then by Lemma \ref{Lem24}
\begin{equation*}
\int_a^b |f^{\prime}(x)-F^{\prime}(x)|^2 \, dx \ \le \
\int_a^b |f^{\prime}(x)|^2 \, dx \ .
\end{equation*}
Summing over all intervals $(a,b)$, we see that
\begin{eqnarray*}
\int_{\Gamma} |f^{\prime}(x)-F^{\prime}(x)|^2 \, dx
& = & \int_K |f^{\prime}(x)-F^{\prime}(x)|^2 \, dx
         + \int_U |f^{\prime}(x)-F^{\prime}(x)|^2 \, dx \\
& \le & (\ell(\Gamma) + 1) \cdot \varepsilon
\end{eqnarray*}
where $\ell(\Gamma)$ is the total length of $\Gamma$.

Here $F(x)$ need not lie in $\CPA_{\mu}(\Gamma)$, but 
$F_{\mu}(x) =  F(x) - \frac{1}{\ell(\Gamma)} \int_{\Gamma} F(x) \, dx$ 
does, and $\<f-F_{\mu},f-F_{\mu}\>_{\Dir} = \<f-F,f-F\>_{\Dir}$.  
Since $\varepsilon$ is arbitrary, the proposition follows.
\end{proof}

\vskip .1 in
Recall that $\Dir_{\mu}(\Gamma)$ is the completion of $\Zh_{\mu}(\Gamma)$ 
under $\| \ \|_{\Dir}$, with norm $\| \ \|_{\Dir}$ extending the Dirichlet norm
on $\Zh_{\mu}(\Gamma)$.
By Proposition \ref{Prop23}, both $\BDV_{\mu}(\Gamma)$
and $V_{\mu}(\Gamma)$ are dense subspaces of 
$\Dir_{\mu}(\Gamma)$. 
In particular, $\Dir_{\mu}(\Gamma)$ can be defined as the
completion of any of the spaces $\CPA_{\mu}(\Gamma)$, $V_{\mu}(\Gamma)$,
$\Zh_{\mu}(\Gamma)$, or $\BDV_{\mu}(\Gamma)$
under the Dirichlet norm.

Recall from Lemma~\ref{lemma:PI}
that there is a constant $C_\infty$ such that for each $f \in \Zh_{\mu}(\Gamma)$,
$\|f\|_{\infty} \ \le \ C_\infty \cdot \|f\|_{\Dir}$.
This implies that the abstract Hilbert space $\Dir_\mu(\Gamma)$ can
be identified with a space of continuous functions on $\Gamma$:

\vskip .1 in
\begin{corollary} \label{Cor26b}
We can uniquely identify $\BDV_{\mu}(\Gamma)$ and $V_{\mu}(\Gamma)$ with
subspaces of $\Dir_{\mu}(\Gamma)$, where the latter is regarded as 
a subset of $\cC(\Gamma)$ via the embedding defined in Corollary~\ref{Cor26}.  
\end{corollary}             

\begin{proof}  We must show that each $g \in \BDV_{\mu}(\Gamma)$ 
(resp. $V_{\mu}(\Gamma)$) is taken to itself under the embedding 
$\iota_{\mu} : \Dir_{\mu}(\Gamma) \hookrightarrow \cC(\Gamma)$ 
defined in the proof of Corollary~\ref{Cor26}.  

If $g \in \BDV_{\mu}(\Gamma)$, 
let $\nu = \Delta(g)$, so $g = \varphi_{\mu}(\nu)$.
Let $\nu_1, \nu_2, \ldots$ be a bounded sequence of
discrete measures of total mass zero converging weakly to $\nu$.

By Proposition \ref{Prop7} the functions $g_n(x) = \varphi_{\mu}(\nu_n)$
converge uniformly to $g(x)$.  Moreover, each $g_n$ belongs to
$\CPA_{\mu}(\Gamma) \subset \Dir_{\mu}(\Gamma)$.  
By Corollary~\ref{Cor26} $\iota_{\mu}(g_n) = g_n$, 
and by Lemma \ref{lemma:PI} the $g_n$ form a Cauchy sequence under 
$\| \cdot \|_{\infty}$.  By definition, the image $G = \iota_{\mu}(g)$ 
is the limit of this Cauchy sequence.  Hence $G = g$.  

Examination of the proof of Lemma \ref{lemma:PI}
shows that the conclusions of that Lemma hold for $V_{\mu}(\Gamma)$ 
as well as for $\Zh_{\mu}(\Gamma)$.  Hence the fact that 
$\iota_{\mu} (g) = g$ for all $g \in V_{\mu}(\Gamma)$ follows by the proof
of Corollary~\ref{Cor26}.  
\end{proof}

\vskip .1 in

Let $\mu_1$ and $\mu_2$ be real, signed Borel measures on $\Gamma$ with 
$\mu_1(\Gamma) = \mu_2(\Gamma) = 1$.  There is a canonical projection  
$\pi_{1,2} : \Dir_{\mu_1}(\Gamma) \rightarrow \Dir_{\mu_2}(\Gamma)$, 
defined by 
\begin{equation*}
\pi_{1,2}(f) \ = \ f(x) - \int_{\Gamma} f(x) \, d\mu_2(x) \ .
\end{equation*}
It is easy to check that $\pi_{2,1} = \pi_{1,2}^{-1}$.  

\begin{proposition} \label{Prop31}
$\pi_{1,2}$ is an isometry from $\Dir_{\mu_1}(\Gamma)$ onto 
$\Dir_{\mu_2}(\Gamma)$ relative to the Dirichlet norms on those spaces, 
which takes $V_{\mu_1}(\Gamma)$ to $V_{\mu_2}(\Gamma)$,
$\CPA_{\mu_1}(\Gamma)$ to $\CPA_{\mu_2}(\Gamma)$, 
and $\BDV_{\mu_1}(\Gamma)$ to $\BDV_{\mu_2}(\Gamma)$.  
\end{proposition}

\begin{proof} Since $\pi_{1,2}$ simply translates each function by a constant,
it is clear that it takes $V_{\mu_1}(\Gamma)$ to $V_{\mu_2}(\Gamma)$,
$\CPA_{\mu_1}(\Gamma)$ to $\CPA_{\mu_2}(\Gamma)$, and $\BDV_{\mu_1}(\Gamma)$
to $\BDV_{\mu_2}(\Gamma)$.  Since the Dirichlet pairings on $V(\Gamma)$ and
$\BDV(\Gamma)$ are invariant under translation by constants, we have
\begin{equation}  \label{bnmq}
\<\pi_{1,2}(f),\pi_{1,2}(f)\>_{\Dir} \ = \ \<f,g\>_{\Dir} 
\end{equation}
for all $f,g \in V_{\mu_1}(\Gamma)$ and all $f,g \in \BDV_{\mu_1}(\Gamma)$.  
These spaces are dense in $\Dir_{\mu_1}(\Gamma)$, so (\ref{bnmq}) holds 
for all $f, g \in \Dir_{\mu_1}(\Gamma)$.  
     This shows that $\pi_{1,2}$ is an injective isometry from 
$\Dir_{\mu_1}(\Gamma)$ into $\Dir_{\mu_2}(\Gamma)$.  The surjectivity 
follows by considering the inverse map $\pi_{2,1}$.  
\end{proof} 
 
\vskip .1 in
This simple proposition has two interesting consequences.

First, if $\{F_n\}$ is a complete orthonormal basis for 
$\Dir_{\mu_1}(\Gamma)$ given by eigenfunctions of 
$\Delta$ in $\Dir_{\mu_1}(\Gamma)$, then $\{\pi_{1,2}(F_n)\}$ 
is a complete orthonormal basis for $\Dir_{\mu_2}(\Gamma)$.  
Note however that the functions $\pi_{1,2}(F_n)$ are in general not
eigenfunctions of $\Delta$ in $\Dir_{\mu_2}(\Gamma)$.  Furthermore,
although the $F_n$ are mutually orthogonal relative to the $L^2$ 
inner product $\<\ , \>$,  the projected functions $\pi_{1,2}(F_n)$
are not in general orthogonal relative to $\<\ ,\ \>$ even though
they remain orthogonal relative to $\<\ ,\ \>_{\Dir}$.

Second, if we define the `Dirichlet space' $\Dir(\Gamma)$ by
\begin{equation*}
\Dir(\Gamma) \ = \ \Dir_{\mu_1}(\Gamma) \oplus \CC \cdot 1 
\end{equation*}
then $\Dir(\Gamma)$ contains the spaces $\Dir_{\mu_2}(\Gamma)$
for all $\mu_2$, and in particular is independent of $\mu_1$.   
It appears to be an important space to study.   

\vskip .1 in

\section{Equivalence between notions of Eigenfunctions}
\label{LaplacianEigenfunctionsSectionII}

\vskip .1 in

We will now show that the eigenfunctions of $\Delta$ in $\Dir_{\mu}(\Gamma)$
constructed in Theorem \ref{MainPairingTheorem} coincide with the 
eigenfunctions of $\varphi_{\mu}$ given by Theorem \ref{Thm9}.  

This equivalence provides new information:  
it implies that the eigenfunctions of $\Delta$ in $\Dir_{\mu}(\Gamma)$ 
are ``smooth enough'' to belong to $\BDV_{\mu}(\Gamma)$, 
and that the eigenfunctions
of $\varphi_{\mu}$, which are known to provide a Hilbert space 
basis for $L^2_{\mu}(\Gamma)$,
also form a basis for $\Dir_{\mu}(\Gamma)$.






\begin{theorem} \label{EquivThm}
Given a nonzero function $f \in L^2(\Gamma)$, the following are equivalent:

$A)$  $f$ is an eigenfunction of $\Delta$ in $\Dir_{\mu}(\Gamma)$ with 
             eigenvalue $\lambda > 0$ 
             $(\S$\ref{EigenfunctionViaPairingsSection},  
                         Definition \ref{DefDirDF}$)$,   
             
$B)$  $f$ is an eigenfunction of $\varphi_{\mu}$ with 
             eigenvalue $\alpha = 1/\lambda > 0$ 
             $(\S$\ref{IntegralOperatorEigenfunctions},
                      Definition \ref{DefVarphiEV}$)$, 
             
$C)$  $f$ is an eigenfunction of $\Delta$ in $\BDV_{\mu}(\Gamma)$
             with eigenvalue $\lambda > 0$ 
             $(\S$\ref{EFLaplacianSection}, Definition \ref{EFDef}$)$.   
\end{theorem}

\begin{proof}  
We know that (B) $\Leftrightarrow$ (C) by Proposition \ref{Prop32}, 
so it suffices to show that (A) $\Leftrightarrow$ (B).

(B) $\Rightarrow$ (A).  
Suppose first that $0 \ne f$ is an eigenfunction of 
$\varphi_{\mu}$ with eigenvalue $\alpha \ne 0$.  By Theorem \ref{Thm9},
$f \in \BDV_{\mu}(\Gamma)$.  Theorem \ref{Thm16} shows that $\alpha > 0$.
By Lemma \ref{Lem10},
\begin{equation} \label{FAZ1}
\alpha \cdot \Delta(f) \ = \ \Delta(\varphi_{\mu}(f)) 
                     \ = \ f(x)\, dx - (\int_{\Gamma} f(x) \, dx) \, \mu \ .
\end{equation}
To see that $f$ is an eigenfunction of $\Delta$ in $\Dir_{\mu}(\Gamma)$
with eigenvalue $1/\alpha$,
we must show that $\<f,g\>_{\Dir} = {\frac{1}{\alpha}} \<f,g\>_{L^2}$ for each
$g \in \Dir_{\mu}(\Gamma)$.  Since $\Zh_{\mu}(\Gamma)$ is dense in 
$\Dir_{\mu}(\Gamma)$, we can assume that $g \in \Zh_{\mu}(\Gamma)$.  
By (\ref{FAZ1}), 
\begin{equation*}
\alpha  \int_{\Gamma} \overline{g(x)} \, d\Delta(f)(x) 
    \ = \ \int_{\Gamma} \overline{g(x)} f(x) \, dx \, - \, 
               (\int_{\Gamma} f(x) \, dx) \cdot 
                      (\int_{\Gamma} \overline{g(x)} \, d\mu(x) ) \ .
\end{equation*}
Here $g(x) \in \Zh_{\mu}(\Gamma)$ and $\mu$ is real, 
so $\overline{g(x)} \in \Zh_{\mu}(\Gamma)$
and $\int_{\Gamma} \overline{g(x)} \, d\mu(x) = 0$.  Thus 
\begin{equation*}
\<f,g\>_{\Dir} \ = \ 
\int_{\Gamma} \overline{g(x)} \, d\Delta(f)(x)   
\ = \ \frac{1}{\alpha} \, \int_{\Gamma} \overline{g(x)} f(x) \, dx  
\ = \ \frac{1}{\alpha} \, \<f,g\>_{L^2} 
\ .
\end{equation*}

(A) $\Rightarrow$ (B).  Conversely, suppose that $0 \ne f$ 
is an eigenfunction of $\Delta$ 
in $\Dir_{\mu}(\Gamma)$ with eigenvalue $\lambda$.  
By Theorem \ref{MainPairingTheorem}, necessarily  $\lambda > 0$.   
Let $g_1, g_2, \ldots$ be an orthonormal basis for $L^2(\Gamma)$ consisting of
eigenfunctions of $\varphi_{\mu}$, as given by Theorem \ref{Thm9}. 
We can expand $f = \sum_n c_n g_n$ in $L^2(\Gamma)$.  
Since $\int f(x) \, d\mu = 0$, Proposition \ref{Prop11} 
shows that for each $n$ with $c_n \ne 0$ the corresponding 
eigenfunction $g_n$ has a nonzero eigenvalue $\alpha_n$, 
and Theorem \ref{Thm9} says that $g_n$ belongs to $\BDV_{\mu}(\Gamma)$.    
 
Since $\BDV_{\mu}(\Gamma) \subset \Dir_{\mu}(\Gamma)$, 
for each such $g_n$ we have 
\begin{equation*}
\<f,g_n\>_{\Dir} \ = \ \lambda \<f,g_n\>_{L^2} \ = \ \lambda \cdot c_n \ .
\end{equation*}
On the other hand, $g_n$ 
is an eigenfunction of $\Delta$ with eigenvalue $1/\alpha_n$, 
and $\alpha_n$ is real by Theorem \ref{Thm9}, so 
\begin{equation*}
\<f,g_n\>_{\Dir} \ = \ \overline{\<g_n,f\>}_{\Dir}
                 \ = \ \frac{1}{\alpha_n} \overline{\<g_n,f\>}_{L^2}
                 \ = \ \frac{1}{\alpha_n} \<f,g_n\>_{L^2}  
                 \ = \ \frac{1}{\alpha_n} \cdot c_n \ . 
\end{equation*}                 
Thus, the only eigenfunctions $g_n$ for which $c_n \ne 0$ are ones
satisfying $1/\alpha_n = \lambda$.  Theorem \ref{Thm9} says that there
are only finitely many such $n$, so 
\begin{equation*}
f \ = \ \sum_{1/\alpha_n = \lambda} c_n g_n \ \in \ \BDV_{\mu}(\Gamma) \ .
\end{equation*}
Since each $g_n$ is an eigenfunction of $\varphi_{\mu}$ with eigenvalue 
$\alpha_n = 1/\lambda$, so is $f$. 
\end{proof}
    
\vskip .1 in

\begin{corollary}  \label{BDVCor}
Each eigenfunction 
of $\Delta$ in $\Dir_{\mu}(\Gamma)$ belongs to $\BDV_{\mu}(\Gamma)$.
\end{corollary}

\begin{corollary} \label{BasisCor}
The eigenfunctions of $\varphi_{\mu}$ in $\BDV_{\mu}(\Gamma)$ contain
a complete orthonormal basis for $\Dir_{\mu}(\Gamma)$.   More precisely,
let $\{f_n\}$ be a complete $L^2$-orthonormal basis of
eigenfunctions for $\varphi_{\mu}$ in $\BDV_{\mu}(\Gamma)$,
and put $F_n = f_n/\sqrt{\lambda_n}$.  Then 
$\{F_n\}$ is a complete orthonormal basis for $\Dir_{\mu}(\Gamma)$.
\end{corollary}

\section{The eigenfunction expansion of $g_{\mu}(x,y)$, and applications}
\label{ExpansionSection} 

\vskip .1 in

Let $\{f_n\}$ be a basis for $\BDV_{\mu}(\Gamma)$ consisting of eigenfunctions
of $\Delta$, normalized so that $\|f_n\|_2 = 1$ for each $n$.  

We now  consider the expansion of $g_{\mu}(x,y)$ in terms of the
eigenfunctions $\{f_n\}$.
For each $y \in \Gamma$, put $G_y(x) = g_{\mu}(x,y)$.  In $L^2(\Gamma)$,
we can write $G_y(x) = \sum_{n=1}^{\infty} a_n f_n(x)$, with 
\begin{equation*}
a_n \ = \ \<G_y(x),f_n(x)\>
    \ = \ \int_{\Gamma} g_{\mu}(x,y) \overline{f_n(x)} \, dx 
    \ = \ \overline{\varphi_{\mu}(f_n)(y)}
    \ = \ \frac{1}{\lambda_n} \overline{f_n(y)}.
\end{equation*}
Thus
\begin{equation} \label{Form2}
g_{\mu}(x,y)
   \ = \ \sum_{n=1}^{\infty} \frac{f_n(x) \overline{f_n(y)}}{\lambda_n} \ ,
\end{equation}
where for each $y$ the series converges to $g_{\mu}(x,y)$ in $L^2(\Gamma)$.

\begin{proposition} \label{Prop18} { \ \ \ }

The series $\sum_{n=1}^{\infty} \frac{f_n(x) \overline{f_n(y)}}{\lambda_n}$ 
converges uniformly to $g_{\mu}(x,y)$ for all $x, y \in \Gamma$.  
\end{proposition} 

\begin{proof}  This follows from a classical theorem of Mercer 
(see \cite{RN}, p.245).  Mercer's theorem asserts that if $K$ is a 
compact measure space and $A(x,y)$ is a continuous, symmetric kernel for which 
the corresponding integral operator $A: L^2(K) \rightarrow L^2(K)$ is 
positive (that is, $\<Af,f\> \ge 0$ for all $f \in L^2(K)$), 
then the $L^2$ eigenfunction expansion of $A(x,y)$ converges 
uniformly to $A(x,y)$.  The proof in (\cite{RN}) is given when 
$K = [a,b]$ is a closed interval, but the argument is general.  

In our case $g_{\mu}(x,y)$ is continuous, real-valued, and symmetric, 
and the integral operator $\varphi_{\mu}$ is positive by Theorem \ref{Thm16}.
\end{proof}

\vskip .1 in
Several important facts follow from this. 

\begin{corollary} \label{Cor19}
$\sum_{n=1}^{\infty} 1/\lambda_n 
\ = \int_{\Gamma} g_{\mu}(x,x) \, dx < \ \infty$ .
\end{corollary}

\begin{proof}
By Proposition \ref{Prop18} 
$\sum_{n=1}^{\infty} f_n(x) \overline{f_n(x)}/\lambda_n$ 
converges uniformly to $g_{\mu}(x,x)$ on $\Gamma$.  
Since $\<f_n,f_n\> = 1$ for each $n$, 
\begin{equation*}
\int_{\Gamma} g_{\mu}(x,x) \, dx
   \ = \ \sum_{n=1}^{\infty} \frac{1}{\lambda_n}
            \int_{\Gamma} f_n(x) \overline{f_n(x)} \, dx 
   \ = \ \sum_{n=1}^{\infty} \frac{1}{\lambda_n} \ .
\end{equation*}
\end{proof}

\vskip .1 in

\begin{corollary} \label{Cor20}
$g_{\mu}(x,x) \ge 0$ for all $x \in \Gamma$.  
Indeed, $g_{\mu}(x,x) > 0$ 
for all $x$ unless $\mu = \delta_{x_0}(x)$ for some $x_0$, 
in which case $g_{\mu}(x,x) > 0$ for all $x \ne x_0$.  
\end{corollary}

\begin{proof}
It follows from 
$g_{\mu}(x,x) = \sum_{n=1}^{\infty} f_n(x) \overline{f_n(x)}/\lambda_n$ 
that $g_{\mu}(x,x) \ge 0$ for all $x$.  If $g_{\mu}(x_0,x_0) = 0$ 
for some $x_0$, then $f_n(x_0) = 0$ for all $n$.  But then 
$g_{\mu}(x_0,y) = \sum_n f_n(x_0) \overline{f_n(y)}/\lambda_n = 0$
for all $y$, which means that 
$0 = \Delta_y(g_{\mu}(x_0,y)) = \delta_{x_0}(y) - \mu(y)$.
Thus $\mu = \delta_{x_0}(y)$.  If there were another point $x_1$ with
$g_{\mu}(x_1,x_1) = 0$, then $\delta_{x_0}(y) = \mu = \delta_{x_1}(y)$,
so $x_1 = x_0$.    
\end{proof}

\vskip .1 in  
\begin{corollary} \label{Cor21}
$\sum_{n=1}^{\infty} |f_n(x) \overline{f_n(y)}/\lambda_n|$  converges
uniformly for all $x, y$.
\end{corollary}

\begin{proof} Note that 
\begin{equation*}
\sum_{n=1}^{\infty} |\frac{f_n(x)}{\sqrt{\lambda_n}}|^2
 \ = \ \sum_{n=1}^{\infty} f_n(x) \overline{f_n(x)}/\lambda_n
  \ = \ g_{\mu}(x,x) 
\end{equation*}
for each $x$, and the convergence is uniform for all $x \in \Gamma$.  
Hence for any $\varepsilon > 0$ there is an $N$ such that
$\sum_{n=N}^{\infty} |\frac{f_n(x)}{\sqrt(\lambda_n)}|^2 \le \varepsilon$.  

By the Cauchy-Schwarz inequality, 
\begin{equation*}
\sum_{n=1}^{\infty} |f_n(x) \overline{f_n(y)}/\lambda_n| \ \le \ 
   \big(\sum_{n=1}^{\infty} |\frac{f_n(x)}{\sqrt{\lambda_n}}|^2 \big)^{1/2}
   \big(\sum_{n=1}^{\infty} |\frac{f_n(y)}{\sqrt{\lambda_n}}|^2 \big)^{1/2} 
\end{equation*}
converges.  Another application of Cauchy-Schwarz gives 
\begin{equation*}
\sum_{n=N}^{\infty} |f_n(x) \overline{f_n(y)}/\lambda_n| \ \le \ \varepsilon \ ,  
\end{equation*} 
so the convergence is uniform.
\end{proof}

\begin{corollary} \label{Cor22}
With the eigenfunctions normalized so that $\|f_n\|_2 = 1$ for
all $n$, we have $\|f_n\|_{\infty} = o(\sqrt{\lambda_n})$ 
as $n \rightarrow \infty$.
\end{corollary}

\begin{proof} As in the proof of the preceding corollary, 
for any $\varepsilon > 0$ there is an $N$ so that 
$\sum_{n=N}^{\infty} |\frac{f_n(x)}{\sqrt{\lambda_n}}|^2 \le \varepsilon$
for all $x$.  Hence for each $n \ge N$, 
$\|f_n\|_{\infty} \le \sqrt{\varepsilon} \cdot \sqrt{\lambda_n}$.
\end{proof}

Let $\{f_n\}_{1 \le n < \infty}$ be an $L^2$-orthonormal basis for 
$\BDV_{\mu}(\Gamma)$ consisting of eigenfunctions for $\Delta$.
We will now show that the $L^2$-expansion of each
$F \in \Dir_{\mu}(\Gamma)$ in terms of the $f_n$ converges uniformly to $F(x)$
on $\Gamma$.

Suppose $F \in \Dir_{\mu}(\Gamma)$.
If $\mu = g(x) dx$ for some $g \in L^2(\Gamma)$, then by
Proposition \ref{Prop11} the eigenfunctions
$f_n$ form an orthonormal basis for $(\CC g)^{\perp}
= \{f \in L^2(\Gamma) : \int_{\Gamma} f(x) \overline{g(x)} \, dx = 0\}$.
Since $F \in \cC(\Gamma)$ and $\int_{\Gamma} F(x) \, d\mu(x) = 0$,
we have $F \in (\CC g)^{\perp}$.  If $\mu$ is not of this form,
then the $f_n$ form an orthonormal basis for $L^2(\Gamma)$.

In either case we can expand
$F(x) = \sum_{n=1}^{\infty} c_n f_n(x)$ in $L^2(\Gamma)$.

\begin{corollary} \label{Cor27} { \ \ \ }

For each $F \in \Dir_{\mu}(\Gamma)$, the series
$\sum_{n=1}^{\infty} c_n f_n(x)$ converges uniformly to $F(x)$
on $\Gamma$.
\end{corollary}

\begin{proof}
Let $\lambda_n$ be the eigenvalue corresponding to $f_n$,
and put $F_n(x) = f_n(x)/\sqrt{\lambda_n}$, so $\<F_n,F_n\>_{\Dir} = 1$.
Since the $F_n$ are eigenfunctions of $\Delta$ belonging
to $\BDV_{\mu}(\Gamma)$ by Corollary~\ref{BasisCor} and Theorem~\ref{EquivThm}, 
they are mutually orthogonal under $\< \ ,\ \>_{\Dir}$.

By assumption, we have
\begin{equation*}
F(x) \ = \ \sum_n c_n f_n(x) \ = \ \sum_n \sqrt{\lambda_n} c_n F_n(x)
\end{equation*}
in $L^2(\Gamma)$.
By Parseval's inequality in $\Dir_{\mu}(\Gamma)$,
\begin{equation} \label{prt1}
\sum_n |\sqrt{\lambda_n} c_n|^2 \ \le \ \<F,F\>_{\Dir} \ .
\end{equation}
By Corollary \ref{Cor21}, the series
\begin{equation*}
\sum_n |F_n(x)|^2 \ = \ \sum_n |f_n(x) \overline{f_n(x)}/\lambda_n|
\end{equation*}
converges uniformly on $\Gamma$.  Hence, for each $\varepsilon > 0$,
there is an $N$ such that $\sum_{n \ge N} |F_n(x)|^2 \le \varepsilon$.
By (\ref{prt1}) and the Cauchy-Schwarz inequality,
\begin{equation*}
\sum_{n \ge N} |c_n f_n(x)| \ \le \
    (\<F,F\>_{\Dir})^{1/2} \cdot \sqrt{\varepsilon}
\end{equation*}
for all $x \in \Gamma$.

Thus the series $\sum_n c_n f_n(x)$ converges uniformly and absolutely
to a function $G(x) \in \cC(\Gamma)$.  Since $F(x)$ and $G(x)$
are continuous, and $F = G$ in $L^2(\Gamma)$,
it follows that $F(x) = G(x)$.
\end{proof}


We can use Proposition~\ref{Prop18} and Corollary~\ref{Cor27} to give another proof of 
Theorem~\ref{Thm17}(B), the energy minimization principle, as follows.

\vskip .1 in

\begin{proof} (of Theorem~\ref{Thm17}(B))
Let $\nu$ be a bounded measure of total mass 1.  Define the 
$n$-th Fourier coefficient of $\nu$ to be
\[
c_n = \int_\Gamma f_n(x) d\nu(x)\ .
\]

By Corollary \ref{Cor27}, the linear space spanned by the functions $\{ f_n \}$
(which lie in $\BDV_{\mu}(\Gamma)$, and in particular have $\int f_n d\mu = 0$)
is dense in $\Dir_{\mu}(\Gamma)$.  Together with 
$f_0 = 1$, the $f_n$ span a dense subspace of $\Dir(\Gamma)$,
which itself is dense in $\cC(\Gamma)$.  
By the Riesz Representation Theorem, it follows that
$\mu = \nu$ if and only if $c_n = 0$ for all $n \ge 1$.  
By uniform convergence, we find that
\begin{equation*}
\begin{aligned}
\iint g_\mu(x,y) d\nu(x) \overline{d\nu(y)}
& =  \iint \left(  \sum_{n \ge 1} f_n(x) \overline{f_n(y)} / \lambda_n \right) 
     d\nu(x) \overline{d\nu(y)} \\
& =  \sum_{n \ge 1} \left( \int f_n(x) d\nu(x) \right)  \left(  \int
  \overline{f_n(y)} \overline{d\nu(y)} \right) / \lambda_n \\
& =  \sum_{n \ge 1} |c_n|^2 / \lambda_n\ \geq 0,
\end{aligned}
\end{equation*}
with equality if and only if $c_n = 0$ for all $n\geq 1$, as desired.
\end{proof}

\vskip .1 in
\noindent{\bf Application: Lower bounds for $g_{\mu}(x,y)$-discriminant sums.}
\vskip .1 in

From the uniform convergence of the expansion (\ref{Form2}),
we obtain Elkies-type lower bounds 
for discriminant sums formed using $g_{\mu}(x,y)$.  
Such sums occur naturally in arithmetic geometry  
(see for example \cite[Lemma~2.1]{HS}, in which $\Gamma$ is a circle).
     
\begin{proposition}
There exists a constant $C = C(\Gamma,\mu) > 0$ such that for all $N\geq 2$ 
and all $x_1,\ldots,x_N \in \Gamma$, then
\begin{equation*}
\frac{1}{N(N-1)} \sum_{i \neq j} g_\mu(x_i,x_j) \geq -\frac{C}{N}.
\end{equation*}
\end{proposition}

\begin{proof}
Put $M = \sup_{x \in \Gamma} g_\mu(x,x)$.  Then

\begin{equation*}
\begin{aligned}
\sum_{i \neq j} g_\mu(x_i,x_j) 
& =  \sum_{i \neq j} \sum_n f_n(x_i) \overline{f_n(x_j)} / \lambda_n \\
& =  \sum_n | \sum_{i=1}^N f_n(x_i) / \sqrt{\lambda_n} |^2 
- \sum_{i=1}^N g_\mu(x_i,x_i) \\
& \geq - M\cdot N, \\
\end{aligned}
\end{equation*}
which yields the desired conclusion.
\end{proof}

\vskip .1 in

\section{The canonical measure and the tau constant}
\label{CanonicalTauSection}

\vskip .1 in

In this section we discuss some properties of a ``canonical measure'' discovered by
Chinburg and Rumely.  This will not only provide us with a nice
application of the energy pairing, it will also
illustrate why it is useful to allow general measures $\mu$ (rather
than just multiples of $dx$, for example) when discussing eigenvalues of the
Laplacian.

\medskip

For $x,y \in \Gamma$, we denote by $r(x,y)$ the ``effective resistance''
between $x$ and $y$, where $\Gamma$ is considered as a resistive
electric circuit as in \S\ref{jfunctionSection}.  Equivalently, we
have
$$
r(x,y) = j_x(y,y) \left( = j_y(x,x) \right).
$$

For example, if $\Gamma = [0,1]$, then $r(x,y)$ is simply $|x-y|$.  Referring
back to \S\ref{BriefExamplesSection}, note that $-\frac{1}{2} |x-y|$
differs by a constant factor of $\frac{1}{4}$ from the Arakelov-Green's function $g_\mu(x,y)$ for the measure
$\mu = \frac{1}{2}\delta_0 + \frac{1}{2}\delta_1$.   

This observation generalizes to arbitrary metrized graphs as follows
(see \cite[Theorem~2.11]{CR} and \cite{BF}):

\begin{theorem} \label{CanonicalMeasureTheorem} { \ }

$A)$ The probability measure $\mu_{\can} = \Delta_x \left(
  \frac{1}{2} r(x,y) \right)  + \delta_y(x)$ is independent of $y \in \Gamma$.  

$B)$ $\mu_\can$ is the unique measure $\mu$ of total mass 1 on $\Gamma$ for 
which $g_\mu(x,x)$ is a constant independent of $x$.

$C)$ There is a constant $\tau(\Gamma) \in \RR$ such that  
$$
g_{\mu_\can}(x,y) = -\frac{1}{2} r(x,y) + \tau(\Gamma) \ .
$$

\end{theorem}

We call $\mu_\can$ the {\it canonical measure} on $\Gamma$, and we
call $\tau(\Gamma)$ the {\it tau constant} of $\Gamma$.

\medskip

In \cite[Theorem 2.11]{CR}, 
the following explicit formula is given:
\begin{equation}
\label{CanonicalMeasureFormula}
\mu_\can \ = \ 
\sum_{\text{vertices $p$}} (1 - \frac{1}{2} \text{valence}(p)) \, \delta_p(x)
+ \sum_{\text{edges $e$}} \frac{dx}{L(e)+R(e)}
\end{equation}
where 
$L(e)$ is the length of edge $e$ and $R(e)$ is the effective
resistance between the endpoints of $e$ 
in the graph $\Gamma \backslash e$, when the graph is regarded as an
electric circuit with resistances equal to the edge lengths.  
In particular, $\mu_\can$ has a point mass of negative
weight at each branch point $p$,
so it is a positive measure only if $\Gamma$ is a segment or a loop.  

\medskip

The following result is an immediate consequence of
Theorem~\ref{Thm17}.

\begin{corollary}
The canonical measure is the unique measure $\nu$ of total mass 1 on
$\Gamma$ maximizing the integral
\begin{eqnarray}
\iint_{\Gamma \times \Gamma} r(x,y) \, 
                   d\nu(x) \overline{d\nu(y)}.  \label{resistanceenergy}
\end{eqnarray}
\end{corollary}

One can therefore think of the canonical measure as being like an
``equilibrium measure'' on $\Gamma$ (in the sense of capacity theory),
and of $\tau(\Gamma)$ as the corresponding ``capacity'' of $\Gamma$
(with respect to the potential kernel $\frac{1}{2} r(x,y)$).  
Note, however, that $\mu_\can$, unlike equilibrium measures 
in capacity theory, is not necessarily a positive measure.

It can also be shown using Theorem~\ref{Thm17}(A) that there is 
a unique {\em probability} measure maximizing
(\ref{resistanceenergy}) over all positive measures $\nu$ of total mass 1.  
However, we do not know an explicit formula for this measure analogous to 
(\ref{CanonicalMeasureFormula}).

\medskip

The next result follows immediately from Corollary~\ref{Cor19}.

\begin{corollary}
\label{cor:tauisatrace}
If $\ell(\Gamma)$ denotes the total length of $\Gamma$, then
$$
\ell(\Gamma) \cdot \tau(\Gamma) = \sum_{n=1}^\infty \frac{1}{\lambda_n},
$$
where the $\lambda_n$'s are the eigenvalues of the Laplacian with
respect to the canonical measure $($i.e., the eigenvalues of $\Delta$ on
$\Dir_{\mu_\can}(\Gamma))$.  

In particular, if  $\ell(\Gamma) = 1$, then $\tau(\Gamma)$ 
is the trace of the operator $\varphi_{\mu_\can}$.  
\end{corollary}

Especially, $\tau(\Gamma) > 0$. 
Another description of $\tau(\Gamma)$, which shows this directly, 
is as follows.

\begin{lemma}
\label{lemma:tauformula}
For any fixed $y \in \Gamma$, we have
$$
\tau(\Gamma) = \frac{1}{4} \int_\Gamma \left( \frac{\partial}{\partial
    x} r(x,y) \right)^2 dx .
$$
\end{lemma}

\begin{proof}
Fix $y \in \Gamma$, and set $f(x) = \frac{1}{2} r(x,y)$.
Since $\int g_{\mu_\can}(x,y) d\mu_{\can}(x) = 0$, we have
$\tau(\Gamma) = \int f(x) d\mu_{\can}(x)$.  Substituting $\mu_\can =
\delta_y(x) + \Delta_x f(x)$ and noting that $f(y) = 0$, we obtain
$$
\tau(\Gamma) = \int f(x) \Delta f(x) = \la f,f \ra_{\Dir} = \int
\left( f'(x) \right)^2 dx.
$$
\end{proof}

Using Lemma~\ref{lemma:tauformula}, it is not hard to see
that if we multiply all lengths on $\Gamma$ by
a positive constant $\beta$ (obtaining a graph $\Gamma(\beta)$ 
of total length $\beta \cdot \ell(\Gamma)$), 
then $\tau(\Gamma(\beta)) = \beta \cdot \tau(\Gamma)$.  Thus it is
natural to consider the ratio $\tau(\Gamma)/\ell(\Gamma)$, which is
scale-independent.

%

One can show using Lemma~\ref{lemma:tauformula} that for
a metrized graph $\Gamma$ with $n$ edges, 
\begin{equation} \label{FMM1}
\frac{1}{16n} \ell(\Gamma) \ \leq \ \tau(\Gamma) 
\ \leq \ \frac{1}{4} \ell(\Gamma) \ ,
\end{equation}
with equality in the upper bound if and only if $\Gamma$ is a tree.
The lower bound is not sharp.

Experience shows that it is hard to construct 
graphs of total length $1$ and small tau constant.    
The smallest known example was  
found by Phil Zeyliger, with $\tau(\Gamma) \cong .021532$.  
We therefore pose the following conjecture:\footnote{Added June 2005:  
Zubeyir Cinkir has recently proved this conjecture with $C = 1/108$.}

\vskip .1 in

\begin{conjecture} \label{TauBound} 
 There is a universal constant $C>0$ 
such that for all metrized graphs $\Gamma$,
$$\tau(\Gamma) \ \geq \ C \cdot \ell(\Gamma)\ .$$
\end{conjecture}

\vskip .1 in

We remark that by Corollary~\ref{cor:tauisatrace}, a universal positive
lower bound for $\tau(\Gamma)$ over all metrized graphs of length 1 would be implied by a universal
upper bound for $\lambda_1$, the smallest nonzero eigenvalue of the
Laplacian with respect to the canonical measure, on such graphs.  

\vskip .1 in
\noindent{\bf Application:  The operator $\varphi_{dx}$ has minimal trace.}
\vskip .1 in

Given a measure $\mu$ on $\Gamma$ of total mass $1$, 
the trace of the operator $\varphi_{\mu}$ is   
$${\rm Tr}(\varphi_{\mu}) \ = \ \int_{\Gamma} g_{\mu}(x,x) dx 
\ = \ \sum_{n=1}^{\infty} \frac{1}{\lambda_n}, $$
where the $\lambda_n$ are the eigenvalues of $\Delta$ on  
$\Dir_{\mu}(\Gamma)$.  
We can compare the traces for different measures using the
energy pairing (see Theorem \ref{Thm17}):    

\begin{proposition}
\label{TraceComparisonProp}
Let $\mu_1,\mu_2 \in \Meas(\Gamma)$ be measures of total mass 1.
Then 
\begin{equation} \label{FTrForm} 
\begin{aligned}
\Tr(\varphi_{\mu_1}) \ = \ 
\Tr(\varphi_{\mu_2}) - \< dx,dx \>_{\mu_2} + 
\< dx - \mu_1, dx - \mu_1 \> \ .
\end{aligned}
\end{equation}
\end{proposition}

\begin{proof}
We first claim that 
\begin{equation}
\label{gmuidentity}
g_{\mu_1}(x,y) = g_{\mu_2}(x,y) - \phi_{\mu_2}(\mu_1)(x) - \phi_{\mu_2}(\mu_1)(y) 
+ \< \mu_1,\mu_1 \>_{\mu_2}.
\end{equation}
Indeed, it follows from Proposition \ref{Prop6} that the Laplacian with 
respect to $x$ of both sides of (\ref{gmuidentity}) is $\delta_y - \mu_1$, 
so the two sides differ by a constant depending only on $x$.  By symmetry, 
the left and right-hand sides of (\ref{gmuidentity}) also differ by a constant 
depending only on $y$.  Thus (\ref{gmuidentity}) is true up to an additive 
constant.  Integrating both sides with respect to 
$d\mu_1(x) d\mu_1(y)$ shows that the constant is zero, proving the claim.

Setting $x = y$ in (\ref{gmuidentity}) 
and integrating both sides with respect to $dx$ now gives
$\Tr(\varphi_{\mu_1}) = 
\Tr(\varphi_{\mu_2}) - 2 \< dx,\mu_1 \>_{\mu_2} + 
\< \mu_1, \mu_1 \>_{\mu_2}$,
which is equivalent to (\ref{FTrForm}).   
\end{proof}

\vskip .1 in

Suppose $\Gamma$ has total length 1.  Setting $\mu_1 = \mu$ and $\mu_2 = dx$ in 
Proposition~\ref{TraceComparisonProp}, 
we find that of all measures $\mu$ on $\Gamma$, 
the $dx$ measure is the one for which $\varphi_{\mu}$ has minimal trace:  

\begin{corollary} \label{MinTrace} 
If $\ell(\Gamma) = 1$ and $\mu \in \Meas(\Gamma)$ has total mass 1, then
\begin{equation*}
\begin{aligned}
\Tr(\varphi_{\mu}) \
&= \ \Tr(\varphi_{dx}) + \< dx - \mu, dx - \mu \> \\
         &\geq \ \Tr(\varphi_{dx}) \ ,
\end{aligned}
\end{equation*}
with equality if and only if $\mu = dx$.  
\end{corollary}

\begin{proof}
By Theorem \ref{Thm17}  $\< dx - \mu, dx - \mu \> \ge 0$, with
equality if and only if $\mu = dx$.
\end{proof}

Since $\Tr(\varphi_{dx})$ is minimal, is useful to have a formula for it:  

\begin{corollary} \label{dxtraceFormula}
If $\ell(\Gamma) = 1$, then 
$\Tr(\varphi_{dx}) 
= \frac{1}{2} \iint_{\Gamma \times \Gamma} r(x,y) \, dx \, dy$ \ .
\end{corollary}

\begin{proof}
Taking $\mu_1 = dx$ and $\mu_2 = \mu_{\can}$ in Proposition 
\ref{TraceComparisonProp} gives 
\begin{equation} \label{FVVV1}
\Tr(\varphi_{dx}) \ = \ \Tr(\varphi_{\mu_\can}) - \< dx,dx \>_{\mu_\can} \ .
\end{equation}
Here $\Tr(\varphi_{\mu_\can}) = \tau(\Gamma)$ 
by Corollary \ref{cor:tauisatrace}, 
while by the definition of $\<dx,dx\>_{\mu_\can}$ and 
Theorem \ref{CanonicalMeasureTheorem}(C), 
\begin{equation} \label{FVVV2}
\<dx,dx\>_{\mu_\can} \ 
= \ \iint_{\Gamma \times \Gamma} g_{\mu_{\can}}(x,y) \, dx \, dy \ 
= \ \iint_{\Gamma \times \Gamma} \tau(\Gamma) - \frac{1}{2} r(x,y) 
                             \, dx \, dy \ .
\end{equation}
Combining (\ref{FVVV1}) and (\ref{FVVV2}) gives the result.  
\end{proof} 

In contrast with Conjecture \ref{TauBound}, which asserts that  
$\Tr(\varphi_{\mu_\can}) \ge C > 0$ for all graphs of total length $1$, 
$\Tr(\varphi_{dx})$ can be arbitrarily small.  
For example, if $\Gamma = B_n$ is the ``banana graph'' with two vertices 
connected by $n$ edges of length $1/n$, then using Corollary \ref{dxtraceFormula}
one computes that
$$
\Tr(\varphi_{dx}) \ = \ \frac{n+2}{12 n^2} \ .
$$

\section{Regularity and Boundedness}

\vskip .1 in

In this section, we will study the smoothness and boundedness of eigenfunctions.


\begin{proposition} \label{Prop33}
Suppose $\mu$ has the form $\mu = g(x) dx + \sum_{i=1}^N c_i \delta_{p_i}(x)$
where $g(x)$ is piecewise continuous and belongs to $L^1(\Gamma)$.  Assume
that $X = \{p_1, \ldots, p_N\}$  contains a vertex set for $\Gamma$,
as well as all points where $g(x)$ is discontinuous.
Then each eigenfunction of $\Delta$ in $\BDV_{\mu}(\Gamma)$ belongs to
$\Zh(\Gamma)$. If $g(x)$ is $\cC^m$ on $\Gamma \backslash X$, 
then each eigenfunction is $\cC^{m+2}$ on $\Gamma \backslash X$, and satisfies
the differential equation
\begin{equation*}
\frac{d^2f}{dx^2} + \lambda f(x) \ = \ \lambda C g(x) \ ,
\end{equation*}
\end{proposition}

\begin{proof} If $0 \ne f \in \BDV_{\mu}(\Gamma)$ is an eigenfunction of
$\Delta$ with eigenvalue $\lambda$, then by Proposition \ref{Prop32},
\begin{equation} \label{fbb1}
\Delta(f) \ = \ \lambda \cdot (f(x) dx - C\mu)
          \ = \ \lambda (f(x) - C g(x)) dx
                   - C  \lambda \sum_{i=1}^N c_i \delta_{p_i}(x)
\end{equation}
with $C = \int_\Gamma f(x) dx$.

Here $f(x)$ is continuous, so $\lambda (f(x)-Cg(x))$ is continuous on
$\Gamma \backslash X_g$ and belongs to $L^1(\Gamma)$.  By Proposition
\ref{Prop3}, $f \in \Zh(\Gamma)$ and
\begin{equation*}
\Delta(f) \ = \ \Delta_{\Zh}(f) \ = \ - f^{\prime \prime}(x) dx -
    \sum_p (\sum_{\vv \in \VEC(p)} d_{\vv}f(p)) \, \delta_p(x) \ .
\end{equation*}
Comparing this with (\ref{fbb1}), we see that $\Delta(f)$ has no point masses 
except at points of $X$, and that on each subinterval of $\Gamma \backslash X$
\begin{equation*}
f^{\prime \prime}(x) = \ \lambda (C g(x) - f(x)) \ .
\end{equation*}
Differentiating this inductively shows that on each such subinterval,
$f(x)$ has continuous derivatives up to order $m+2$.  Furthermore, $f(x)$
satisfies the differential equation
 $f^{\prime \prime}(x) + \lambda f(x) =  \lambda C g(x)$.  
\end{proof}

\vskip .1 in
Under somewhat stronger hypotheses, the eigenfunctions are uniformly bounded.  

\begin{proposition} \label{Prop34}
Suppose $\mu = g(x) dx + \sum_{i=1}^N c_i \delta_{p_i}(x)$, where $g(x)$ is
piecewise $\cC^1$ on $\Gamma$ and $g^{\prime}(x) \in L^1(\Gamma)$.
Let $\{f_n\}_{1 \le n < \infty}$ be the eigenfunctions of $\Delta$
in $\BDV_{\mu}(\Gamma)$, 
normalized so that $\|f_n\|_2 = 1$ for each $n$.  Then there is a constant
$B$ such that $\|f_n\|_{\infty} \le B$ for all $n$.
\end{proposition}

\begin{proof} Write $X = \{p_1, \ldots, p_N\}$.  After enlarging $X$, we
can assume that $X$ contains a vertex set for $\Gamma$, as well as all points
where $g^{\prime}(x)$ fails to exist.  Then $\Gamma \backslash X$
is a finite union of open segments;  let their closures be denoted
$e_1, \ldots, e_m$.  Without loss of generality, assume $e_{\ell}$ is isometrically 
parametrized by $[0,L_{\ell}]$, and identify it with that interval.

Let $f = f_n$ be an eigenfunction of $\Delta$ in $\BDV_{\mu}(\Gamma)$,
normalized so that $\|f\|_2 = 1$.  Let $\lambda = \lambda_n$ be the
corresponding eigenvalue.  By Proposition \ref{Prop32},
$f \in \Zh(\Gamma)$ and on each segment $e_{\ell}$ it satisfies the
differential equation
\begin{equation} \label{fm1a}
\frac{d^2 f}{dx^2} + \lambda f \ = \ \lambda C g \ .
\end{equation}
where $C = C_f = \int_{\Gamma} f(x) dx$.
Let $\gamma > 0$ be such that $\gamma^2 = \lambda$,
and write $f_{\ell}$ for the restriction of $f$ to $e_{\ell}$.

By the method of variation of parameters (see \cite{Sim}, p.92),
\begin{equation*}
f_{\ell}(x) \ = \ A_{\ell} e^{i \gamma x} + B_{\ell} e^{-i \gamma x}
                  + C\cdot  h_{\gamma,\ell}(x)
\end{equation*}
for certain constants $A_{\ell}, B_{\ell} \in \CC$, where
\begin{equation*}
h_{\gamma,\ell}(x) \ = \
  \lambda e^{i \gamma x} 
        \int_{x_{\ell}}^x \frac{ - e^{- i \gamma t} g(t)}{-2i \gamma} \, dt
  + \lambda e^{-i \gamma x} 
        \int_{x_{\ell}}^x \frac{e^{i \gamma t} g(t)}{-2i \gamma} \, dt \ 
\end{equation*}
for some fixed $x_{\ell} \in e_{\ell}$. Integrating by parts gives
\begin{equation*}
h_{\gamma,\ell}(x) \ = \ g(x) - g(x_{\ell}) 
       + \frac{1}{2} e^{ i \gamma x} 
               \int_{x_{\ell}}^x e^{-i \gamma t} g^{\prime}(t) \, dt
       - \frac{1}{2} e^{ -i \gamma x} 
               \int_{x_{\ell}}^x e^{i \gamma t} g^{\prime}(t) \, dt  \ .
\end{equation*}
Here $g(x)$ is bounded since
$g^{\prime}(x) \in L^1(\Gamma)$, and for the same reason the last two terms
are also bounded, independent of $\gamma$ and $\ell$.  
Hence the $h_{\gamma,\ell}(x)$ are uniformly bounded.  Let $M$ be such that 
$|h_{\gamma,\ell}(x)| \le M$, for all $\ell$, $\gamma$, and $x$.  
                    
Note that $C = C_f$ is bounded, independent of $\gamma$:  indeed 
\begin{equation*}
|C_f| \ = \ |\int_{\Gamma} f(x) dx|
      \ \le \ \|f\|_1 \ \le \ c(\Gamma) \|f\|_2 \ = \ c(\Gamma) \ .
\end{equation*}
We now will use the fact that $\|f\|_2 = 1$ 
to bound the $A_{\ell}$ and $B_{\ell}$.
Write $\<f,f\>_{\ell} = \int_{e_{\ell}} f(x) \overline{f(x)} \, dx$.  Then
\begin{eqnarray}
1 & = & \<f,f\> \ = \ \sum_{\ell=1}^m \<f,f\>_{\ell} \notag \\
& = & \sum_{\ell=1}^m
   \big( |A_{\ell}|^2 L_{\ell} + |B_{\ell}|^2 L_{\ell} 
           + |C|^2 \<h_{\gamma,\ell},h_{\gamma,\ell}\>_{\ell} \label{ggg1} \\
& & \qquad + 
     2 \Re(A_{\ell} \overline{B_{\ell}} \frac{e^{2i \gamma x}-1}{2 i \gamma})
     + 2 \Re(A_{\ell} \overline{C} 
            \<e^{ i \gamma x}, h_{\gamma,\ell}(x)\>_{\ell})
           \notag  \\
& & \qquad \qquad
   + 2 \Re(B_{\ell} \overline{C} 
              \<e^{- i \gamma x}, h_{\gamma,\ell}(x)\>_{\ell}) \big)\ .
           \notag
\end{eqnarray}
The fourth term is $O(|A_{\ell} B_{\ell}|/\gamma)$, 
and by the Cauchy-Schwarz inequality, 
for sufficiently large $\gamma$ it is bounded by 
$\frac{1}{3}(|A_{\ell}|^2 L_{\ell} + |B_{\ell}|^2 L_{\ell})$. 
The fifth term is bounded by $2 |A_{\ell}| \cdot c(\Gamma) M L_{\ell} $, 
and the sixth term is bounded by $2 |B_{\ell}| \cdot c(\Gamma) M L_{\ell}$. 

Suppose that for some $\ell$ we have 
$\max(|A_{\ell}|,|B_{\ell}|) \ge  12 c(\Gamma) M$. 
Then
\begin{eqnarray*}
\lefteqn{\frac{1}{3} (|A_{\ell}|^2 L_{\ell} + |B_{\ell}|^2 L_{\ell})
         \ \ge \ 
   4 \max(|A_{\ell}|,|B_{\ell}|) \cdot c(\Gamma) M L_{\ell}} \qquad & &  \\   
   & \ge & 
    |2 \Re(A_{\ell} \overline{C} 
           \<e^{ i \gamma x}, h_{\gamma,\ell}(x)\>_{\ell})| +
    |2 \Re(B_{\ell} 
           \overline{C} \<e^{- i \gamma x}, h_{\gamma,\ell}(x)\>_{\ell})| \ .
\end{eqnarray*}       
If $\gamma$ is large enough that 
$2 \Re(A_{\ell} \overline{B_{\ell}} \frac{e^{2i \gamma x}-1}{2 i \gamma})
\le  \frac{1}{3}(|A_{\ell}|^2 L_{\ell} + |B_{\ell}|^2 L_{\ell})$, then 
(\ref{ggg1}) shows that 
$\frac{1}{3}(|A_{\ell}|^2 L_{\ell} + |B_{\ell}|^2 L_{\ell}) \le 1$, from
which it follows that $\max(|A_{\ell}|,|B_{\ell}|) \le \sqrt{3/L_{\ell}}$.
Thus for all sufficiently large $\gamma$ and all $\ell$,
\begin{equation*}
|A_{\ell}|, |B_{\ell}| \ \le \ \max(12 c(\Gamma) M, \sqrt{3/L_{\ell}}) \ .
\end{equation*}

It follows that $\|f_n\|_{\infty}$ is uniformly bounded for all sufficiently 
large $n$, say $n > n_0$.
Trivially $\|f_n\|_{\infty} < \infty$ for each $n \le n_0$, 
and the Proposition follows.
\end{proof}

\vskip .1 in

\section{Examples} \label{ExampleSection} 

\vskip .1 in

In this section, we will explain how the theory developed above 
can be used to compute eigenvalues and eigenfunctions in certain cases, 
and we will illustrate this with several examples.  
When $\mu = dx$, the method we 
are about to describe is well-known (see \cite{Be1}, \cite{Nicaise1}).

\vskip  .1 in
Throughout this section, we suppose that $\mu$ has the same form as in 
Proposition \ref{Prop33}:  
\begin{equation*}
\mu \ = \ g(x) dx + \sum_{i=1}^N c_i \delta_{p_i}(x) \ , 
\end{equation*}
where $g(x)$ is piecewise continuous and belongs to $L^1(\Gamma)$.  
Let $X_g \subset \Gamma$ be a finite set containing a vertex set for 
$\Gamma$, all points where $g(x)$ is not continuous, and all points  
where $\mu$ has a nonzero point mass.  Then $\Gamma \backslash X_g$ 
is a finite union of open segments;  let their 
closures in $\Gamma$ be denoted $e_{\ell}$, for $\ell = 1, \ldots, m$.  
Without loss of generality, assume $e_{\ell}$ is isometrically parametrized by 
$[a_{\ell},b_{\ell}]$, and identify it with that interval.
For each $p \in X_g$, let $v(p)$ be the valence of $\Gamma$ at $p$, 
that is, the number of segments $e_{\ell}$ with $p$ as an endpoint.   

\vskip .1 in
If $f$ is an eigenfunction of $\Delta$ in $\BDV_{\mu}(\Gamma)$, then
by Proposition \ref{Prop33}, $f$ belongs to $\Zh(\Gamma)$ and satisfies
the differential equation
\begin{equation}  \label{tw3e}
\frac{d^2f}{dx^2} + \lambda f(x) \ = \ \lambda C g(x) \ 
\end{equation}
on each $e_{\ell}$.  Necessarily $\lambda > 0$;  write $\lambda = \gamma^2$.   
Suppose that for each $e_{\ell}$ we are able to find a particular solution 
$h_{\gamma,\ell}(x)$ to  $\frac{d^2f}{dx^2} + \lambda f(x) \ = \ \lambda g(x)$.  
Then the general solution to (\ref{tw3e}) has the form
\begin{equation*}
f_{\ell}(x) \ = \ 
A_{\ell} e^{i \gamma x} + B_{\ell} e^{- i \gamma x} +  C h_{\gamma,\ell}(x)  
\end{equation*}
for some constants $A_{\ell}$, $B_{\ell}$.  

Proposition \ref{Prop32} shows that for a collection of solutions 
$\{f_{\ell}(x)\}_{1 \le \ell \le m}$ to be the restriction 
of an eigenfunction $f(x)$,  it is necessary and sufficient that  
$\{A_1, \ldots, A_m, B_1, \ldots, B_m, C\}$ satisfy the 
following system of linear equations: 

A)  For each $p \in X_g$, there are $v(p)-1$ `continuity' conditions,
each of the form
\begin{equation*}
f_i(q_i) \ = \ f_j(q_j), 
\end{equation*}
where $e_i$ is a fixed edge containing $p$, $e_j$ runs over the remaining
edges containing $p$, and $q_i$ (or similarly $q_j$) is the endpoint 
$a_i$ or $b_i$ of $e_i \cong [a_i,b_i]$ corresponding to $p$.

B)  For each $p \in X_g$, there is one `derivative' condition, corresponding to
the fact that $\Delta_{Zh}(f)$ and $-\lambda C \mu$ must have the same point
mass at $p$.  This condition reads 
\begin{equation*}
\sum_i \pm f_i^{\prime}(q_i) \ = \ \lambda C \cdot c_p,
\end{equation*} 
where the sum runs over the edges $e_i$ containing $p$, $q_i$ is the endpoint
of $e_i$ corresponding to $p$, the sign $\pm$ is $+$ if $q_i = a_i$ and
$-$ if $q_i = b_i$, and  $c_p = \mu(\{p\})$. 

C)  There is a global `integral' condition coming from the fact that $f$ belongs
to $\BDV_{\mu}(\Gamma)$:
\begin{equation*}
\int_{\Gamma} f(x) \, d\mu(x) \ = \ 0 \ .
\end{equation*} 

Each condition gives a homogeneous 
linear equation in the $A_{\ell}$, $B_{\ell}$ and $C$, whose
coefficients depend on $\gamma$.   There are $1 + \sum_{p \in X_g} v(p)$
such equations;  since each edge has two endpoints, 
$\sum_{p \in X_g} v(p) = 2m$.  Hence the number of equations is 
the same as the number of variables, $2m+1$.  
For an eigenfunction to exist, 
it is necessary and sufficient that this system of equations should
have a nontrivial solution.  

Let $M(\gamma)$ be the $(2m+1) \times (2m+1)$ matrix of coefficients. 
It follows that $\lambda = \gamma^2$ is an eigenvalue of $\Delta$ in 
$\BDV_{\mu}(\Gamma)$ if and only if $\gamma > 0$ 
satisfies the characteristic equation
\begin{equation} \label{fkj5}
\det(M(\gamma)) \ = \ 0 \ .
\end{equation} 

In all examples known to the authors, the multiplicity of $\gamma$ as a root coincides
with the dimension of the nullspace of $M(\gamma)$, which is the multiplicity
of the corresponding eigenvalue.\footnote{It might be possible to explain this observation using 
the results of \cite{Be1}, which we learned about from the referee.}

\vskip .1 in

\noindent{\bf Example 1}.  Let $\Gamma$ be the segment $[0,1]$, 
and take $\mu = dx$.  
\vskip .05 in

In this case, the eigenvalues are $\lambda_n = n^2 \pi^2$ for $n = 1,2,3, \ldots$
and the eigenfunctions (normalized to have $L^2$-norm $1$) are
\begin{equation*}
f_n(x) \ = \ \sqrt{2} \cos(n \pi x) \ .
\end{equation*}
One can show that 
\begin{equation*}
g_{\mu}(x,y) \ = \ \left\{ \begin{array}{ll}
   \frac{1}{2} x^2 + \frac{1}{2} (1-y)^2 - \frac{1}{6} & \text{if $x < y$} \\ \\
   \frac{1}{2} (1-x)^2 + \frac{1}{2} y^2 - \frac{1}{6} & \text{if $x \ge y$} 
                            \end{array} \right. \ ,
\end{equation*}
\vskip .05 in
\noindent{and the} eigenfunction expansion in Proposition \ref{Prop18} reads 
\begin{equation*}
g_{\mu}(x,y) \ = \ 2 \sum_{n=1}^{\infty} 
          \frac{\cos(n \pi x) \cos(n \pi y)} {n^2 \pi^2} \ .
\end{equation*}
We have
\begin{equation*}
\sum \frac{1}{\lambda_n} \ = \ \sum_{n=1}^{\infty} \frac{1}{\pi^2 n^2}
    \ = \ \int_0^1 g_{\mu}(x,x) \, dx \ = \ \frac{1}{6} \ .
\end{equation*}    
\vskip .1 in                             

\noindent{\bf Example 2}.  Let $\Gamma = [0,1]$, 
and take $\mu = \delta_0(x)$.  
\vskip .05 in

In this case, the eigenvalues are $\lambda_n =  \frac{n^2 \pi^2}{4}$ for $n = 1,3,5,  \ldots$
and the corresponding normalized eigenfunctions are
\begin{equation*}
f_n(x) \ = \ \sqrt{2} \sin(\frac{n \pi}{2} x) \ .
\end{equation*}


One can show that 
\begin{equation*}
g_{\mu}(x,y) \ = \ \min(x,y),
\end{equation*}
and its eigenfunction expansion is 
\begin{equation*}
g_{\mu}(x,y) \ = \ 8 \sum_{\substack{ n \ge 1 \\ \text{$n$ odd} }}
          \frac{\sin(n \pi x) \sin(n \pi y)} {n^2 \pi^2} \ .
\end{equation*}
The trace of $\varphi_\mu$ is
\begin{equation*}
\sum \frac{1}{\lambda_n} \ = \
   \sum_{\substack{ n \ge 1 \\ \text{$n$ odd} }} \frac{4}{\pi^2 n^2} 
    \ = \ \int_0^1 g_{\mu}(x,x) \, dx \ = \ \frac{1}{2} \ .
\end{equation*}    
Note that the upper bound (\ref{FMM1}) for the trace need not hold since 
$\mu \ne \mu_\can$.
\vskip .1 in

\noindent{\bf Example 3}.  

Let $\Gamma = [0,1]$, 
and take $\mu = \frac{1}{2}\delta_0(x)+\frac{1}{2} \delta_1(x)$.  
This is the `canonical measure' for $\Gamma$, as defined in 
Section~\ref{CanonicalTauSection}.

By a calculation similar to that in Examples 1 and 2, we see that
the eigenvalues are $\lambda_n =  n^2 \pi^2$ for $n = 1,3,5,  \ldots$
and the eigenspace corresponding to each eigenvalue is two-dimensional;  
a normalized basis is given by
\begin{equation*}
f_n^{+}(x) \ = \ e^{i n \pi x}, \qquad  f_n^{-}(x) = e^{-in \pi x} \ .
\end{equation*}
One can show that 
\begin{equation*}
g_{\mu}(x,y) \ = \ \frac{1}{4}-\frac{1}{2}|x-y| 
\end{equation*}
and its eigenfunction expansion is 
\begin{equation*}
g_{\mu}(x,y) \ = \  \sum_{\substack{ \text{$n$ odd} }} 
          \frac{e^{i n \pi x} \overline{e^{i n \pi y}}} {n^2 \pi^2} \ .
\end{equation*}
The tau constant $\tau(\Gamma)$ is
\begin{equation*}
\sum \frac{1}{\lambda_n} \ = \
   \sum_{\substack{ n \ge 1 \\ \text{$n$ odd} }} \frac{2}{\pi^2 n^2} 
    \ = \ \int_0^1 g_{\mu}(x,x) \, dx \ = \ \frac{1}{4} \ .
\end{equation*}    
Here $\Gamma$ is a tree, and the upper bound for $\tau(\Gamma)$ 
in (\ref{FMM1}) is achieved.  
\vskip .1 in

\noindent{\bf Example 4}.  
Let $\Gamma = [0,1]$,
and take $\mu = \delta_0(x) + \delta_1(x) - dx$.
\vskip .05 in

This example illustrates phenomena which occur with more
complicated measures.  We have $g(x) = -1$.  It will be convenient to
write the general solution to (\ref{tw3e}) as
\begin{equation*}
f(x) \ = \ A \cos(\gamma x) + B \sin(\gamma x) - C \ .
\end{equation*}
The derivative conditions say
\begin{eqnarray*}
f^{\prime}(0) & = & B \gamma \ = \ \gamma^2 C \ , \\
f^{\prime}(1) & = &
   - A \gamma \sin(\gamma) + B \gamma \cos(\gamma) \ = \ - \gamma^2 C  \ ,
\end{eqnarray*} 
and the integral condition says
\begin{equation*}
A \big(1 + \cos(\gamma) + \frac{\sin(\gamma)}{\gamma}\big)
    + B\big(\sin(\gamma) + \frac{\cos(\gamma)-1}{\gamma}\big) - C \ = \ 0 \ .
\end{equation*}
One has 
\begin{equation*}
\det(M(\gamma))
    \ = \ 2 \gamma^3 (1+ \cos(\gamma)) - 3 \gamma^2 \sin(\gamma) \ .
\end{equation*}

It is easy to show that for each odd integer $2k-1 > 0$, there are two solutions
to $M(\gamma)= 0$ near $(2k-1) \pi$:  
one of these, which we will denote by $\gamma_{2k-1}$, 
is slightly less than $(2k-1) \pi$; the other is $\gamma_{2k} = (2k-1) \pi$.
The first six eigenvalues $\lambda_n = \gamma_n^2$ are
\begin{eqnarray*}
\lambda_1 & \cong & 2.854280792 \ , \\
\lambda_2 & = & \pi^2 \ \cong \ 9.869604401 \ , \\
\lambda_3 & \cong & 82.77313456 \ , \\
\lambda_4 & = & 9 \pi^2 \ \cong \ 88.82643963 \ , \\
\lambda_5 & \cong & 240.7215434 \ , \\
\lambda_6 & = & 25 \pi^2 \ \cong \ 246.7401101 \ .
\end{eqnarray*}
The normalized eigenfunctions corresponding to the eigenvalues
$\lambda_{2k-1} = \gamma^2$ are
\begin{equation*}
f_{2k-1}(x) \ = \
    c_{2k-1} ( \cos(\gamma (x-\frac{1}{2})) -
            \frac{1}{\gamma} \sin(\frac{\gamma}{2}))
\end{equation*}
for an appropriate constant $c_{2k-1}$, and 
the normalized eigenfunctions corresponding to the eigenvalues
$\lambda_{2k} = (2k-1)^2 \pi^2$ are
\begin{equation*}
f_{2k}(x) \ = \ \sqrt{2} \cos((2k-1)\pi x) \ .
\end{equation*}

One has
\begin{equation*}
g_{\mu}(x,y) \ = \ \frac{7}{12} - \frac{1}{2} |x-y|
           -\frac{1}{2} (x-y)^2  - (x-\frac{1}{2}) (y-\frac{1}{2}) \ ,
\end{equation*}
and the trace of $\varphi_\mu$ is
\begin{equation*}
\sum \frac{1}{\lambda_n} \ = \
   \int_0^1 g_{\mu}(x,x) \, dx \ = \ \frac{1}{2} \ .
\end{equation*}    
\vskip .1 in

\noindent{\bf Example 5}.  
Let $\Gamma$ be a circle of length $L$,
and take $\mu = \frac{1}{L} dx$ to be the canonical measure on $\Gamma$.
\vskip .05 in

This is the setting of classical Fourier analysis, and 
our formalism leads to the familiar Fourier expansion
of a periodic function.

In this case, the eigenvalues are $\lambda_n =  4 n^2 \pi^2/L^2$ for $n = 1,2,3, \ldots$\ ;
each eigenspace is two-dimensional, with normalized basis
\begin{equation*}
f_n^{+}(x) \ = \ \frac{1}{\sqrt{L}} e^{2 \pi i n x/L},
\qquad  f_n^{-}(x) \ = \frac{1}{\sqrt{L}} e^{-2 \pi i n x/L} \ .
\end{equation*}
One has
\begin{equation*}
g_{\mu}(x,y) \ = \ \frac{1}{2L} |x-y|^2 - \frac{1}{2} |x-y| + \frac{L}{12}
\end{equation*}
and its eigenfunction expansion is 
\begin{equation*}
g_{\mu}(x,y) \ = \  \sum_{n \ne 0} 
       \frac{e^{2 \pi i n x/L} \overline{e^{2 \pi i n y/L}}} {4 n^2 \pi^2/L}
           \ .
\end{equation*}
We have
\begin{equation*}
\sum \frac{1}{\lambda_n} \ = \
   \sum_{n \ne 0} \frac{L^2}{4 \pi^2 n^2}
    \ = \ \int_0^L g_{\mu}(x,x) \, dx \ = \ \frac{L^2}{12},
\end{equation*}    
and $\tau(\Gamma) = \frac{L}{12}$.

\vskip .1 in
During summer 2003, the authors lead a Research Experience for
Undergraduates (REU) on metrized graphs at the 
University of Georgia.  The student participants were
Maxim Arap, Jake Boggan, Rommel Cortez, Crystal Gordan, Kevin Mills,
Kinsey Rowe, and Phil Zeyliger.  Among other things,
the students wrote a MAPLE package implementing the above algorithm,
and used it to investigate numerical examples, some of which are given
in the following table.  In the table, 
each graph has total length $1$, with edges of equal length.
The graph's name is followed by its tau constant, and then by
the two smallest eigenvalues for the measures $\mu = dx$ and $\mu = \mu_{\can}$. 
The multiplicity of each eigenvalue is given in parentheses.  

\vskip .1 in  

\begin{center} 
\begin{tabular}{| l | r | r | r | r | r | } 
\hline
\text{Name}  
        & $\tau(\Gamma)$ & $\lambda_{1,dx} \quad $ & $\lambda_{2,dx} \quad $ 
        & $\lambda_{1,\can}$ \quad& $\lambda_{2,\can} \quad$
                 \\ \hline
$K_{3,3}$
        & .0442 & 199.86(4) & 799.44(5)  & 105.63(1) & 199.86(4) \\ \hline 
$K_5$   
        & .0460 & 332.51(4) & 986.96(5)  &  47.62(1) & 332.51(4) \\ \hline
\text{Petersen}
        & .0353 &  340.93(5) & 1190.79(4) & 107.14(1) & 340.93(5) \\ \hline        
\text{Tetrahedron}
        & .0521 & 131.42(3) & 355.31(2)  & 102.75(1) & 131.42(3) \\ \hline
\text{Cube}
        & .0396 & 218.20(3) & 525.67(3)  & 106.66(1) & 218.20(3) \\ \hline
\text{Octahedron}
        & .0434 & 355.31(3) & 631.65(2) &  47.73(1)  & 355.31(3) \\ \hline
\text{Dodecahedron} 
        & .0264 & 479.25(3) & 1363.73(5) & 107.78(1) & 479.25(3) \\ \hline
\text{Icosahedron}
        & .0399 & 1103.20(3) & 2826.48(5) & 33.31(1) & 1103.20(3) \\ \hline
\end{tabular}
\end{center}

%
%

\vskip .1 in

Here $K_{3,3}$ denotes the complete bipartite graph on 6 vertices, $K_5$ denotes the complete graph on
5 vertices, and `Petersen' denotes the well-known Petersen graph, which is a 3-regular graph on 10 vertices.
The other graph names are hopefully self-explanatory.  For more data, together with some of the MAPLE routines used to
calculate them, see \verb+www.math.uga.edu/~mbaker/REU/REU.html+.

In all examples known to us, $\lambda_{1,\can}$ has multiplicity $1$.
For the graphs in the table, 
$\lambda_{1,dx}$ coincides with $\lambda_{2,\can}$.  
Factoring the characteristic equation 
for the corresponding eigenvalue problems shows that in fact for these graphs, all 
eigenvalues for $dx$ are eigenvalues for $\mu_{\can}$,  
with multiplicities differing by at most 1 
($\mu_{\can}$ has infinitely many other eigenvalues as well).  
The graphs in the table are highly symmetric: all their 
edges have equal weight under the canonical measure, 
and all their vertices have the same valence.  
For more general graphs, it seems that some 
but not all eigenvalues for $dx$ 
are eigenvalues for $\mu_{\can}$.  It would be interesting to 
understand the reasons for this.  
It would also be very interesting to understand 
the geometric/combinatorial significance of the the eigenvalues of the
canonical measure.

\end{document}